\documentclass[12pt,reqno]{amsart}
\usepackage{amscd}
\usepackage{amssymb}
\usepackage{epsfig}
\usepackage{float}
\usepackage{url}
\usepackage[all]{xy}
\usepackage{graphicx}

\newcommand{\N}{{\mathbb N}}
\newcommand{\Z}{{\mathbb Z}}
\newcommand{\Q}{{\mathbb Q}}
\newcommand{\C}{{\mathbb C}}
\newcommand{\R}{{\mathbb R}}

\newcommand{\CC}{{\mathcal C}}

\newcommand{\PP}{{\mathcal P}}

\newcommand{\www}{\widetilde}

\newcommand{\paa}{\partial}

\newcommand{\Or}{{\rm Or}}

\DeclareMathOperator{\Aut}{Aut}

\DeclareMathOperator{\divv}{div}
\DeclareMathOperator{\discr}{Discr}

\DeclareMathOperator{\id}{id}

\DeclareMathOperator{\lcm}{lcm}
\DeclareMathOperator{\modd}{mod}

\DeclareMathOperator{\norm}{Norm}
\DeclareMathOperator{\ord}{ord}

\DeclareMathOperator{\rk}{rk}

\DeclareMathOperator{\Res}{Res}

\DeclareMathOperator{\supp}{supp}

\DeclareMathOperator{\tr}{tr}

\begin{document}

\theoremstyle{plain}
\newtheorem{lemma}{Lemma}[section]
\newtheorem{definition/lemma}[lemma]{Definition/Lemma}
\newtheorem{theorem}[lemma]{Theorem}
\newtheorem{proposition}[lemma]{Proposition}
\newtheorem{corollary}[lemma]{Corollary}
\newtheorem{conjecture}[lemma]{Conjecture}
\newtheorem{conjectures}[lemma]{Conjectures}

\theoremstyle{definition}
\newtheorem{definition}[lemma]{Definition}
\newtheorem{withouttitle}[lemma]{}
\newtheorem{remark}[lemma]{Remark}
\newtheorem{remarks}[lemma]{Remarks}
\newtheorem{example}[lemma]{Example}
\newtheorem{examples}[lemma]{Examples}
\newtheorem{notation}[lemma]{Notation}
\newtheorem{notations}[lemma]{Notations}

\title
[Integral monodromy of quasihomogeneous singularities]
{The integral monodromy of isolated quasihomogeneous singularities}

\author{Claus Hertling \and Makiko Mase}

\address{Claus Hertling\\Universit\"at Mannheim\\ 
Lehrstuhl f\"ur algebraische Geometrie\\
B6, 26\\
68159 Mannheim, Germany}
\email{hertling@math.uni-mannheim.de}

\address{Makiko Mase\\Universit\"at Mannheim\\ 
Lehrstuhl f\"ur algebraische Geometrie\\
B6, 26\\
68159 Mannheim, Germany}
\email{mmase@mail.uni-mannheim.de}

\thanks{This work was funded by the Deutsche 
Forschungsgemeinschaft (DFG, German Research Foundation) 
-- 242588615}



\date{September 17, 2020}

\begin{abstract}
The integral monodromy on the Milnor lattice of an 
isolated quasihomogeneous singularity is subject of an
almost untouched conjecture of Orlik from 1972.
We prove this conjecture for all iterated Thom-Sebastiani
sums of chain type singularities and cycle type
singularities. The main part of the paper is purely
algebraic. It provides tools for dealing with sums 
and tensor products of $\Z$-lattices with automorphisms 
of finite order and with cyclic generators. 
The calculations are involved. They use fine properties
of unit roots, cyclotomic polynomials, their resultants 
and discriminants. 
\end{abstract}

\maketitle

\tableofcontents

\setcounter{section}{0}

\section{Introduction and main results}\label{c1}
\setcounter{equation}{0}

\noindent
Matrices in $GL(n,\Z)$ arise in algebraic geometry as monodromy matrices.
Usually it is not so difficult to control their eigenvalues and Jordan blocks,
so their conjugacy classes with respect to $GL(n,\C)$, 
but very difficult to understand their conjugacy classes with respect to
conjugation by $GL(n,\Z)$. 

This paper gives some general algebraic tools which deal 
with block diagonal matrices whose blocks are companion 
matrices. And it shows their usefulness and applies them 
in a special situation in algebraic geometry, 
namely in the case of integral monodromy matrices of 
isolated quasihomogeneous singularities. 
We prove an old and beautiful, but almost untouched 
conjecture of Orlik \cite[Conjecture 3.1]{Or72} 
in a number of cases. 
They contain all invertible polynomials, so all iterated
Thom-Sebastiani sums of chain type singularities and 
cycle type singularities.

We start with some notions. Then we formulate Orlik's 
conjecture, and then the results for quasihomogeneous
singularities. The algebraic tools are described
only informally in the introduction.
Instead of conjugacy classes of matrices, we work
now with $\Z$-lattices with endomorphisms.

\begin{definition}\label{t1.1}
Let $H$ be a $\Z$-lattice of rank $n\in\N$, 
and let $h:H\to H$ be an endomorphism. 
The characteristic polynomial of $h$ is called $p_{H,h}$.

\medskip
(a) The pair $(H,h)$ is a {\it companion block}  
if an element $a_0\in H$ exists such that
\begin{eqnarray}\label{1.1}
H=\bigoplus_{j=0}^{n-1}\Z\cdot h^j (a_0).
\end{eqnarray}
Such an element $a_0$ is called a {\it generating element}.

\medskip
(b) A sublattice $H^{(1)}\subset H$ is a {\it companion block} if
it is $h$-invariant and the pair $(H^{(1)},h)$ is a companion block
(here and below, we write $h$ instead of $h|_{H^{(1)}}$).

\medskip
(c) A {\it decomposition of $H$ into companion blocks} is a decomposition
$H=\bigoplus_{i=1}^k H^{(i)}$ such that each $H^{(i)}$ is a companion block.

\medskip
(d) A decomposition as in (c) is a {\it standard} decomposition into 
companion blocks if 
\begin{eqnarray}\label{1.2}
p_{H^{(k)},h} | p_{H^{(k-1)},h} | ... | p_{H^{(2)},h} | p_{H^{(1)},h}.
\end{eqnarray}

\medskip
(e) A companion block $(H,h)$ is called an {\it Orlik block} if $h$ is an
automorphism of finite order. 
Specializing (c) and (d), one obtains the notions 
{\it decomposition into Orlik blocks} and {\it standard decomposition
into Orlik blocks}.
\end{definition}

In \eqref{1.2}, the tuple $(p_{H^{(1)},h},...,p_{H^{(k)},h})$
of characteristic polynomials is unique, see Remark \ref{t2.5} (iv).

A polynomial $f\in\C[x_1,...,x_n]$ is called 
quasihomogeneous if for some weight system 
$(w_1,...,w_n)$ with $w_i\in(0,1)\cap\Q$
each monomial in $f$ has weighted degree 1. 
It is called an isolated quasihomogeneous singularity 
if it is quasihomogeneous and the functions 
$\frac{\paa f}{\paa x_1},...,\frac{\paa f}{\paa x_n}$
vanish simultaneously only at $0\in\C^n$. 
Then the Milnor lattice 
$H_{Mil}:=H_{n-1}^{(red)}(f^{-1}(1),\Z)$
(here $H_{n-1}^{(red)}$ means the reduced homology in the case $n=1$
and the usual homology in the cases $n\geq 2$)
is a $\Z$-lattice of some finite rank $\mu\in\N$, which is called
the {\it Milnor number} \cite{Mi68}. It comes equipped with
an automorphism $h_{Mil}$ of finite order, the {\it monodromy}.

Orlik conjectured the following.

\begin{conjecture}\label{t1.2}
(Orlik's conjecture \cite[Conjecture 3.1]{Or72})
For any isolated quasihomogeneous singularity, the pair
$(H_{Mil},h_{Mil})$ admits a standard decomposition into Orlik blocks.
\end{conjecture}

Conjecture \ref{t1.2} should not be confused with the 
weaker Conjecture 3.2 in \cite{Or72}, which deals with
the homology of the link $f^{-1}(1)\cap S^{2n-1}$ of the
singularity and which would follow from Conjecture 
\ref{t1.2}.
Section \ref{c13} discusses Conjecture 3.2 
in \cite{Or72} and applications of both conjectures.

As an application of our algebraic results,
we prove Conjecture \ref{t1.2} in the following cases.
They surpass all known cases.

\begin{theorem}\label{t1.3}
(a) Conjecture \ref{t1.2} holds for the chain type 
singularities.

\medskip
(b) \cite{HM20-1} 
Conjecture \ref{t1.2} holds for the cycle type 
singularities.

\medskip
(c) If Conjecture \ref{t1.2} holds for a singularity $f$ 
and a singularity $g$, 
then Conjecture \ref{t1.2} holds also for the 
Thom-Sebastiani sum
$f+g=f(x_1,...,x_{n_f})+g(x_{n_f+1},...,x_{n_f+n_g})$.

\medskip
(d) Conjecture \ref{t1.2} holds for all iterated 
Thom-Sebastiani sums of chain type singularities and
cycle type singularities. 
\end{theorem}

Part (a) follows from the combination of Theorem (2.11) in 
the paper \cite{OR77} of Orlik and Randell, 
which we cite as Theorem \ref{t10.1}, 
and our algebraic result Theorem \ref{t6.2}.
Theorem (2.11) in \cite{OR77} gives for a chain type
singularity an automorphism $h:H_{Mil}\to H_{Mil}$ 
such that $h_{Mil}=h^\mu$ and such that $(H_{Mil},h)$
is a single Orlik block. 
The algebraic Theorem \ref{t6.2} starts with a single 
Orlik block $(H,h)$ and a number $\mu\in\N$ and gives 
a sufficient (and probably also necessary) 
condition for $(H,h^\mu)$ 
to admit a standard decomposition into Orlik blocks.

Part (b) follows from \cite{HM20-1}. It builds on the paper
\cite{Co82} of Cooper, who worked on the conjecture,
but made two serious mistakes, see section \ref{c11}. 

Part (c) follows from the basic result
\begin{eqnarray}\label{1.3}
(H_{Mil},h_{Mil})(f+g) &\cong& (H_{Mil},h_{Mil})(f)\otimes (H_{Mil},h_{Mil})(g)
\end{eqnarray}
of Sebastiani and Thom \cite{ST71}, from our algebraic
result Theorem \ref{t9.10} and from Theorem \ref{t12.1}.
Theorem \ref{t9.10} states conditions under which
the tensor product of two standard decompositions
into Orlik blocks admits again a standard decomposition
into Orlik blocks. Theorem \ref{t12.1} verifies
these conditions in the cases of the quasihomogeneous
singularities. 

Part (d) is an immediate consequence of the parts
(a), (b) and (c). The case of curve singularities
(the case $n=1$) is contained in Theorem \ref{t1.3} (d).
For this case Michel and Weber claimed in the introduction
of \cite{MW86} to have a proof of Conjecture \ref{t1.2}. 
In \cite{He92} a few other cases, which are also contained
in Theorem \ref{t1.3} (d), were checked by hand
(using Coxeter-Dynkin diagrams). So, part (d) surpasses
all cases in which Conjecture \ref{t1.2} was known before. 

Theorem \ref{t9.10} builds on other algebraic
results in the sections \ref{c3}, \ref{c5}, \ref{c7},
\ref{c8} and \ref{c9}. 

Theorem \ref{t3.1} allows to improve a certain $\Z$-basis 
of a pair $(H,h)$ to a $\Z$-basis of a standard decomposition 
into companion blocks. Its proof is short and elementary.
Theorem \ref{t3.1} is itself some evidence that standard 
decompositions into companion blocks are natural and 
arise more often than one might expect at first sight.
It is used in the sections \ref{c6} and \ref{c8}.

Theorem \ref{t5.1} compares different decompositions into
companion blocks. 
The proof of Theorem \ref{t5.1} uses results on resultants
which are recalled in section \ref{c4}.
Section \ref{c4} also collects results on cyclotomic
polynomials (and their resultants and discriminants) 
which are partly proved here and partly in \cite{He20}. 

These results are used in the proofs of Theorem
\ref{t6.2} and Theorem \ref{t7.4}. Their
proofs are long, especially the one of Theorem \ref{t7.4}.
It takes the whole section \ref{c8}.
Theorem \ref{t7.4} starts with two single Orlik blocks
$(H^{(1)},h^{(1)})$ and $(H^{(2)},h^{(2)})$ and gives 
a sufficient (and probably also necessary) 
condition for their tensor product
$(H^{(1)}\otimes H^{(2)},h^{(1)}\otimes h^{(2)})$ 
to admit a standard decomposition into Orlik blocks.
The proofs of the Theorems \ref{t6.2} and \ref{t7.4} 
combine Theorem \ref{t3.1} with the calculation of a 
certain determinant. 
The determinant calculations deal a lot with unit roots, 
cyclotomic polynomials, their resultants and their discriminants.

Section \ref{c9} builds on Theorem \ref{t7.4}.
It provides a condition on one Orlik block such that
for any two Orlik blocks which satisfy such conditions,
their tensor product admits a standard decomposition
into Orlik blocks, see Theorem \ref{t9.9}. 
Theorem \ref{t9.10} builds on this. 

In section \ref{c10}, chain type singularities are
introduced, and Theorem \ref{t1.3} (a) is proved.
In section \ref{c11}, cycle type singularities are 
introduced and remarks on \cite{HM20-1} and \cite{Co82}
are made. 
In section \ref{c12}, the conditions in Theorem \ref{t9.10}
are verified in the case of quasihomogeneous singularities.
Section \ref{c13} recalls the second conjecture 3.2
in \cite{Or72}, it tells about applications of 
both conjectures of Orlik, and it formulates an open problem.

\section{Some notations and basic observations}\label{c2}
\setcounter{equation}{0}

\noindent
First, we fix some basic general notations.

\begin{notations}\label{t2.1}
$\N:=\{1,2,3,...\}$, $\N_0:=\N\cup\{0\}$. 
The subset of $\N$ of all prime numbers is denoted by 
$\PP\subset\N$. 

For any $m\in\N$ and any prime number $p$ denote
\begin{eqnarray*}
v_p(m)&:=& \max(k\in\N_0\, |\, p^k\textup{ divides }m).
\end{eqnarray*}
Thus $m=\prod_{p\textup{ prime number}}p^{v_p(m)}.$

For any subset $I\subset\R$ denote 
(especially for an interval $[r_1,r_2]\subset\R$)
\begin{eqnarray*}
\Z_I&:=& \Z\cap I.
\end{eqnarray*}

For a case discussion, whether or not a certain condition
$(Cond)$ holds, the Kronecker delta is generalized as
follows,
\begin{eqnarray*}
\delta_{(Cond)}&:=& \left\{\begin{array}{ll}
1 &\textup{ if }(Cond)\textup{ holds},\\
0 &\textup{ if }(Cond)\textup{ does not hold.}
\end{array}\right. 
\end{eqnarray*}
\end{notations}
Beyond this section, some general notations are given in the
notations and definitions \ref{t4.1}, \ref{t4.4}, 
\ref{t6.1}, \ref{t7.1}, \ref{t9.1}, \ref{t9.3} 
and \ref{t9.4}. 
The next notations fix our way to deal with matrices 
and bases.

\begin{notation}\label{t2.2}
Let $R$ be a principal ideal domain, 
and let $V$ a free $R$-module of rank $n\in\N$.
Let ${\bf a}=(a_1,...,a_n)$ be an ordered $R$-basis of $V$,
and let ${\bf b}=(b_1,...,b_k)\in R^k$ for some $k\in\N$.
Let $f:V\to V$ be an $R$-linear endomorphism. 
Then $M({\bf a},f,{\bf b})$ denotes the matrix which expresses
the elements $f(b_1),...,f(b_k)$ as linear combinations of $a_1,...,a_n$,
namely $M({\bf a},f,{\bf b})=(r_{ij})_{i=1,...,n;j=1,...,k}
\in M_{n\times k}(R)$ with
\begin{eqnarray}\label{2.1}
f(b_j)&=& \sum_{i=1}^n r_{ij}\cdot a_i.
\end{eqnarray}
We write this simultaneously for all $j$ as follows,
\begin{eqnarray}\label{2.2}
f({\bf b}) &=& {\bf a}\cdot M({\bf a},f,{\bf b}).
\end{eqnarray}
If ${\bf b}$ is also an $R$-basis of $V$, and if ${\bf c}\in R^l$
is a tuple and $g:V\to V$ is a second endomorphism, then
the calculation
\begin{eqnarray*}
f(g({\bf c})) &=& f({\bf b}\cdot M({\bf b},g,{\bf c})) 
 = f({\bf b})\cdot M({\bf b},g,{\bf c})\\
&=& {\bf a}\cdot M({\bf a},f,{\bf b})\cdot M({\bf b},g,{\bf c})
\end{eqnarray*}
shows
\begin{eqnarray}\label{2.3}
M({\bf a},f\circ g,{\bf c})&=& M({\bf a},f,{\bf b})\cdot
M({\bf b},g,{\bf c}).
\end{eqnarray}
Especially, in the case $f=\id$, we write $M({\bf a},{\bf b}):=
M({\bf a},\id,{\bf b})$. If ${\bf b}$ is also an $R$-basis of $V$, 
this is in $GL_n(R)$. 
\end{notation}

\begin{remark}\label{t2.3}
Let $R$ be a principal ideal domain. Because of the Notation \ref{t2.2},
the following is clear.
The isomorphism class of a pair $(V,h)$ with $V$ a free $R$-module
of rank $n\in\N$ and $h:V\to V$ an endomorphism is
equivalent to the conjugacy class 
$$\{M({\bf a},h,{\bf a})\, |\, {\bf a}\textup{ an }R\textup{-basis of }V\}$$
of matrices in $M_{n\times n}(R)$ with respect to 
$GL_n(R)$.
\end{remark}

\begin{notation}\label{t2.4}
If $H$ is a $\Z$-lattice of some finite rank and $R$ is a principal
ideal domain which contains $\Z$, then $H_R:=H\otimes_\Z R$ is an $R$-lattice
of the same rank. Then a $\Z$-linear endomorphism $h$ of $H$ extends to
an $R$-linear endomorphism of $H_R$.
In the case of $R=\C$, the generalized eigenspace with eigenvalue 
$\lambda\in\C$ of $h:H_\C\to H_\C$ is denoted by $H_\lambda$,
so that $H_\C=\bigoplus_{\lambda\textup{ eigenvalue}}H_\lambda$.
\end{notation}

\begin{remarks}\label{t2.5}
(i) Let $R$ be a principal ideal domain, and let $V$ be a finitely generated
$R$-module. Then it is a basic theorem on such $R$-modules that $V$ is
isomorphic to a direct sum of quotients
\begin{eqnarray*}
R^{k_1}\oplus \bigoplus_{j=1}^{k} \frac{R}{p_jR}\quad\textup{with}\quad
p_1,...,p_k\in R-(R^*\cup\{0\}), \quad p_k|p_{k-1}|...|p_2|p_1,
\end{eqnarray*}
($R^*:=\{\textup{units in }R\})$.
The numbers $k_1$ and $k$ are unique, 
and the elements $p_1,...,p_k$ are unique up to multiplication by units, 
and they are called {\it elementary divisors}.
$R^{k_1}$ is the {\it free part}, and 
$\bigoplus_{j=1}^kR/p_jR$ is the {\it torsion part} 
of the sum.

\medskip
(ii) Let $(H,h)$ be a $\Z$-lattice of rank $n\in\N$ and $h:H\to H$
an endomorphism. Then $H$ is a $\Z[t]$-module, where $t$ acts as $h$ on $H$.
The ring $\Z[t]$ is not a principal ideal domain, but $\Q[t]$ is.
$H_\Q$ is a $\Q[t]$-module. Part (i) applies. $H_\Q$ is a torsion module
of $\Q[t]$, it is isomorphic to $\bigoplus_{j=1}^k\Q[t]/p_j\Q[t]$
for unique unitary polynomials $p_1,...,p_k\in\Q[t]$ of degrees $\geq 1$, 
which satisfy $p_k|p_{k-1}|...|p_2|p_1$. They are the elementary divisors.
In fact, 
\begin{eqnarray}\label{2.4}
p_{H,h}&=& p_1\cdot ...\cdot p_k,
\end{eqnarray}
and as this is unitary and in $\Z[t]$, all $p_j$ are in $\Z[t]$.

\medskip
(iii) $(H,h)$ as in (ii) admits a standard decomposition into 
companion blocks if and only if a decomposition of $H_\Q$ as in (ii)
lifts from $\Q[t]$ to $\Z[t]$, so that $H$ is isomorphic to
$\bigoplus_{j=1}^k\Z[t]/p_j\Z[t]$ as a $\Z[t]$-module.
Orlik \cite{Or72} gave this fact as a motivation for Conjecture \ref{t1.2}.

\medskip
(iv) The tuple of characteristic polynomials $p_{H^{(j)},h}$ in \eqref{1.2} is unique,
because over $\Q$ they become the elementary divisors of the $\Q[t]$-module
$H_\Q$. They can also be understood in terms of the Jordan block structure
of $h$ on $H_\C$. If $\lambda_1,...,\lambda_m\in\C$ are the different
eigenvalues and if for $\lambda_i$ there are Jordan blocks of sizes
$b_{i,1},...,b_{i,m(i)}$ with $b_{i,1}\geq b_{i,2}\geq ...\geq b_{i,m(i)}$
$(m(i)\geq 1)$ then
\begin{eqnarray}
k&=& \max(m(i)\, |\, i=1,...,m),\nonumber\\
p_{H^{(j)},h}&=& \prod_{i:\, m(i)\geq j} (t-\lambda_i)^{b_{i,j}}.
\label{2.5}
\end{eqnarray}

\medskip
(v) Choose a unitary polynomial $p\in \C[t]$ of some degree $n\geq 1$.
Define $H_\C:=\C[t]/p\C[t]$ and $h:H_\C\to H_\C$ as multiplication by $t$.
Consider the $\C$-basis ${\bf a}=(1,[t],...,[t^{n-1}])$ of $H_\C$.
The matrix $M({\bf a},h,{\bf a})$ is the {\it companion matrix
with characteristic polynomial $p$}, 
\begin{eqnarray}\label{2.6}
M({\bf a},h,{\bf a})&=& \begin{pmatrix}
0 &        & & -p_0 \\
1 & \ddots & & \vdots \\
  & \ddots & 0 & -p_{n-2} \\
  &        & 1 & -p_{n-1}  
\end{pmatrix} \\
\textup{where}\quad p(t)&=& t^n+p_{n-1}t^{n-1}+...+p_1t+p_0.\nonumber
\end{eqnarray}
It is a regular matrix, that means, it has for each eigenvalue only
one Jordan block. Of course, the characteristic polynomial of $h$ on $H_\C$ 
is $p_{H_\C,h}=p$, and it is also the minimal polynomial of $h$ on $H_\C$.

\medskip
(vi) If in (v) $p\in \Z[t]$ and $H:=\Z[t]/p\Z[t]$, then $(H,h)$ is a
companion block, and $p_{H,h}=p$.

\medskip
(vii) Vice versa, if $(H,h)$ is a companion block of rank $n\in\N$ 
and $a_0\in H$ is a generating element of it as in \eqref{1.1}, 
then the matrix $M({\bf a},h,{\bf a})$ of $h$ with respect to the 
$\Z$-basis ${\bf a}=(a_0,h(a_0),...,h^{n-1}(a_0))$ is the
companion matrix with characteristic polynomial $p_{H,h}$,
and the map $H\to \Z[t]/p\Z[t],\ h^j(a_0)\mapsto [t^j]$, induces
an isomorphism of $\Z$-lattices and $\Z[t]$-modules. 
\end{remarks}

\begin{definition}\label{t2.6}
(a) 
Let $p\in\Z[t]$ be a unitary polynomial of some degree $n\geq 0$.
We denote $H^{[p]}:=\Z[t]/p\Z[t]$, and $h_{[p]}:H^{[p]}\to H^{[p]}$ is the
multiplication by $t$. If $n=0$ then $H^{[p]}=0$. 
If $n\geq 1$ then $(H^{[p]},h_{[p]})$ is a companion block, because of (vi) above,  
Up to isomorphism, it is the unique companion block
with characteristic polynomial $p$.

\medskip
(b) Let $M\subset\N$ be a finite subset. It defines
a unitary polynomial $p_M:=\prod_{m\in M}\Phi_m$
(where $\Phi_m$ is the cyclotomic polynomial, see
the Notation \ref{t4.4} (iii)). Then
$(H^{[p]},h_{[p]})$ is an Orlik block. 
It is also denoted $\Or(M):=(H^{[p]},h_{[p]})$. 
$M$ is the set of orders of eigenvalues of its monodromy.
The Orlik block and $M$ determine one another. 
\end{definition}

\begin{remarks}\label{t2.7}
(i) Let $H$ be a $\Z$-lattice of rank $n\in\N$.
A sublattice $H^{(1)}\subset H$ is {\it primitive}
if the quotient $H/H^{(1)}$ has no torsion, or, equivalently,
if $H^{(1)}=H\cap H^{(1)}_\Q$, where $H^{(1)}_\Q\subset H_\Q$
and $H\subset H_\Q$, or, equivalently, 
if a sublattice $H^{(2)}\subset H$ with $H=H^{[1]}\oplus H^{[2]}$ exists.
For any sublattice $H^{(3)}\subset H$,
a unique primitive sublattice $H^{(4)}\subset H$ with
$H^{(4)}_\Q=H^{(3)}_\Q$ exists, namely $H^{(4)}=H\cap H^{(3)}_\Q\supset H^{(3)}$.

\medskip
(ii) Let $(H,h)$ be a $\Z$-lattice and $h:H\to H$ an endomorphism
with characteristic polynomial $p=p_1\cdot p_2$ with $p_1,p_2\in\Z[t]$
unitary of degrees $\geq 1$. Then $p_1(h)p_2(h)=0$ and
\begin{eqnarray}\label{2.7}
p_2(h)(H)\subset \ker(p_1(h):H\to H).
\end{eqnarray}
The second sublattice is a kernel, so it is a primitive sublattice.
If $h$ is regular (only one Jordan block for each eigenvalue) then $p$
is also the minimal polynomial, and then the two sublattices
$p_2(h)(H)$ and $\ker(p_1(h):H\to H)$ have the same rank, so the first
has a finite index in the second. Often this index is $>1$.
But by Lemma \ref{t2.8} (b), equality holds if $(H,h)$ is a companion block. 
\end{remarks}

Part (a) of Lemma \ref{t2.8} was stated and proved in \cite{He20}.

\begin{lemma}\label{t2.8}
Let $p_1,p_2\in \Z[t]$ be unitary polynomials of degrees $\geq 1$.

\medskip
(a) $H^{[p_1p_2]}$ contains a unique primitive sublattice which is
$h_{[p_1p_2]}$-invariant and such that the characteristic polynomial
of $h_{[p_1p_2]}$ on it is $p_1$. It is $(p_2)/(p_1p_2)\subset H^{[p_1p_2]}$,
and $(p_2)/(p_1p_2)\cong H^{[p_1]}$.

\medskip
(b) Let $(H,h)$ be a companion block with characteristic polynomial
$p_1p_2$. Then 
\begin{eqnarray}\label{2.8}
p_2(h)(H)=\ker(p_1(h):H\to H).
\end{eqnarray}
Especially, $p_2(h)(H)$ is a primitive sublattice of $H$.
\end{lemma}

{\bf Proof:} (a) \cite[Lemma 6.1]{He20} (not difficult).

(b) We can choose $(H,h)=(H^{[p_1p_2]},h_{[p_1p_2]})$. Then
\begin{eqnarray*}
p_2(h)(H)&=&(p_2)/(p_1p_2) \\
&=&\ker((\textup{multiplication by }p_1):H^{[p_1p_2]}
\to H^{[p_1p_2]}) \\
&=& \ker(p_1(h):H\to H). \hspace*{4cm}\Box 
\end{eqnarray*}

\section{A frame for constructing a standard 
decomposition into companion blocks}\label{c3}
\setcounter{equation}{0}

\noindent
Let $(H,h)$ be a $\Z$-lattice of rank $n\in\N$ with an 
endomorphism $h:H\to H$. By Remark \ref{t2.5}, the 
$\Q$-vector space $H_\Q$ as a $\Q[t]$-module is isomorphic to
$\bigoplus_{j=1}^k \Q[t]/p_j\Q[t]=\bigoplus_{j=1}^k H^{[p_j]}_\Q$
(see Definition \ref{t2.6} and Notation \ref{t2.4} for $H^{[p_j]}_\Q$), 
where $p_1,...,p_k\in\Z[t]$ are the elementary divisors.
They are unitary of degrees $\geq 1$ and satisfy $p_k|p_{k-1}|...|p_2|p_1$
and $p_{H,h}=p_1\cdot ...\cdot p_k$. 
If $(H,h)$ admits a standard decomposition into companion blocks, 
that is isomorphic to $\bigoplus_{j=1}^k H^{[p_j]}$.

Theorem \ref{t3.1} gives a frame for constructing a standard decomposition
into companion blocks. It builds on Lemma \ref{t2.8} (b).

\begin{theorem}\label{t3.1}
Let $(H,h)$ and $p_1,...,p_k$ be as above. 
Consider $k$ elements $a^{[j]}_0\in H$ for $j\in\{1,...,k\}$,
and consider the elements $a^{[j]}_i:=h^i(a^{[j]}_0)$ for $i\in\N_0$.
If the tuple of elements
\begin{eqnarray}
(a^{[1]}_0,a^{[1]}_1,...,a^{[1]}_{\deg p_1-1}, \nonumber \\
a^{[2]}_0,a^{[2]}_1,...,a^{[2]}_{\deg p_2-1},\nonumber \\
...,a^{[k]}_0,a^{[k]}_1,...,a^{[k]}_{\deg p_k-1})\nonumber \\
=(a^{[j]}_i\, |\, j=1,...,k,\ i=0,...,\deg p_j-1)\label{3.1}
\end{eqnarray}
is a $\Z$-basis of $H$, then $(H,h)$ admits a standard decomposition
$H=\bigoplus_{j=1}^k B^{[j]}$
into companion blocks $B^{[j]}$. The companion blocks can be chosen as follows.
The first block $B^{[1]}$ with characteristic polynomial $p_1$ 
is generated by $a^{[1]}_0$.
The $j$-th block $B^{[j]}$ for $j\in\{2,...,k\}$ with characteristic polynomial $p_j$ 
is generated by $a^{[j]}_0+b^{[j]}$ where $b^{[j]}$ is a certain element of the 
sum of the first $j-1$ blocks $B^{[1]},...,B^{[j-1]}$.
\end{theorem}

{\bf Proof:} Suppose that the tuple in \eqref{3.1} is a $\Z$-basis of $H$.
Then its first $\deg p_1$ elements 
$a^{[1]}_0,a^{[1]}_1,...,a^{[1]}_{\deg p_1-1}$ generate a primitive
sublattice $B^{[1]}$ of rank $\deg p_1$. It is $h$-invariant as 
$p_1$ is the minimal polynomial of $h$ on $H$ and $a^{[1]}_i=h^i(a^{[1]}_0)$.
So it is a companion block with characteristic polynomial $p_1$ and 
generator $a^{[1]}_0$. We define $b^{[1]}:=0$.

Now we proceed by induction on $j$ and suppose that for some $j\in\{2,...,k\}$
the following holds. For any $l\in\{1,...,j-1\}$ an element
$b^{[l]}\in H$ has been constructed such that $a^{[l]}_0+b^{[l]}$ is 
a generator of a companion block 
\begin{eqnarray}\label{3.2}
B^{[l]}&=&\bigoplus_{i=0}^{\deg p_l-1}\Z\cdot h^i(a_0^{[l]}+b^{[l]})
\end{eqnarray}
with characteristic polynomial $p_l$
and such that $b^{[l]}$ is in the sublattice generated by
the $l-1$ companion blocks $B^{[1]},...,B^{[l-1]}$ constructed before.

Then the sublattice generated by the first $j-1$ companion blocks is
\begin{eqnarray}\label{3.3}
\bigoplus_{l=1}^{j-1} B^{[l]}
&=& 
\bigoplus_{l=1}^{j-1}\bigoplus_{i=0}^{\deg p_l-1}\Z \cdot a^{[l]}_i .
\end{eqnarray}
It is indeed a direct sum, and it is a primitive sublattice of $H$,
both because the tuple in \eqref{3.1} is a $\Z$-basis of $H$.

We want to find an element 
$b^{[j]}\in \bigoplus_{l=1}^{j-1}B^{[l]}$ 
such that $a^{[j]}_0+b^{[j]}$ generates a companion block 
with characteristic polynomial $p_j$.
We claim that it is sufficient to find an element $b^{[j]}$ with
\begin{eqnarray}\label{3.4}
p_j(h)(a^{[j]}_0+b^{[j]})&=& 0
\end{eqnarray}
That the tuple in \eqref{3.1} is a $\Z$-basis shows that we have
a direct sum
\begin{eqnarray}\label{3.5}
\left(\bigoplus_{l=1}^{j-1}B^{[l]}\right)
\oplus\left( \bigoplus_{i=0}^{\deg p_j-1}\Z\cdot h^i(a^{[j]}_0+b^{[j]})\right) .
\end{eqnarray}
\eqref{3.4} and \eqref{3.5} show that the second summand in \eqref{3.5} 
is a companion block $B^{[j]}$ with characteristic polynomial $p_j$ 
and that it extends the sum of the earlier constructed 
companion blocks to a bigger primitive sublattice.

So it remains to find $b^{[j]}\in \bigoplus_{l=1}^{j-1}B^{[l]}$ with \eqref{3.4}.

First we consider $H_\Q$ and $p_j(h)(H_\Q)\subset H_\Q$. 
By Remark \ref{t2.5} (ii), $H_\Q$ has a decomposition
$H_\Q=\bigoplus_{l=1}^k \www B^{[l]}_\Q$ into $h$-invariant blocks $\www B^{[l]}_\Q$
with characteristic and minimal polynomials $p_l$.
For $l\geq j$ we have $p_l|p_j$ and thus $p_j(h)(\www B^{[l]}_\Q)=0$.
Therefore
\begin{eqnarray}\label{3.6}
p_j(h)(H_\Q)&=& \bigoplus_{l=1}^{j-1} p_j(h)(\www B^{[l]}_\Q)\\
&=& \bigoplus_{l=1}^{j-1}\ker(\frac{p_l}{p_j}(h):\www B^{[l]}_\Q
\to \www B^{[l]}_\Q),\nonumber
\end{eqnarray}
and this has dimension $\sum_{l=1}^{j-1}(\deg p_l-\deg p_j)$. 
The subspace
\begin{eqnarray}\label{3.7}
\bigoplus_{l=1}^{j-1} p_j(h)(B^{[l]}_\Q)
= \bigoplus_{l=1}^{j-1}\ker(\frac{p_l}{p_j}(h):B^{[l]}_\Q
\to B^{[l]}_\Q) 
\end{eqnarray}
has the same dimension, thus it is equal to $p_j(h)(H_\Q)$. 
By Lemma \ref{t2.8} (b) $p_j(h)(B^{[l]})$ is for any $l\in\{1,...,j-1\}$
a primitive sublattice of $B^{[l]}$. Therefore the sum
$\bigoplus_{l=1}^{j-1} p_j(h)(B^{[l]})$ is a primitive sublattice of
$\bigoplus_{l=1}^{j-1} B^{[l]}$. As this is a primitive sublattice of $H$,
the sum $\bigoplus_{l=1}^{j-1} p_j(h)(B^{[l]})$ is a primitive sublattice of $H$.
As its extension to $H_\Q$ is equal to $p_j(h)(H_\Q)$, it is itself
equal to $p_j(h)(H_\Q)\cap H$. We obtain
\begin{eqnarray}
p_j(h)(H)&\supset& p_j(h)(\bigoplus_{l=1}^{j-1} B^{[l]})
= \bigoplus_{l=1}^{j-1} p_j(h)(B^{[l]}) \nonumber \\
&= & p_j(h)(H_\Q)\cap H\supset p_j(h)(H), \nonumber\\
\textup{thus}\quad 
p_j(h)(H) &=& p_j(h)(\bigoplus_{l=1}^{j-1} B^{[l]}).\label{3.8}
\end{eqnarray}
Therefore $b^{[j]}\in \bigoplus_{l=1}^{j-1} B^{[l]}$ with \eqref{3.4} exists.
\hfill$\Box$

\begin{remark}\label{t3.2}
(i) Theorem \ref{t3.1} and its proof are evidence that 
{\it standard} decompositions into companion blocks are 
natural and arise more often than one might
expect at first sight.

\medskip
(ii) For any $j\in \{1,...,k\}$, the $\Z$-lattice
$V_j:=\langle\textup{the first }j\textup{ lines in 
\eqref{3.1}}\rangle_\Z$ is an $h$-invariant 
primitive sublattice of $H$. The proof of Theorem \ref{t3.1}
shows this and provides a $\Z$-sublattice $B_j$
with $V_j=B_j\oplus V_{j-1}=\bigoplus_{l=1}^jB_l$. 
\end{remark}

\section{Resultants and discriminants of cyclotomic polynomials}\label{c4}
\setcounter{equation}{0}

\noindent
This section reviews resultants and discriminants 
in general and in the case of cyclotomic polynomials. 
This is needed in the determinant calculations 
in the proofs of Theorem \ref{t6.2} and Theorem \ref{t7.4}. 
Essentially everything in this section is well-known.
The review here is close to section 2 in \cite{He20}.

The resultant of two polynomials and the discriminant of 
one polynomial are very classical objects. One reference 
for the following definition is \cite[\S 34]{vW71}.

\begin{definition}\label{t4.1}
(a) The {\it resultant} of two polynomials 
$f=\sum_{i=0}^m f_it^i\in\C[t]-\{0\}$
and $g=\sum_{j=0}^ng_jt^j\in\C[t]-\{0\}$ 
of degrees $\deg f=m,\deg g=n$
with $m+n\geq 1$ is $\Res(f,g):=\det A(f,g)\in\C$ 
where $A(f,g)\in M_{(m+n)\times (m+n)}(\C)$ is
the matrix
\begin{eqnarray}\label{4.1}
A(f,g)=\begin{pmatrix}
f_0 & 0 & \dots & 0 & g_0 & 0 & \dots & 0 \\
f_1 & f_0 & \ddots & \vdots & g_1 & g_0 & \ddots & \vdots \\
\vdots & f_1 & \ddots & 0 & \vdots & g_1 & \ddots & \vdots \\
\vdots & \ddots & \ddots & f_0 & \vdots & \ddots & \ddots & g_0 \\
f_m & \ddots & \ddots & f_1 & g_n & \ddots & \ddots & g_1 \\
0 & \ddots & \ddots & \vdots & 0 & \ddots & \ddots & \vdots \\
\vdots & \ddots & \ddots & \vdots & \vdots & \ddots & \ddots & \vdots \\
0 & \dots & 0 & f_m & 0 & \dots & 0 & g_n 
\end{pmatrix}
\end{eqnarray}
whose first $n$ columns contain the coefficients of $f$
and whose last $m$ columns contain the coefficients of $g$.
In other words, it is the matrix with
\begin{eqnarray}\label{4.2}
(f,tf,...,t^{n-1}f,g,tg,...,t^{m-1}g) = (1,t,...,t^{m+n-1})\cdot A(f,g).
\end{eqnarray}
In the case $m+n=0$ one defines $\Res(f,g):=1$.

\medskip
(b) The {\it discriminant} of a polynomial 
$f=\sum_{i=0}^m f_it^i\in\C[t]$ of degree $\deg f=m\geq 1$
is $\discr(f):=\Res(f,\frac{d}{dt}f)\in\C$. 
\end{definition}

The basic properties of the resultant and the discriminant are well known.

\begin{proposition}\label{t4.2}
(a) Let $f$ and $g\in\C[t]$ be as in definition \ref{t4.1} (a).
Let $a_1,...,a_m\in\C$ and $b_1,...,b_n\in\C$ be the zeros of $f$ and $g$, so
$$f=f_0\cdot\prod_{i=1}^m(t-a_i),\quad g=g_0\cdot\prod_{j=1}^n(t-b_j).$$
Then
\begin{eqnarray}\label{4.3}
\Res(f,g)&=&f_0^n g_0^m\cdot \prod_{i=1}^m\prod_{j=1}^n(a_i-b_j)\\
&=& f_0^n\cdot\prod_{i=1}^m g(a_i) 
= (-1)^{m\cdot n}g_0^m\cdot\prod_{j=1}^n f(b_j) \label{4.4}\\
&=&(-1)^{m\cdot n}\Res(g,f),\label{4.5}\\
\Res(f,g)\neq 0&\iff& \gcd(f,g)_{\C[t]}=1.\label{4.6}
\end{eqnarray}

(b) If $f,g,h\in\C[t]-\{0\}$ then
\begin{eqnarray}\label{4.7}
\Res(f,gh)=\Res(f,g)\cdot \Res(f,h).
\end{eqnarray}
If $f^{(1)},...,f^{(r)},g^{(1)},...,g^{(s)}\in\C[t]-\{0\}$ then
\begin{eqnarray}\label{4.8}
\Res(\prod_{i=1}^r f^{(i)},\prod_{j=1}^s g^{(j)})
=\prod_{i=1}^r\prod_{j=1}^s \Res(f^{(i)},g^{(j)}).
\end{eqnarray}

(c) If $f=f_0\cdot \prod_{i=1}^m(t-a_i)\in\C[t]$ has degree $m\geq 1$ 
and zeros $a_1,...,a_m\in\C$ then
\begin{eqnarray}\label{4.9}
\discr(f)&=& f_0^{2m-1}\cdot\prod_{(i,j)\in\{1,...,m\}^2
\textup{ with }i\neq j}(a_i-a_j).
\end{eqnarray}
\end{proposition}

\eqref{4.3} is proved for example in \cite[\S 35]{vW71}. \eqref{4.4}, \eqref{4.5}, 
\eqref{4.6} and \eqref{4.7} follow from \eqref{4.3}. \eqref{4.8}
follows from \eqref{4.7} and \eqref{4.5}. And \eqref{4.9} follows from
\eqref{4.4} and $\frac{d}{dt}f=\sum_{i=1}^m f_0\cdot \prod_{j\neq i}(t-a_j)$.

We are mainly interested in $\Res(f,g)$ where $f$ and $g$ are unitary
polynomials. 
We denote for $k\in\Z_{\geq -1}$
\begin{eqnarray}\label{4.10}
\C_{k}[t]&:=&\{h\in\C[t]\, |\, \deg h\leq k\},\\
\Z_{k}[t]&:=&\C_{k}[t]\cap \Z[t]\nonumber
\end{eqnarray}
(so that $\C_{-1}[t]=\Z_{-1}[t]=\{0\}$).
The following lemma is proved for example in \cite{He20}.

\begin{lemma}\label{t4.3} \cite[Lemma 2.3]{He20}
Let $f,g\in\Z[t]$ be unitary polynomials of degrees
$m=\deg f,n=\deg g$. They generate
an ideal $(f,g)\subset \Z[t]$ (here $\Z[t]$ is also considered
as an ideal).

(a) 
\begin{eqnarray}\label{4.11}
\Z_{ n-1}[t]\cdot f + \Z_{ m-1}[t]\cdot g 
=(f,g)\cap \Z_{ m+n-1}[t].
\end{eqnarray}

(b) The $\Z$-lattice in \eqref{4.11} has rank $m+n$ if and
only if $\Res(f,g)\neq 0$, and then
\begin{eqnarray}\label{4.12}
|\Res(f,g)| = | \frac{\Z_{ m+n-1}[t]}{(f,g)\cap \Z_{ m+n-1}[t]}|
\in\Z_{>0}.
\end{eqnarray}

(c) 
\begin{eqnarray}\label{4.13}
|\Res(f,g)|=1\iff (f,g)=\Z[t].
\end{eqnarray}
\end{lemma}

Now we turn to unit roots and cyclotomic polynomials.

\begin{notations}\label{t4.4}
(i) In the proof of Theorem \ref{t4.5} (and only there), 
we will use the notation $[a]_m$ for the class of $a\in\Z$ in $\Z/m\Z$.

(ii) The order $\ord(\lambda)\in\N$ of a unit root 
$\lambda\in S^1\subset\C$ is the minimal $k\in\N$ with $\lambda^k=1$.
In the rest of this section, $\lambda$ denotes always a unit root.
In the rest of this paper, $e(z)$ for $z\in\C$ denotes 
$e^{2\pi i z}\in\C$, so for example $e(r)$ for $r\in\Q$ is a unit root.

(iii) For $m\in\N$, the cyclotomic polynomial $\Phi_m$ is the polynomial
\begin{eqnarray}\label{4.14}
\Phi_m(t) := \prod_{\lambda:\ \ord(\lambda)=m}(t-\lambda),
\end{eqnarray}
whose zeros are the $m$-th primitive unit roots.
It is a unitary and irreducible polynomial in $\Z[t]$ of degree 
$\deg\Phi_m=\varphi(m)\in\N$, where $\varphi:\N\to
\N$ is the Euler phi-function
(see e.g. \cite[Ch 1,2]{Wa82}). Except for the irreducibility, this follows
easily inductively from the formula
\begin{eqnarray}\label{4.15}
t^m-1=\prod_{k|m}\Phi_k.
\end{eqnarray}
Using this formula, one can compute the $\Phi_k$ inductively.

(iv) Recall (see e.g. \cite[Ch 1,2]{Wa82}) that $\Z[e(\frac{1}{m})]$
is the ring of the algebraic integers within $\Q[e(\frac{1}{m})]$
and that 
\begin{eqnarray}\label{4.16}
\Z[e(\frac{1}{m})]\cap S^1 = \{\pm e(\frac{k}{m})\, |\, k\in\Z\}.
\end{eqnarray}
We will also use the norm
\begin{eqnarray}\label{4.17}
\norm_m:\Z[e(\frac{1}{m})]\to\Z,\quad g(e(\frac{1}{m}))\mapsto
\prod_{\lambda:\ \ord(\lambda)=m}g(\lambda).
\end{eqnarray}
An element of $\Z[e(\frac{1}{m})]$ has norm in $\{\pm 1\}$ if and 
only if it is a unit in $\Z[e(\frac{1}{m})]$.
This and the calculation 
\begin{eqnarray}\label{4.18}
\norm_{\ord(\lambda)}(1-\lambda)=\prod_{\kappa:\, \ord(\kappa)=\ord(\lambda)}
(1-\kappa)=\Phi_{\ord(\lambda)}(1)
\end{eqnarray}
and Theorem \ref{t4.5} (a) imply Theorem \ref{t4.5} (b).
\end{notations}

The following theorem collects relevant facts on unit roots,
and it gives formulas for the resultants and the discriminants
of cyclotomic polynomials.
The parts (a), (b) and (d) and a part of part (c) 
are proved in \cite[Theorem 3.1]{He20}.
Part (d) gives the resultants of the cyclotomic polynomials.  
It is the main result of \cite{Ap70}. The proof of it in \cite{He20}
is shorter than that in \cite{Ap70}.
We do not know a reference for part (e) and the rest of 
part (c), although they are
certainly known. Therefore we provide proofs for them.
Part (e) gives the discriminants of the cyclotomic polynomials.

\begin{theorem}\label{t4.5}
(a) $\Phi_m(1)=1$ if $m\geq 2$ and $m$ is not a power of a prime number.
$\Phi_{p^k}(1)=p$ if $p$ is a prime number and $k\in\N$.

(b) $1-\lambda$ is a unit in $\Z[\lambda]$ if and only if $\ord(\lambda)$
is not a power of a prime number and not equal to $1$.

(c) Fix $m,n\in\Z_{\geq 2}$, $k\in\N$, a prime number $p$, and 
denote
\begin{eqnarray}\label{4.19}
\Lambda(m,n,p,k)&:=& \varphi(p^k)^{-1}\cdot
|\{(a,b)\in (\Z/m\Z)^*\times (\Z/n\Z)^* \ :\hspace*{1cm}\\
&&\hspace*{2cm}\ord(e(\frac{a}{m}-\frac{b}{n}))=p^k\}|\in\N_0 . \nonumber
\end{eqnarray}
It is the multiplicity with which $e(\frac{a}{m}-\frac{b}{n})$
gives a fixed unit root of order $p^k$ if $(a,b)$ runs through
$(\Z/m\Z)^*\times (\Z/n\Z)^*$.

(i) If neither $\frac{m}{n}$ nor $\frac{n}{m}$ is a power of a prime number, then
\begin{eqnarray}\label{4.20}
\Lambda(m,n,p,k)&=&0.
\end{eqnarray}

(ii) Suppose $\frac{m}{n}=q^l$ for a prime number $q$ and some $l\in\N$.
Then 
\begin{eqnarray}\label{4.21}
\Lambda(n,m,p,k)=\Lambda(m,n,p,k)=\left\{\begin{array}{ll}
0 & \textup{if }(p,k)\neq (q,v_p(m)),\\
\varphi(n) & \textup{if }(p,k)=(q,v_p(m)).\end{array}\right. 
\end{eqnarray}

(iii) Suppose $m=n$. Then
\begin{eqnarray}\label{4.22}
\Lambda(m,m,p,k)&=& \left\{\begin{array}{ll}
0 & \textup{if }v_p(m)<k,\\
\varphi(m)\cdot\frac{p-2}{p-1} & \textup{if }v_p(m)=k,\\
\varphi(m) & \textup{if }v_p(m)>k.\end{array}\right. 
\end{eqnarray}

(d) \cite{Ap70} For $m,n\in\N$,
\begin{eqnarray}
\Res(\Phi_m,\Phi_n)&=& 0\qquad\textup{if }m=n.\label{4.23}\\
\Res(\Phi_m,\Phi_n)&=& 1\qquad\textup{if neither }\frac{m}{n}\textup{ nor }
\frac{n}{m} \nonumber\\
&& \quad \textup{ is a power of a prime number}.\label{4.24}\\
\Res(\Phi_{p^kn},\Phi_n)&=& \Res(\Phi_n,\Phi_{p^kn})=p^{\varphi(n)}
\quad\textup{if }p\textup{ is a prime number}\nonumber\\
&&\quad \textup{ and }k\in\N
\textup{ and }(p,k,n)\neq (2,1,1).\label{4.25}\\
\Res(\Phi_1,\Phi_2)&=&-\Res(\Phi_2,\Phi_1)=2.\label{4.26}
\end{eqnarray}

(e) For $m\in\N$ $\discr(\Phi_m)\in\N$. 
For any prime number $p$
\begin{eqnarray}\label{4.27} 
v_p(\discr(\Phi_m)) &=& \left\{\begin{array}{ll}
0 & \textup{if }v_p(m)=0,\\
(v_p(m)-\frac{1}{p-1})\cdot\varphi(m) & \textup{if }v_p(m)\geq 1.
\end{array}\right. 
\end{eqnarray}
\end{theorem}

{\bf Proof:}
(a) and (b) and (d) \cite[Theorem 3.1]{He20}.

(c) The parts (i) and (ii) are proved in \cite{He20} as part of 
the proof of part (d). Therefore here we prove only part (iii). 

Consider $m\in\Z_{\geq 2}$, $k\in\N$, and a prime number $p$.
If $v_p(m)<k$ (for example if $v_p(m)=0$) then the reduced denominator
$\beta$ of $\frac{a-b}{m}=\frac{\alpha}{\beta}$
with $\gcd(\alpha,\beta)=1$ is never equal to $p^k$. Thus then
$\Lambda(m,m,p,k)=0$. 

Suppose now $k\leq v_p(m)$. If $a$ and $b$ run through $(\Z/m\Z)^*$,
then the classes $[a]_{m/p^k}\in\Z/(m/p^k)\Z$ and $[b]_{m/p^k}\in\Z/(m/p^k)\Z$
run both through $(\Z/(m/p^k)\Z)^*$ with multiplicity 
$\varphi(m)/\varphi(m/p^k)$. Therefore the class $[a-b]_{m/p^k}$
in $\Z/(m/p^k)\Z$ is zero in 
\begin{eqnarray*}
\frac{\varphi(m)}{\varphi(m/p^k)}
\cdot \frac{\varphi(m)}{\varphi(m/p^k)} \cdot \varphi(m/p^k)
&=& \frac{\varphi(m)^2}{\varphi(m/p^k)}
\end{eqnarray*}
cases. This is the number of cases where the reduced denominator 
of $\frac{a-b}{m}$ divides $p^k$.
Therefore the reduced denominator of $\frac{a-b}{m}$ is 
equal to $p^k$ in 
\begin{eqnarray}
&&\frac{\varphi(m)^2}{\varphi(m/p^k)}-\frac{\varphi(m)^2}{\varphi(m/p^{k-1})}
\nonumber \\
&=& \left\{\begin{array}{ll}
\varphi(p^k)\cdot\varphi(m) & \textup{if }v_p(m)>k\\
\varphi(p^k)\cdot\varphi(m)\cdot \frac{p-2}{p-1} & \textup{if }v_p(m)=k
\end{array} \right. \label{4.28}
\end{eqnarray}
cases. In these cases, there are $\varphi(p^k)$ possibilities 
for the reduced numerator of $\frac{a-b}{m}$. This shows \eqref{4.22}.

(e) $\discr(\Phi_m)\in\Z$ because $\Phi_m\in\Z[t]$. 
And $\discr(\Phi_m)>0$ because the zeros $a_i$ in formula \eqref{4.9}
are here the primitive $m$-th unit roots, and with $\lambda$ also
$\overline\lambda$ is such a unit root.
By formula \eqref{4.9}
\begin{eqnarray}\label{4.29}
\discr(\Phi_m) &=& \prod_{(a,b)\in ((\Z/m\Z)^*)^2}
e(\frac{b}{m})\cdot(1-e(\frac{a-b}{m})).
\end{eqnarray}
Recall \eqref{4.18} $\norm_{\ord(\lambda)}(1-\lambda)=\Phi_{\ord(\lambda)}(1)$
for any unit root $\lambda$, and recall Theorem \ref{t4.5} (a).
The right hand side of \eqref{4.29} can be seen as a product of
unit roots and such norms $\norm_{\ord(\lambda)}(1-\lambda)$ for suitable $\lambda$.
Only $\ord(\lambda)=p^k$ for $k\geq 1$ contributes to $v_p(\discr(\Phi_m))$.
The precise amount of this contribution can be read off from formula
\eqref{4.22}. Thus
\begin{eqnarray}
v_p(\discr(\Phi_m)) &=& 0\hspace*{2cm} \textup{if }v_p(m)=0,\nonumber \\
v_p(\discr(\Phi_m)) &=& \sum_{k=1}^{v_p(m)-1}\varphi(m) + 
\varphi(m) \cdot \frac{p-2}{p-1} \label{4.30}\\
&=& (v_p(m)-\frac{1}{p-1})\cdot \varphi(m)\quad\textup{if }v_p(m)>0.
\hspace*{1cm}\Box \nonumber 
\end{eqnarray}

\section{Different decompositions into companion blocks}\label{c5}
\setcounter{equation}{0}

Theorem \ref{t5.1} is after Theorem \ref{t3.1} our second
structural result about decompositions of a $\Z$-lattice $H$
with endomorphism $h$ into companion blocks. 
Now the focus is on arbitrary decompositions, not just the standard decomposition.
Theorem \ref{t5.1} has some similarity with the chinese remainder theorem. 
We work with the companion blocks $H^{[p]}$ for $p\in\Z[t]$ unitary
from Definition \ref{t2.6}.
Part (a) of Theorem \ref{t5.1} is Lemma 6.2 (b) in \cite{He20},
part (b) is new.

\noindent
\begin{theorem}\label{t5.1}
(a) Let $f,g\in\Z[t]$ be unitary polynomials with 
$\Res(f,g)\neq 0$ (equivalent is $gcd(f,g)_{\Q[t]}=1$).
Then 
\begin{eqnarray}\label{5.1}
H^{[fg]}\cong H^{[f]}\oplus H^{[g]}\iff |\Res(f,g)|=1.
\end{eqnarray}
(b) Let $f_1,f_2,f_3,f_4\in\Z[t]$ be unitary polynomials with
$\Res(f_i,f_j)\neq 0$ for all $i\neq j$. Then
\begin{eqnarray}\label{5.2}
&&H^{[f_1f_3f_4]}\oplus H^{[f_2f_3]}\cong H^{[f_1f_3]}\oplus H^{[f_2f_3f_4]}\\
&\iff& |\Res(f_1,f_4)| = |\Res(f_2,f_4)|=1. \nonumber 
\end{eqnarray}
\end{theorem}

{\bf Proof:}
(a) Recall \eqref{4.13}, $|\Res(f,g)|=1\iff(f,g)=\Z[t]$. 
By Lemma \ref{t2.8} (a), the ideals $(g)/(fg)$ and $(f)/(fg)$ in $H^{[fg]}$ 
are isomorphic to the companion blocks $H^{[f]}$ and $H^{[g]}$, respectively,
and they are the unique primitive sublattices in $H^{[fg]}$ 
which are monodromy invariant and have the characteristic polynomials $f$ and $g$.
Because of $\gcd(f,g)_{\Q[t]}=1$, their intersection is 0, so they form a direct sum
within $H^{[fg]}$. 
This sum is $(f,g)/(fg)\subset \Z[t]/(fg)=H^{[fg]}$. 
Therefore $H^{[fg]}\cong H^{[f]}\oplus H^{[g]}$ is equivalent
to $(f,g)/(fg)=H^{[fg]}$, and this is equivalent to $(f,g)=\Z[t]$.

\medskip
(b) In the case of a sum of companion blocks, 
the monodromy $h$ of the sum is defined to be the sum of the monodromies 
of the single companion blocks.

Suppose that the isomorphism in the first line of \eqref{5.2} holds.
Divide both sides by the kernel $\ker(f_2(h)f_3(h))$. Then one obtains an 
isomorphism
$$H^{[f_1f_4]}\cong H^{[f_1]}\oplus H^{[f_4]}.$$
Part (a) shows $|\Res(f_1,f_4)|=1$. The necessity of $|\Res(f_2,f_4)|=1$ is
obtained analogously.

It remains to show that these two conditions are sufficient for the isomorphism.
So suppose  $|\Res(f_1,f_4)|=1=|\Res(f_2,f_4)|$. 

We identify $H^{[f_1f_3f_4]}\oplus  H^{[f_2f_3]}$ 
with $H^{[f_1f_3f_4]}\times H^{[f_2f_3]}$. 
We denote its monodromy by $h=(h_1,h_2)$, where 
$h_1=h_{[f_1f_3f_4]}$ is the monodromy of $H^{[f_1f_3f_4]}$, 
and $h_2=h_{[f_2f_4]}$ is the monodromy of $H^{[f_2f_4]}$.
Let $a_1\in H^{[f_1f_3f_4]}$ be a cyclic generator of it, 
and let $a_2\in H^{[f_2f_3]}$ be a cyclic generator of it.

$|\Res(f_1,f_4)|=1$ and $|\Res(f_2,f_4)|=1$ imply $|\Res(f_1f_2,f_4)|=1$,
and this implies the existence of polynomials $g_1,g_4\in\Z[t]$
with $g_1f_1f_2-g_4f_4=1$. Observe that the following matrix has determinant 1 
and thus is in $GL(2,\Z[t])$, and that its inverse is as follows,
\begin{eqnarray*}
\begin{pmatrix}g_1f_1 & 1 \\ g_4f_4 & f_2\end{pmatrix}\in GL(2,\Z[t]),
\qquad \begin{pmatrix}g_1f_1 & 1 \\ g_4f_4 & f_2\end{pmatrix}^{-1}
=\begin{pmatrix}f_2 & -1 \\ -g_4f_4 & g_1f_1\end{pmatrix} .
\end{eqnarray*}
Consider the elements
\begin{eqnarray*}
b_1&:=&((g_1f_1)(h_1)(a_1), a_2)\in H^{[f_1f_3f_4]}\times H^{[f_2f_3]},\\
b_2&:=&((g_4f_4)(h_1)(a_1),(f_2)(h_2)(a_2))\in H^{[f_1f_3f_4]}\times H^{[f_2f_3]}.
\end{eqnarray*}

$b_1$ is the generator of an Orlik block $B_1$  
whose characteristic polynomial divides $f_2f_3f_4$.
And $b_2$ is the generator of an Orlik block $B_2$ 
whose characteristic polynomial divides $f_1f_3$.

It remains to show $B_1+B_2=H^{[f_1f_3f_4]}\times H^{[f_2f_3]}$
because then by comparison of ranks one obtains $B_1+B_2=B_1\oplus B_2$ 
and that the characteristic polynomial of $B_1$ is $f_2f_3f_4$
and the characteristic polynomial of $B_2$ is $f_1f_3$.
Thus it remains to show that $(a_1,0)$ and $(0,a_2)$ are in 
$B_1+B_2$.
Calculate 
\begin{eqnarray*}
&& f_2(h)(b_1) - b_2\\
&=& ((f_2g_1f_1)(h_1)(a_1),(f_2)(h_2)(a_2))\\
&& - ((g_4f_4)(h_1)(a_1),(f_2)(h_2)(a_2))\\
&=& (a_1,0)\in B_1+B_2
\end{eqnarray*}
and
\begin{eqnarray*} 
&& (-g_4f_4)(h)(b_1) + (g_1f_1)(h)(b_2)\\
&=& ((-g_4f_4g_1f_1)(h_1)(a_1),(-g_4f_4)(h_2)(a_2))\\
&& + ((g_1f_1g_4f_4)(h_1)(a_1),(g_1f_1f_2)(h_2)(a_2))\\
&=& (0,a_2)\in B_1+B_2.\hspace*{4cm} \Box 
\end{eqnarray*}

\begin{remark}\label{t5.2}
We expect (but we did not prove it) that the following holds.
Let $f_1,...,f_a,g_1,...,g_b\in\Z[t]$ 
be products of cyclotomic polynomials with no multiple roots.
If an isomorphism 
\begin{eqnarray*}
\bigoplus_{i=1}^a H^{[f_i]} \cong \bigoplus_{j=1}^b H^{[g_j]}
\end{eqnarray*} 
holds, then it can be deduced by repeated application of the
rule \eqref{5.2} and adding to both sides the same Orlik blocks.
This property would say that the equivalence in \eqref{5.2} 
would be the most general rule for getting isomorphisms of
sums of Orlik blocks.
\end{remark}

\section{When does a power of an Orlik block admit a standard
decomposition into Orlik blocks?}\label{c6}
\setcounter{equation}{0}

\noindent
Theorem \ref{t6.2} starts with one Orlik block $(H,h)$ and a number 
$\mu\in\N$ and gives a sufficient criterion for $(H,h^\mu)$
to admit a standard decomposition into Orlik blocks. 
It will be crucial for proving Orlik's Conjecture \ref{t1.2} 
in the case of the chain type singularities. 
The condition will work with a graph whose set of vertices is the
set $M\subset \N$ of orders of the eigenvalues of $h:H_\C\to H_\C$.
Theorem \ref{t6.2} is preceded by some definitions and observations.

\begin{definition}\label{t6.1}
(a) Recall that $\PP\subset \N$ denotes the set of prime 
numbers. Consider the infinite directed graph $(\N,E)$ 
whose set of vertices is $\N$ und whose set 
$E\subset \N^2$ of directed edges is defined as follows,
\begin{eqnarray}
E_p&:=& \{(m,n)\in\N^2\, |\, 
\frac{m}{n}=p^k\textup{ for some }k\in\N\} 
\textup{ for any }p\in\PP,\nonumber\\
E&:=& \bigcup_{p\in \PP} E_p.\label{6.1}
\end{eqnarray}
An edge in $E_p$ is called a {\it $p$-edge}.

\medskip
(b) For any finite set $M\subset\N$ consider the directed graph
$(M,E(M))$ which is the restriction of $(\N,E)$ to $M$, so
its set of directed edges is $E(M)=E\cap M^2$.

\medskip
(c) For any $\mu\in\N$ define the map
\begin{eqnarray}\label{6.2}
\gamma_\mu:\N\to\N,\\
m\mapsto \frac{m}{\gcd(m,\mu)}.\nonumber
\end{eqnarray}

\medskip
(d) Consider a finite set $M\subset\N$ and a number 
$\mu\in\N$. The pair $(M,\mu)$ is called 
{\it sdiOb-sufficient} (sdiOb for {\it standard decomposition
into Orlik blocks}) if for any prime number $p$ and 
any $p$-edge $(n_a,n_b)$ in
$E_p(\gamma_\mu(M))$ at least one of the following 
two conditions holds.
\begin{eqnarray}\label{6.3}
E_p\cap \left((\gamma_\mu^{-1}(n_a)\cap M)\times 
\gamma_\mu^{-1}(n_b)\right)
&\subset& E_p(M),\\
E_p\cap \left(\gamma_\mu^{-1}(n_a)\times 
(\gamma_\mu^{-1}(n_b)\cap M)\right)
&\subset& E_p(M)\label{6.4}
\end{eqnarray}
(these two conditions are discussed in Remark \ref{t6.3} (iv)).
\end{definition}

\begin{theorem}\label{t6.2}
Consider a finite non-empty set $M\subset\N$, the
corresponding Orlik block $\Or(M)=H$ with monodromy $h$,
and a number $\mu\in\N$. Then $(H,h^\mu)$ admits
a standard decomposition into Orlik blocks 
if $(M,\mu)$ is sdiOb-sufficient.
\end{theorem}

Before proving this theorem, we make some elementary observations.

\begin{remarks}\label{t6.3}
(i) We expect that {\it if and only if} holds 
in Theorem \ref{t6.2}.


\medskip
(ii) For any $n\in\N$, the fiber 
$\gamma_\mu^{-1}(n)\subset\N$ is finite and nonempty. It is 
\begin{eqnarray}\label{6.5}
\gamma_\mu^{-1}(n)= \{n\cdot c\cdot \prod_{p\in\PP:\, v_p(n)>0}p^{v_p(\mu)}\,|\, 
c\textup{ divides }\prod_{p\in\PP:\, v_p(n)=0}p^{v_p(\mu)}\}.
\end{eqnarray}
Especially, if $v_p(n)>0$ for some prime number $p$, then
$v_p(m)=v_p(n)+v_p(\mu)$ for any $m\in \gamma_\mu^{-1}(n)$.

\medskip
(iii) Consider $(n_a,n_b)\in\N^2$. 
If $v_p(n_a)>v_p(n_b)$ for some prime number $p$, then 
for any $(m_c,m_d)\in \gamma_\mu^{-1}(n_a)\times \gamma_\mu^{-1}(n_b)$
\begin{eqnarray}\label{6.6}
v_p(m_c)= v_p(n_a)+v_p(\mu) > v_p(n_b)+v_p(\mu)\geq v_p(m_d).
\end{eqnarray}

Thus, if $(n_a,n_b)$ is no edge, then there is no edge in 
$\gamma_\mu^{-1}(n_a)\times \gamma_\mu^{-1}(n_b)$.
And if $(n_a,n_b)$ is a $p$-edge, then any edge in
$\gamma_\mu^{-1}(n_a)\times \gamma_\mu^{-1}(n_b)$ 
is a $p$-edge. 

If $(n_a,n_b)$ is a $p$-edge and $v_p(n_b)>0$, then
the map
\begin{eqnarray}\label{6.7}
\gamma_\mu^{-1}(n_a)&\to& \gamma_\mu^{-1}(n_b)\\
m&\mapsto& m\cdot \frac{n_b}{n_a}=m\cdot p^{-v_p(n_a)+v_p(n_b)},\nonumber
\end{eqnarray} 
is a bijection and the set of $p$-edges in $\gamma_\mu^{-1}(n_a)\times \gamma_\mu^{-1}(n_b)$
is 
\begin{eqnarray}\label{6.8}
\{(m,m\cdot \frac{n_b}{n_a})\, |\, m\in\gamma_\mu^{-1}(n_a)\}.
\end{eqnarray}

If $(n_a,n_b)$ is a $p$-edge and $v_p(n_b)=0$, then
\begin{eqnarray}\label{6.9}
\gamma_\mu^{-1}(n_b)&=& \bigcup_{m\in\gamma_\mu^{-1}(n_a)}
\{m\cdot p^{-v_p(m)+k}\, |\, k\in\{0,...,v_p(\mu)\}\},
\end{eqnarray} 
and the set of $p$-edges in $\gamma_\mu^{-1}(n_a)\times \gamma_\mu^{-1}(n_b)$
is 
\begin{eqnarray}\label{6.10}
\bigcup_{m\in\gamma_\mu^{-1}(n_a)}\{m\}\times 
\{(m\cdot p^{-v_p(m)+k}\, |\, k\in\{0,...,v_p(\mu)\}\}.
\end{eqnarray}

\medskip
(iv) In the case $v_p(n_b)>0$ \eqref{6.3} says that for any 
$m\in\gamma_\mu^{-1}(n_a)\cap M$ the number $m\cdot\frac{n_b}{n_a}$
is in $\gamma_\mu^{-1}(n_b)\cap M$, as then
$(m,m\cdot\frac{n_b}{n_a})$ is the only $p$-edge in
$\{m\}\times \gamma_\mu^{-1}(n_b)$.

In the case $v_p(n_b)=0$ \eqref{6.3} says that for any 
$m\in\gamma_\mu^{-1}(n_a)\cap M$ the set 
$\{m\cdot p^{-v_p(m)+k}\, |\, k\in\{0,...,v_p(\mu)\}\}$
is a subset of $\gamma_\mu^{-1}(n_b)\cap M$, as then
$\{m\}\times \{m\cdot p^{-v_p(m)+k}\, |\, 
k\in\{0,...,v_p(\mu)\}\}$ is the set of $p$-edges in
$\{m\}\times \gamma_\mu^{-1}(n_b)$.

\eqref{6.4} says that for any $m\in \gamma_\mu^{-1}(n_b)\cap M$
the number $m\cdot p^{-v_p(m)+v_p(n_a)+v_p(\mu)}$
is in $\gamma_\mu^{-1}(n_a)\cap M$, as then 
$(m\cdot p^{-v_p(m)+v_p(n_a)+v_p(\mu)},m)$ is the only $p$-edge
in $\gamma_\mu^{-1}(n_a)\times\{m\}$.

\medskip
(v) In \cite{He20} the group of automorphisms with eigenvalues
in $S^1$ of an Orlik block $(H,h)$ was studied.
Theorem 1.2 in \cite{He20} gives a necessary and sufficient criterion for 
this group to be only $\{\pm h^k\, |\, k\in\Z\}$. The criterion also uses the
graph $(M,E(M))$. In fact, a $p$-edge $(m_a,m_b)$ here is called a $p$-edge
there only if no $m_c\notin\{m_a,m_b\}$ with $m_b|m_c|m_a$ exists.
The purpose of that restriction in \cite{He20} was mainly to have graphs with not
too many edges. The conditions in Theorem 1.2 in \cite{He20} work
also with $(M,E(M))$.
\end{remarks}

{\bf Proof of Theorem \ref{t6.2}:}
Consider the pair $(H,h^\mu)$. Let $p_1,...,p_k\in\Z[t]$ be the unique
unitary polynomials with $p_{H,h^\mu}=p_1\cdot ...\cdot p_k$
and $p_k|p_{k-1}|...|p_2|p_1$. If $(H,h^\mu)$ admits a standard decomposition
into Orlik blocks, that is isomorphic to $\bigoplus_{j=1}^k H^{[p_j]}$.

We want to apply Theorem \ref{t3.1} to $(H,h^\mu)$ instead of $(H,h)$.

Let $a_0\in H$ be a generating element of the Orlik block $(H,h)$,
so $H=\bigoplus_{i=0}^{\rk H-1}\Z\cdot h^i(a_0)$. Define 
\begin{eqnarray}\label{6.11}
a_i&:=& h^i(a_0)\quad\textup{for }i\geq 0,\\
a_0^{[j]}&:=& a_{j-1}\quad \textup{for }j\in\{1,...,k\},\label{6.12}\\
a_i^{[j]}&:=& (h^\mu)^i(a_0^{[j]})=a_{j-1+\mu\cdot i}
\quad \textup{for }i\geq 0\textup{ and }j\in\{1,...,k\},\label{6.13}\\
{\bf a}^{[j]}&:=& (a_0^{[j]},a_1^{[j]},...,a_{\deg p_j-1}^{[j]}).\label{6.14}
\end{eqnarray}

We will show that the tuple 
\begin{eqnarray}\label{6.15}
{\bf a}^{dec}&:=& (a_i^{[j]}\, |\, j\in\{1,...,k\},i\in\{0,...,\deg p_j-1\})\\
&=& ({\bf a}^{[1]},{\bf a}^{[2]},...,{\bf a}^{[k]})\nonumber
\end{eqnarray}
in \eqref{3.1} is a $\Z$-basis of $H$ if and only if $(M,\mu)$ is sdiOb-sufficient. 
This and Theorem \ref{t3.1} show that $(H,h^\mu)$ admits a 
standard decomposition into Orlik blocks  if $(M,\mu)$ is sdiOb-sufficient.

It remains to show that the tuple ${\bf a}^{dec}$
is a $\Z$-basis of $H$
if and only if $(M,\mu)$ is sdiOb-sufficient.
The tuple
\begin{eqnarray}\label{6.16}
{\bf a}^{st}&:=& (a_0,a_1,...,a_{\rk H-1})
\end{eqnarray} 
is a $\Z$-basis of $H$. The tuple ${\bf a}^{dec}$ is a $\Z$-basis
if and only if the matrix $M({\bf a}^{st},{\bf a}^{dec})$ 
with ${\bf a}^{dec}={\bf a}^{st}\cdot M({\bf a}^{st},{\bf a}^{dec})$
(see the Notation \ref{t2.2}) has determinant $\pm 1$.
It remains to calculate this determinant. In order to do so, we consider
also certain tuples of eigenvectors of $h$ and $h^\mu$.

Let $\{\kappa_1,\kappa_2,...,\kappa_{\rk H}\}\subset\C$ be the set
of eigenvalues of $h$, ordered such that
\begin{eqnarray}\label{6.17}
\ord(\kappa_{\alpha})\leq \ord(\kappa_{\beta})
\quad\textup{if }\alpha<\beta,
\end{eqnarray} 
and let 
\begin{eqnarray}\label{6.18}
{\bf v}^{I}&=& (v_1,v_2,...,v_{\rk H})
\end{eqnarray}
be the tuple of eigenvectors $v_\alpha\in H_\C$ with 
\begin{eqnarray}\label{6.19}
h(v_\alpha)&=&\kappa_\alpha\cdot v_\alpha,\\
a_0 &=&\sum_{\alpha=1}^{\rk H}v_\alpha. \label{6.20}
\end{eqnarray}
Then $M({\bf v}^{I},{\bf a}^{st})$ with ${\bf a}^{st}={\bf v}^{I}
\cdot M({\bf v}^{I},{\bf a}^{st})$ is the Vandermonde matrix
\begin{eqnarray}\label{6.21}
M({\bf v}^{I},{\bf a}^{st})
&=& \begin{pmatrix}1&\kappa_1^1&\cdots&\kappa_1^{\rk H-1}\\
\vdots & \vdots & & \vdots \\
1 & \kappa_{\rk H}^1 & \cdots & \kappa_{\rk H}^{\rk H-1} 
\end{pmatrix} .
\end{eqnarray}

Let $\{\lambda_1,...,\lambda_{\deg p_1}\}$ be the set of 
eigenvalues of $h^\mu$, ordered such that
\begin{eqnarray}\label{6.22}
p_j(t)&=& \prod_{l=1}^{\deg p_j}(t-\lambda_l) 
\quad\textup{for }j\in\{1,...,k\} .
\end{eqnarray}
Define for $\beta\in\{1,2,...,\deg p_1\}$
the index set $A(\beta)$ by 
\begin{eqnarray}\label{6.23}
A(\beta)&:=&\{\alpha\in\{1,2,...,\rk H\}\, |\, \kappa_\alpha^\mu=\lambda_\beta\}\\
&=:& \{\alpha(\beta,1),\alpha(\beta,2),...,\alpha(\beta,|A(\beta)|)\}
\nonumber\\
&&\textup{ with }\alpha(\beta,1)<\alpha(\beta,2)<...<
\alpha(\beta,|A(\beta)|).\nonumber
\end{eqnarray}
The space $\bigoplus_{\alpha\in A(\beta)}\C\cdot v_\alpha\subset H_\C$
is the eigenspace with eigenvalue $\lambda_\beta$ of $h^\mu$.
For any $j\in\{1,...,k\}$, the vector 
\begin{eqnarray}\label{6.24}
v^{III,j,\beta}&:=& \sum_{\alpha\in A(\beta)}\kappa_\alpha^{j-1}
\cdot v_\alpha\in H_\C
\end{eqnarray}
is an eigenvector with eigenvalue $\lambda_\beta$ of $h^\mu$.
It is useful as for $i\geq 0$ and $j\in\{1,...,k\}$
\begin{eqnarray}
a_i^{[j]}&=& a_{j-1+\mu\cdot i} = \sum_{\alpha=1}^{rk H}
\kappa_\alpha^{j-1+\mu\cdot i}\cdot v_\alpha \nonumber\\
&=& \sum_{\beta=1}^{\deg p_1}\lambda_\beta^i\cdot 
\sum_{\alpha\in A(\beta)}\kappa_\alpha^{j-1}\cdot v_\alpha
= \sum_{\beta=1}^{\deg p_1}\lambda_\beta^i 
\cdot v^{III,j,\beta}.\label{6.25}
\end{eqnarray}
Consider for $j\in\{1,...,k\}$ and 
$\beta\in\{1,2,...,\deg p_1\}$
the following tuples of eigenvectors of $h$ and/or $h^\mu$,
\begin{eqnarray}\label{6.26}
{\bf v}^{II,\beta}&:=& (v_{\alpha(\beta,1)},
v_{\alpha(\beta,2)},...,
v_{\alpha(\beta,|A(\beta)|)}),\\
{\bf v}^{II}&:=& ({\bf v}^{II,1},...,{\bf v}^{II,\deg p_1}),
\label{6.27}\\
{\bf v}^{III,\beta}&:=& (v^{III,1,\beta},...,
v^{III,|A(\beta)|,\beta}),\label{6.28}\\
{\bf v}^{III}&:=& ({\bf v}^{III,1},...,{\bf v}^{III,\deg p_1}),
\label{6.29}\\
{\bf v}^{V,j}&:=& (v^{III,j,1},v^{III,j,2},...,v^{III,j,\deg p_1}),
\label{6.30}\\
{\bf v}^{IV,j}&:=& (v^{III,j,1},v^{III,j,2},...,
v^{III,j,\deg p_j}),
\label{6.31}\\
{\bf v}^{IV}&:=& ({\bf v}^{IV,1},...,{\bf v}^{IV,k}).\label{6.32}
\end{eqnarray}
${\bf v}^{II,\beta}$ and ${\bf v}^{III,\beta}$ are 
$\C$-bases of the eigenspace with eigenvalue $\lambda_\beta$
of $h^\mu$. The base change matrix 
$M({\bf v}^{II,\beta},{\bf v}^{III,\beta})$
rewrites the relation \eqref{6.24}. It is a Vandermonde matrix,
\begin{eqnarray}\label{6.33}
M({\bf v}^{II,\beta},{\bf v}^{III,\beta})
&=& \begin{pmatrix} 1&\kappa_{\alpha(\beta,1)}^1& ... &
\kappa_{\alpha(\beta,1)}^{|A(\beta)|-1}\\
\vdots & \vdots & & \vdots\\
 1&\kappa_{\alpha(\beta,|A(\beta)|)}^1& ...&
 \kappa_{\alpha(\beta,|A(\beta)|)}^{|A(\beta)|-1}
\end{pmatrix}.
\end{eqnarray}
${\bf v}^{I}$, ${\bf v}^{II}$, ${\bf v}^{III}$ and 
${\bf v}^{IV}$ are $\C$-bases of $H_\C$. 
The base change matrices $M({\bf v}^{I},{\bf v}^{II})$ 
and $M({\bf v}^{III},{\bf v}^{IV})$ are just permutation matrices. 
The entries of ${\bf v}^{V,j}$ are linearly independent.
Therefore the following rectangular 
$(\deg p_1)\times(\deg p_j)$-matrix is well defined.
It rewrites the relations \eqref{6.25}.
\begin{eqnarray}\label{6.34}
M({\bf v}^{V,j},{\bf a}^{[j]})
&=& \begin{pmatrix} 1&\lambda_1^1& ...&\lambda_1^{\deg p_j-1}\\
\vdots & \vdots & & \vdots\\
 1&\lambda_{\deg p_1}^1& ...&\lambda_{\deg p_1}^{\deg p_j-1}
\end{pmatrix}.
\end{eqnarray}
Now we want to describe the matrix
$M({\bf v}^{IV},{\bf a}^{dec})$. 
Observe that $v^{III,j,\beta}$ is in the case
$j>|A(\beta)|$ a linear combination of the entries of
${\bf v}^{III,\beta}$. Therefore the matrix
$M({\bf v}^{IV},{\bf a}^{dec})$ is a block upper triangular
matrix whose diagonal blocks are obtained from the
matrices in \eqref{6.34} by cutting off the lower lines.
The diagonal blocks are the Vandermonde matrices
\begin{eqnarray}\label{6.35}
M^{[j]}&:=&
\begin{pmatrix} 1&\lambda_1^1& ...&\lambda_1^{\deg p_j-1}\\
\vdots & \vdots & & \vdots\\
 1&\lambda_{\deg p_j}^1& ...&\lambda_{\deg p_j}^{\deg p_j-1}
\end{pmatrix}.
\end{eqnarray}
Now the matrix $M({\bf a}^{st},{\bf a}^{dec})$ can be written
as the product of matrices
\begin{eqnarray}\label{6.36}
M({\bf a}^{st},{\bf v}^{I})
M({\bf v}^{I},{\bf v}^{II})
M({\bf v}^{II},{\bf v}^{III})
M({\bf v}^{III},{\bf v}^{IV})
M({\bf v}^{IV},{\bf a}^{dec}).
\end{eqnarray}
The absolute value of its determinant is
\begin{eqnarray}
&&|\det M({\bf a}^{st},{\bf a}^{dec})| \nonumber\\
&=&|\frac{\det M({\bf v}^{II},{\bf v}^{III})\cdot
\det M({\bf v}^{IV},{\bf a}^{dec})}
{\det M({\bf v}^I,{\bf a}^{st})}|\nonumber\\
&=& |\frac{\prod_{\beta=1}^{\deg p_1}
\det M({\bf v}^{II,\beta},{\bf v}^{III,\beta})\cdot
\prod_{j=1}^{k}\det M^{[j]}}
{\det M({\bf v}^I,{\bf a}^{st})}| .\label{6.37}
\end{eqnarray}
It is the absolute value of a quotient of determinants
of Vandermonde matrices. The determinants are
\begin{eqnarray}\label{6.38}
\det M({\bf v}^{I},{\bf a}^{st})
&=& \prod_{1\leq \alpha_1<\alpha_2\leq \rk H}
(\kappa_{\alpha_2}-\kappa_{\alpha_1}),\\
\det M({\bf v}^{II,\beta},{\bf v}^{III,\beta})
&=& \prod_{\alpha_1,\alpha_2\in A(\beta):\, \alpha_1<\alpha_2}
(\kappa_{\alpha_2}-\kappa_{\alpha_1}),\label{6.39}\\
\det M^{[j]}
&=& \prod_{1\leq \beta_1<\beta_2\leq \deg p_j}
(\lambda_{\beta_2}-\lambda_{\beta_1}).\label{6.40}
\end{eqnarray}

Only now we use that $\lambda_\beta$ and $\kappa_\alpha$
are unit roots. For $n\in \gamma_\mu(M)$ we denote the multiplicity
as a zero of $p_{H,h^\mu}$ of any unit root $\lambda$ 
with $\ord(\lambda)=n$ by $\psi(n)\in\N_0$. It is 
\begin{eqnarray}\label{6.41}
\psi(n)&=& \varphi(n)^{-1}\cdot\sum_{m\in\gamma_\mu^{-1}(n)\cap M}
\varphi(m)\\
&=& \max(j\in\{1,...,k\}\, |\, \Phi_n|p_j)\leq k.\label{6.42}
\end{eqnarray}
In view of the formulas \eqref{4.3} and \eqref{4.9},
the determinants in \eqref{6.38} and \eqref{6.40}
can be written as follows as products of resultants
and discriminants,
\begin{eqnarray}
|\det M({\bf v}^I,{\bf a}^{st})|&=& 
\prod_{m_c,m_d\in M:\, m_c>m_d}|\Res(\Phi_{m_c},\Phi_{m_d})|
\nonumber\\
&&\cdot \prod_{m\in M}\sqrt{|\discr(\Phi_m)|}, \label{6.43}\\
|\prod_{j=1}^k\det M^{[j]}| &=& 
\prod_{n_a,n_b\in \gamma_\mu(M):\, n_a>n_b}
|\Res(\Phi_{n_a},\Phi_{n_b})|^{\min(\psi(n_a),\psi(n_b))}
\nonumber\\
&&\cdot \prod_{n\in \gamma_\mu(M)}\sqrt{|\discr(\Phi_n)|}^{\psi(n)}.
\label{6.44}
\end{eqnarray}
Thus we can rearrange $|\det M({\bf a}^{st},{\bf a}^{dec})|$
as a product of the following factors
in \eqref{6.45} and \eqref{6.46}:
For each pair $(n_a,n_b)\in \gamma_\mu(M)^2$ with $n_a>n_b$
\begin{eqnarray}\label{6.45}
\frac{|\Res(\Phi_{n_a},\Phi_{n_b})|^{\min(\psi(n_a),\psi(n_b))}}
{\prod_{(m_c,m_d)\in \gamma_\mu^{-1}(n_a)\times \gamma_\mu^{-1}(n_b)
\cap M^2} |\Res(\Phi_{m_c},\Phi_{m_d})|}.
\end{eqnarray}
And for each $n\in \gamma_\mu(M)$
\begin{eqnarray}\label{6.46}
\frac{\sqrt{|\discr(\Phi_n)|}^{\psi(n)}\cdot 
\prod_{\beta:\, \ord(\lambda_\beta)=n}
\prod_{\alpha_1,\alpha_2\in A(\beta):\, \alpha_1<\alpha_2} 
|\kappa_{\alpha_2}-\kappa_{\alpha_1}|}
{\prod_{m\in\gamma_\mu^{-1}(n)\cap M}\sqrt{|\discr(\Phi_m)|}
\cdot \prod_{m_c,m_d\in \gamma_\mu^{-1}(n)\cap M:\, m_c>m_d}
|\Res(\Phi_{m_c},\Phi_{m_d})|} .
\end{eqnarray}
We will now prove the following {\bf claims:}
\begin{list}{}{}
\item[(A)]
Each factor of type \eqref{6.45} with 
$(n_a,n_b)\notin E(\gamma_\mu(M))$ is equal to 1.
\item[(B)]
Each factor in \eqref{6.45} with 
$(n_a,n_b)\in E_p(\gamma_\mu(M))$ for some prime number $p$ 
is a positive integer, and it is equal to 1 
if and only if \eqref{6.3} or \eqref{6.4} is satisfied.
\item[(C)]
Each factor of type \eqref{6.46} is equal to 1.
\end{list}

Together (A), (B) and (C) give that $|\det M({\bf a}^{st},
{\bf a}^{dec})|$ is equal to 1 if and only if 
$(M,\mu)$ is sdiOb-sufficient.
Therefore the tuple in \eqref{3.1} with $a_0^{[j]}$ as in 
\eqref{6.12} is a $\Z$-basis
if and only if $(M,\mu)$ is sdiOb-sufficient.

It remains to prove the claims (A), (B) and (C).

{\bf Claim (A):} If $(n_a,n_b)$ is not an edge, then
$\Res(\Phi_{n_a},\Phi_{n_b})=1$ because of \eqref{4.24}.
By Remark \ref{t6.3} (iii), 
then any pair $(m_c,m_d)\in 
\gamma_\mu^{-1}(n_a)\times \gamma_\mu^{-1}(n_b)$
is also not an edge, and again
$\Res(\Phi_{m_c},\Phi_{m_d})=1$ because of \eqref{4.24}.

{\bf Claim (B):} Let $(n_a,n_b)$ be a $p$-edge of $\gamma_\mu(M)$.
Then by Remark \ref{t6.3} (iii),
also any edge in $\gamma_\mu^{-1}(n_a)\cap \gamma_\mu^{-1}(n_b)$ 
is a $p$-edge. Therefore and by Theorem \ref{t4.5} (d),
the numerator and the denominator in \eqref{6.45} are 
powers of $p$,
\begin{eqnarray}\label{6.47}
|\Res(n_a,n_b)|&=& p^{\varphi(n_b)},\\
|\Res(m_c,m_d)|&=& \left\{\begin{array}{ll}
1 & \textup{if }(m_c,m_d))\notin E_p(M),\\
p^{\varphi(m_d)} & \textup{if }(m_c,m_d)\in E_p(M).
\end{array}\right. \label{6.48}
\end{eqnarray}
We claim
\begin{eqnarray}\label{6.49}
\frac{|\Res(\Phi_{n_a},\Phi_{n_b})|^{\psi(n_a)}}
{\prod_{(m_c,m_d)\in E_p\cap((\gamma_\mu^{-1}(n_a)\cap M)
\times \gamma_\mu^{-1}(n_b))} 
|\Res(\Phi_{m_c},\Phi_{m_d})|}&=&1,\\
\frac{|\Res(\Phi_{n_a},\Phi_{n_b})|^{\psi(n_b)}}
{\prod_{(m_c,m_d)\in E_p(\gamma_\mu^{-1}(n_a)\times 
(\gamma_\mu^{-1}(n_b)\cap M))} 
|\Res(\Phi_{m_c},\Phi_{m_d})|}&=&1.\label{6.50}
\end{eqnarray}
The denominator in \eqref{6.49} is a multiple of the
denominator in \eqref{6.45}, and they are equal if and only if 
\eqref{6.3} holds. 
The denominator in \eqref{6.50} is a multiple of the
denominator in \eqref{6.45}, and they are equal if and only if 
\eqref{6.4} holds. 
Therefore the quotient in \eqref{6.45} is equal to 1 if and
only if \eqref{6.3} or \eqref{6.4} hold.
It remains to prove \eqref{6.49} and \eqref{6.50}.

We use the Remarks \ref{t6.3} (iii) and (iv). 
First suppose $v_p(n_b)>0$. Then
\begin{eqnarray}
&&E_p\cap((\gamma_\mu^{-1}(n_a)\cap M)\times \gamma_\mu^{-1}(n_b))
\nonumber\\
&=& \{(m,m\cdot\frac{n_b}{n_a})\,|\, 
m\in \gamma_\mu^{-1}(n_a)\cap M\},\label{6.51}
\end{eqnarray}
\begin{eqnarray}
&&v_p(\textup{denominator of \eqref{6.49}}) 
= \sum_{m\in \gamma_\mu^{-1}(n_a)\cap M}
\varphi(m\cdot\frac{n_b}{n_a}))\nonumber \\
&=& \sum_{m\in \gamma_\mu^{-1}(n_a)\cap M}
\varphi(m)\cdot\frac{\varphi(n_b)}{\varphi(n_a)}
= \psi(n_a)\cdot\varphi(n_b)\nonumber \\
&=& v_p(\textup{numerator of \eqref{6.49}}).\label{6.52}
\end{eqnarray}
This shows \eqref{6.49} in the case $v_p(n_b)>0$.

Now suppose $v_p(n_b)=0$. Then 
\begin{eqnarray}
&&E_p\cap((\gamma_\mu^{-1}(n_a)\cap M)\times \gamma_\mu^{-1}(n_b))
\label{6.53} \\
&=& \bigcup_{m\in \gamma_\mu^{-1}(n_a)\cap M}
\{m\}\times\{m\cdot p^{-v_p(m)+k}\, |\, 
k\in\{0,...,v_p(\mu)\}\},\nonumber
\end{eqnarray}
\begin{eqnarray}
&&v_p(\textup{denominator of \eqref{6.49}}) \nonumber \\ 
&=& \sum_{m\in \gamma_\mu^{-1}(n_a)\cap M}
\sum_{k=0}^{v_p(\mu)}\varphi(m\cdot p^{-v_p(m)+k})\nonumber\\
&=& \sum_{m\in \gamma_\mu^{-1}(n_a)\cap M} 
\varphi(m\cdot p^{-v_p(m)})\cdot
\left(1+(p-1)\sum_{k=1}^{v_p(\mu)}p^{k-1}\right)\nonumber \\
&=& \sum_{m\in \gamma_\mu^{-1}(n_a)\cap M} 
\frac{\varphi(m)}{\varphi(p^{v_p(m)})}\cdot p^{v_p(\mu)}
= \frac{\psi(n_a)\cdot\varphi(n_a)}
{\varphi(p^{v_p(n_a)+v_p(\mu)})}\cdot p^{v_p(\mu)}\nonumber\\
&=& \psi(n_a)\cdot \varphi(n_b) 
= v_p(\textup{numerator of \eqref{6.49}}).\label{6.54}
\end{eqnarray}
This shows \eqref{6.49} in the case $v_p(n_b)=0$.

In both cases, $v_p(n_b)>0$ or $v_p(n_b)=0$, we have 
\begin{eqnarray}
&&E_p\cap(\gamma_\mu^{-1}(n_a)\times(\gamma_\mu^{-1}(n_b)\cap M))
\nonumber\\
&=& \{(m\cdot p^{-v_p(m)+v_p(n_a)+v_p(\mu)},m)\,|\, 
m\in \gamma_\mu^{-1}(n_b)\cap M\},\label{6.55}
\end{eqnarray}
\begin{eqnarray}
&&v_p(\textup{denominator of \eqref{6.50}}) 
= \sum_{m\in \gamma_\mu^{-1}(n_b)\cap M}
\varphi(m)\nonumber \\
&=& \psi(n_b)\cdot\varphi(n_b)
= v_p(\textup{numerator of \eqref{6.50}}).\label{6.56}
\end{eqnarray}
This shows \eqref{6.50}.

{\bf Claim (C):} The squares of the numerator and of the 
denominator of \eqref{6.46} are positive integers.
Fix a prime number $p$. We will show 
\begin{eqnarray}\label{6.57}
v_p((\textup{numerator of \eqref{6.46}})^2)
&=& v_p((\textup{denominator of \eqref{6.46}})^2).\hspace*{1cm}
\end{eqnarray}

The second factor in the numerator is the most
difficult part. It is the product of $|\kappa_{\alpha_2}-\kappa_{\alpha_1}|$
over the pairs $(\alpha_2,\alpha_1)$ in the set
\begin{eqnarray}\label{6.58}
\{(\alpha_2,\alpha_1)\, |\, 1\leq \alpha_1<\alpha_2\leq \rk H,\,  
\ord(\kappa_{\alpha_i}^\mu)=n,\, (\kappa_{\alpha_2}/\kappa_{\alpha_1})^\mu=1\}.
\end{eqnarray}
For pairs $(\alpha_2,\alpha_1)$ in \eqref{6.58} we denote
$m_2:=\ord(\kappa_{\alpha_2})$ and $m_1:=\ord(\kappa_{\alpha_1})$.
In order to understand their contribution to 
$v_p((\textup{numerator of \eqref{6.46}})^2)$, 
we have to consider Theorem \ref{t4.5} (c) and 
\begin{eqnarray}\label{6.59}
\norm_{\ord(\kappa)}(1-\kappa) =
\left\{\begin{array}{ll}
q&\textup{if }\ord(\kappa)=q^l\textup{ for some }q\in\PP,l\geq 1,\\
1&\textup{else} \end{array}\right. 
\end{eqnarray}
(see \eqref{4.18} and Theorem \ref{t4.5} (a)).
Only pairs with 
$\ord(\kappa_{\alpha_2}/\kappa_{\alpha_1})=p^k$ for some
$k\geq 1$ and with $p^k|\mu$  
give a contribution to $v_p((\textup{numerator of \eqref{6.46}})^2)$.
Its size for any $k\geq 1$ is  $2\Lambda(m_2,m_1,p,k)$ if $m_2\neq m_1$
and $\Lambda(m_2,m_2,p,k)$ if $m_2=m_1$.

If $v_p(n)>0$, pairs $(\alpha_2,\alpha_1)$ in \eqref{6.58} 
with $m_2\neq m_1$ give no contribution as $v_p(m_2)=v_p(m_1)(=v_p(n)+v_p(\mu))$
by Remark \ref{t6.3} (ii), 
and thus neither $m_2/m_1$ nor $m_1/m_2$ is a power of $p$.
Pairs $(\alpha_2,\alpha_1)$ in \eqref{6.58} with $m_2= m_1$
and $\ord(\kappa_{\alpha_2}/\kappa_{\alpha_1})=p^k$ 
satisfy $k\leq v_p(\mu)< v_p(m_2)$ because of 
$(\kappa_{\alpha_2}/\kappa_{\alpha_1})^\mu=1$. Therefore 
and because of the third line of \eqref{4.22}, the contribution
of the square of the second factor in the numerator is in the case $v_p(n)>0$
\begin{eqnarray}\label{6.60}
\sum_{m\in\gamma_\mu^{-1}(n)\cap M}\sum_{k=1}^{v_p(\mu)}\varphi(m)
=v_p(\mu)\cdot \psi(n)\cdot\varphi(n).
\end{eqnarray}
If $v_p(n)=0$, the restriction $(\kappa_{\alpha_2}/\kappa_{\alpha_1})^\mu=1$
for pairs $(\alpha_2,\alpha_1)$ in \eqref{6.58} with $m_2=m_1$ 
gives no restriction on $k$, as anyway $k\leq v_p(m_2)\leq v_p(\mu)$. 
Here any $k\in\{1,...,v_p(m_2)\}$ arises. 
Thus the pairs $(\alpha_2,\alpha_1)$ in \eqref{6.58} with $m_2= m_1$
give in the case $v_p(n)=0$ the contribution 
\begin{eqnarray} 
&&\sum_{m\in\gamma_\mu^{-1}(n)\cap M:\, v_p(m)>0}\left(\sum_{k=1}^{v_p(m)-1}\varphi(m)
+\varphi(m)\cdot\frac{p-2}{p-1}\right)  \nonumber \\
&=&\sum_{m\in\gamma_\mu^{-1}(n)\cap M:\, v_p(m)>0}
(v_p(m)-\frac{1}{p-1})\cdot \varphi(m).\label{6.61}
\end{eqnarray}
If $v_p(n)=0$, the pairs $(\alpha_2,\alpha_1)$ in \eqref{6.58}
with $m_2\neq m_1$ give a contribution only if $m_2/m_1$
is a power of $p$. This contribution for all $(\alpha_2,\alpha_1)$ with
fixed $m_2$ and $m_1$ is $\varphi(m_1)$. It is the same as the contribution
of the part with $m_c=m_2$ and $m_d=m_1$ of the second factor in the
denominator of \eqref{6.46}. Thus these contributions cancel.

In the case $v_p(n)>0$ by Remark \ref{t6.3} (ii), 
any $m\in\gamma_\mu^{-1}(n)$ satisfies $v_p(m)=v_p(n)+v_p(\mu)$.
Therefore $\gamma_\mu^{-1}(n)^2\cap E_p(M)=\emptyset$, 
and because of \eqref{6.48} the second factor in the denominator
gives no contribution at all,
\begin{eqnarray}\label{6.62}
v_p\left(\prod_{m_c,m_d\in\gamma_\mu^{-1}(n)\cap M:\, m_c>m_d} 
Res(\Phi_{m_c},\Phi_{m_d})\right)=0.
\end{eqnarray}

We are left with the contributions of the first factors of the 
numerator and the denominator of \eqref{6.46} and with
\eqref{6.60} in the case $v_p(n)>0$ and with \eqref{6.61} 
in the case $v_p(n)=0$. 

Consider the case $v_p(n)>0$. Then \eqref{4.27} gives 
\begin{eqnarray}
v_p(\discr(\Phi_n)^{\psi(n)})&=& 
(v_p(n)-\frac{1}{p-1})\cdot \psi(n)\cdot\varphi(n)\label{6.63}
\end{eqnarray} 
and 
\begin{eqnarray}
&&v_p(\prod_{m\in\gamma_\mu^{-1}(n)\cap M}\discr(\Phi_m))\nonumber \\
&=& \sum_{m\in\gamma_\mu^{-1}(n)\cap M} 
(v_p(m)-\frac{1}{p-1})\cdot \varphi(m) \nonumber\\
&=& \sum_{m\in\gamma_\mu^{-1}(n)\cap M} 
(v_p(n)+v_p(\mu)-\frac{1}{p-1})\cdot \varphi(m) \nonumber\\
&=& (v_p(n)+v_p(\mu)-\frac{1}{p-1})\cdot \psi(n)\cdot \varphi(n) \nonumber\\
&=&\textup{(the contributions in \eqref{6.60} and \eqref{6.63})}.\label{6.64}
\end{eqnarray}
This shows \eqref{6.57} in the case $v_p(n)>0$.

Consider the case $v_p(n)=0$. Then \eqref{4.27} gives 
$v_p(\discr(\Phi_n)))=0$ and  
\begin{eqnarray}
&&v_p(\prod_{m\in\gamma_\mu^{-1}(n)\cap M}\discr(\Phi_m))\nonumber \\
&=& \sum_{m\in\gamma_\mu^{-1}(n)\cap M:\, v_p(m)>0} 
(v_p(m)-\frac{1}{p-1})\cdot \varphi(m) \nonumber\\
&=&\textup{(the contribution in \eqref{6.61})}.\label{6.65}
\end{eqnarray}
This shows \eqref{6.57} in the case $v_p(n)=0$.
\hfill$\Box$

\section{When does the tensor product of two Orlik blocks admit
a standard decomposition into Orlik blocks?}\label{c7}
\setcounter{equation}{0}

\noindent
Theorem \ref{t7.4} starts with two Orlik blocks
$(G,g)$ and $(H,h)$ and gives a sufficient criterion
for $(G\otimes H,g\otimes h)$ to admit a standard
decomposition into Orlik blocks.
It will be crucial for the Thom-Sebastiani sums
of singularities. The condition will work with the
sets $M\subset\N$ and $N\subset\N$
of orders of eigenvalues of $g:G_\C\to G_\C$ 
and $h:H_\C\to H_\C$.
Theorem \ref{t7.4} is preceded by some definitions
and observations. 
This section has similarities with section \ref{c6}.
Though the statement and the proof are more involved.

\begin{definition}\label{t7.1}
(a) Denote by $\mu(\C)\subset S^1$ the group of all
unit roots. Denote by $\Z[\mu(\C)]$ the group ring
with elements $\sum_{j=1}^la_j[\lambda_j]$ where
$a_j\in\Z$ and $\lambda_j\in\mu(\C)$ and with 
multiplication $[\lambda_1][\lambda_2]=[\lambda_1\lambda_2]$.
The unit element is $[1]$. The trace and the degree
of an element are 
\begin{eqnarray}\label{7.1}
\tr\left(\sum_{j=1}^la_j[\lambda_j]\right)
&:=& \sum_{j=1}^la_j\lambda_j\in\C,\\
\deg\left(\sum_{j=1}^la_j[\lambda_j]\right)
&:=& \sum_{j=1}^la_j\in\Z.\label{7.2}
\end{eqnarray}
The trace map $\tr:\Z[\mu(\C)]\to\C$ and the degree map
$\deg:\Z[\mu(\C)]\to\Z$ are ring homomorphisms.

\medskip
(b) The divisor of a unitary polynomial
$f=\prod_{j=1}^l(t-\lambda_j)$ with $\lambda_j\in\mu(\C)$
is
\begin{eqnarray}\label{7.3}
\divv(f)&:=& \sum_{j=1}^l [\lambda_j]\in\Z[\mu(\C)].
\end{eqnarray}
The divisor of an endomorphism $F:H_\C\to H_\C$ 
of a finite dimensional complex vector space $H_\C$
with characteristic polynomial $f$ is
\begin{eqnarray*}
\divv(F)&:=& \divv(f).
\end{eqnarray*}
Then $\deg(f)=\deg(\divv(f))$. 

\medskip
(c) For two unitary polynomials $f=\prod_{j=1}^l(t-\lambda_j)$
and $g=\prod_{i=1}^k(t-\kappa_i)$ with 
$\lambda_j,\kappa_i\in\mu(\C)$, define the new unitary 
polynomial $f\otimes g$ with zeros in $\mu(\C)$ by 
\begin{eqnarray}\label{7.4}
f\otimes g &:=&\prod_{j=1}^l\prod_{i=1}^k(t-\lambda_j\kappa_i).
\end{eqnarray}
Then 
\begin{eqnarray}\label{7.5}
\divv(f\otimes g)&=& \divv(f)\cdot\divv(g),\\
\tr(\divv(f\otimes g))&=& \tr(\divv(f))\cdot \tr(\divv(g)),
\label{7.6}\\
\deg(f\otimes g)&=& \deg(f)\cdot \deg(g).\label{7.7}
\end{eqnarray}

(d) For $m\in\N$ define
\begin{eqnarray}\label{7.8}
\Lambda_m&:=& \divv(t^m-1),\quad 
\Psi_m:= \divv(\Phi_m).
\end{eqnarray}
Of course, then $\Lambda_m=\sum_{d|m}\Psi_d$,
$\deg(\Lambda_m)=m$, $\deg(\Psi_m)=\varphi(m)$.

\medskip
(e) Define two maps $\beta$ and $\delta$, 
\begin{eqnarray}\label{7.9}
\beta:\N\times\N&\to&\N,
\quad \beta(m,n):=\prod_{p\in\PP:v_p(m)=v_p(n)>0}
p^{v_p(m)},
\end{eqnarray}
\begin{eqnarray}
&&\delta:\N\times\N\times\N
\to\N_0,\label{7.10}\\
&&\delta(m,n,l):=\left\{\begin{array}{l}
\varphi(\gcd(m,n))\cdot\prod_{p\in\PP:
0=v_p(c)<v_p(\beta(m,n))}\frac{p-2}{p-1} \\
\hspace*{1cm} \textup{ if }\lcm(m,n)=c\cdot l
\textup{ with }c|\beta(m,n),\\
0 \hspace*{1cm} \textup{ else.}\end{array}\right. \nonumber
\end{eqnarray}
\end{definition}

\begin{lemma}\label{t7.2}
\begin{eqnarray}\label{7.11}
\Lambda_m\cdot \Lambda_n 
&=& \gcd(m,n)\cdot\Lambda_{\lcm(m,n)}\quad
\textup{for }m,n\in\N,\\
{}[\lambda]\cdot \Lambda_m&=& \Lambda_m\quad
\textup{for }\lambda\in\mu(\C)\textup{ with }\ord(\lambda)|m,
\label{7.12}\\
\Psi_m\cdot\Psi_n&=& \sum_{l\in\N}
\delta(m,n,l)\cdot \Psi_l.\label{7.13}
\end{eqnarray}
\end{lemma}

{\bf Proof:} \eqref{7.11} and \eqref{7.12} are obvious.
\eqref{7.13} follows from the special cases
\begin{eqnarray}\label{7.14}
\Psi_m\cdot\Psi_n&=& \Psi_{m\cdot n}\quad
\textup{ if }\gcd(m,n)=1,\\
\Psi_{p^k}\cdot \Psi_{p^l}&=& \varphi(p^l)\cdot \Psi_{p^k}
\quad\textup{ for }p\in\PP\textup{ and }k>l\geq 0,
\label{7.15}\\
\Psi_{p^k}\cdot\Psi_{p^k}&=& 
\frac{p-2}{p-1}\cdot\varphi(p^k)\cdot \Psi_{p^k} + 
\varphi(p^k)\cdot \sum_{l=0}^{k-1}\Psi_{p^l}\\
&&\hspace*{2cm}\textup{for }p\in\PP\textup{ and }k>0,\nonumber
\label{7.16}
\end{eqnarray}
which follow easily from Theorem \ref{t4.5} (c).\hfill$\Box$

\begin{definition}\label{t7.3}
(a) For a prime number $p$ define the projection
\begin{eqnarray}\label{7.17}
\pi_p:\N\to\N,\quad 
m\mapsto m\cdot p^{-v_p(m)}.
\end{eqnarray}
Then $\pi_p(\N)=\{m\in\N\, |\, 
v_p(m)=0\}$.

\medskip
(b) Now fix two finite non-empty sets $M,N\subset\N$. 
Then 
\begin{eqnarray}\nonumber
\left(\sum_{m\in M}\Psi_m\right)
\cdot \left(\sum_{n\in N}\Psi_n\right)
&=& \sum_{(m,n,l)\in M\times N\times\N}
\delta(m,n,l)\cdot\Psi_l \nonumber\\
&=& \sum_{l\in \N}\chi(l)\cdot \Psi_l
= \sum_{l\in L}\chi(l)\cdot \Psi_l\label{7.18}
\end{eqnarray}
where 
\begin{eqnarray}\label{7.19}
\chi(l)&:=&\sum_{(m,n)\in M\times N}
\delta(m,n,l)\quad \textup{ for }l\in\N
\end{eqnarray}
is the multiplicity of $\Psi_l$ in the product 
$\left(\sum_{m\in M}\Psi_m\right)
\cdot \left(\sum_{n\in N}\Psi_n\right)$ and 
\begin{eqnarray}\label{7.20}
L&:=& L(M,N):=\{l\in\N\,|\, 
\chi(l)>0\}
\end{eqnarray}
is the set of numbers $l\in\N$ such that $\Psi_l$ 
turns up in this product.

\medskip
(c) Fix two finite non-empty sets $M,N\subset\N$. 
For each choice of a prime number $p$,
it will be useful to decompose $\chi(l)$ into 
pieces as follows. 
For $p\in\PP$ and $(m_0,n_0)\in\pi_p(M)\times \pi_p(N)$
and $l\in L$ define
\begin{eqnarray}\label{7.21}
\chi_{p,m_0,n_0}(l)&:=& \sum_{(m,n)\in (M\cap\pi_p^{-1}(m_0))
\times (N\cap\pi_p^{-1}(n_0))}\delta(m,n,l).
\end{eqnarray}
Then
\begin{eqnarray}\label{7.22}
\chi(l)&=& \sum_{(m_0,n_0)\in \pi_p(M)\times \pi_p(N)}
\chi_{p,m_0,n_0}(l).
\end{eqnarray}

(d) The pair $(M,N)$ of finite non-empty 
subsets of $\N$
is called {\it sdiOb-sufficient} (sdiOb for
{\it standard decomposition into Orlik blocks})
if for any prime number $p$ and any $p$-edge
$(l_a,l_b)\in E_p(L)$ 
at least one of the following two conditions holds:
\begin{eqnarray}\label{7.23}
\chi_{p,m_0,n_0}(l_b)\leq \chi_{p,m_0,n_0}(l_a)
\textup{ for any }
(m_0,n_0)\in \pi_p(M)\times \pi_p(N),\\
\chi_{p,m_0,n_0}(l_b)\geq \chi_{p,m_0,n_0}(l_a)
\textup{ for any }
(m_0,n_0)\in \pi_p(M)\times \pi_p(N)\label{7.24}
\end{eqnarray}
(these conditions will be discussed in Lemma \ref{t7.6}).
\end{definition}

\begin{theorem}\label{t7.4}
Consider two Orlik blocks $(G,g)$ and $(H,h)$.
Let $M$ and $N\subset\N$ be the finite sets
of orders of $g:G_\C\to G_\C$ respectively
$h:H_\C\to H_\C$. Then $(G\otimes H,g\otimes h)$
admits a standard decomposition into Orlik blocks
if $(M,N)$ is sdiOb-sufficient.
\end{theorem}

Theorem \ref{t7.4} will be proved in Section \ref{c8}.
Here in Section \ref{c7} we make the Remarks \ref{t7.5}, 
we make the property {\it sdiOb-sufficient}
explicit in  Lemma \ref{t7.6}, and 
we give the Examples \ref{t7.7}.

\begin{remarks}\label{t7.5}
(i) We expect that {\it if and only if} holds
in Theorem \ref{t7.4}.

\medskip
(ii) Of course,
\begin{eqnarray}
\divv(g)&=&\sum_{m\in M}\Psi_m,\quad
\divv(h)= \sum_{n\in N}\Psi_n,\label{7.25}\\
\divv(g\otimes h)&=& \divv(g)\cdot \divv(h)
=\sum_{l\in L}\chi(l)\cdot\Psi_l.\label{7.26}
\end{eqnarray}

\medskip
(iii) For any $l\in\N$, the set
of pairs $(m,n)\in\N^2$ with $\delta(m,n,l)>0$
is infinite and is as follows,
\begin{eqnarray}\label{7.27}
&&\{(m,n)\in\N^2\, |\, \delta(m,n,l)>0\}\\
&=& \{(m,n)\in\N^2\,|\, 
\textup{one has for any prime number }
p: \nonumber\\
&& \hspace*{0.5cm} \textup{either }
v_p(n)<v_p(m)=v_p(l), \nonumber\\
&& \hspace*{0.5cm}  \textup{or } v_p(m)< v_p(n)=v_p(l),
\nonumber\\
&& \hspace*{0.5cm}  \textup{or }
v_p(m)=v_p(n)\left\{\begin{array}{ll}
\geq v_p(l)&\textup{ if }p\geq 3\\
 & \textup{ or if }p=2\textup{ and }v_p(m)=0,\\
> v_p(l)&\textup{ if }p=2\textup{ and }v_p(m)>0.
\end{array}\right\} \}.\nonumber
\end{eqnarray}

(iv) In section \ref{c6}, the two conditions \eqref{6.3} and 
\eqref{6.4} for sdiOb-sufficiency of a pair $(M,\mu)$ 
were formulated only in terms of existence of $p$-edges. 
This was made explicit in Remark \ref{t6.3} (iv).
The two conditions \eqref{7.23} and \eqref{7.24} for
sdiOb-sufficiency of a pair $(M,N)$ can also be made 
explicit by necessary and sufficient conditions.
But these are more involved. We do not give all details.
We will need only the sufficient conditions in part (a)
of Lemma \ref{t7.6} and the special case where
everything is vanishing in part (b). 
\end{remarks}

\begin{lemma}\label{t7.6}
Consider the data in Remark \ref{t7.3} (d),
so two finite non-empty sets $M,N\subset\N$,
a prime number $p$
and numbers $m_0\in\pi_p(M)$, $n_0\in\pi_p(N)$,
$l_a,l_b\in L$ with $(l_a,l_b)\in E_p(L)$.
Write $k_a:=v_p(l_a)>k_b:=v_p(l_b)\geq 0$. 
Define two finite sets of exponents 
$K_{M,p,m_0}$ and $K_{N,p,n_0}\subset\N_0$ 
(they are non-empty because of $m_0\in \pi_p(M)$ and 
$n_0\in \pi_p(N)$) by
\begin{eqnarray}\label{7.28}
M\cap\pi_p^{-1}(m_0)&=& \{p^km_0\,|\,k\in K_{M,p,m_0}\},\\
N\cap\pi_p^{-1}(n_0)&=& \{p^kn_0\,|\,k\in K_{N,p,n_0}\}.\nonumber
\end{eqnarray}

(a) Then 
\begin{eqnarray}
&&k_a\in K_{M,p,m_0}-K_{N,p,n_0}\textup{ or }
k_a\in K_{N,p,n_0}-K_{M,m_0} \nonumber\\
&&\hspace*{2cm}\Rightarrow 
\chi_{p,m_0,n_0}(l_b)\leq \chi_{p,m_0,n_0}(l_a).\label{7.29}\\
&&k_a\notin K_{M,p,m_0}\cup K_{N,p,n_0}\textup{ or }
k_a\in K_{M,p,m_0}\cap K_{N,p,n_0}\nonumber\\
&&\hspace*{2cm}\Rightarrow
\chi_{p,m_0,n_0}(l_b)\geq \chi_{p,m_0,n_0}(l_a).\label{7.30}
\end{eqnarray}

(b) If $\delta(m_0,n_0,\pi_p(l_a))=0$, then
$\chi_{p,m_0,n_0}(l_a)=\chi_{p,m_0,n_0}(l_b)=0$. 
\end{lemma}

{\bf Proof:}
Write $l_0:=\pi_p(l_a)=\pi_p(l_b)$. Then
$l_a=l_0\cdot p^{k_a}$, $l_b=l_0\cdot p^{k_b}$. 
For $m=m_0\cdot p^{k_1}\in M\cap\pi_p^{-1}(m_0)$
(so $k_1\in K_{M,p,m_0}$) and 
$n=n_0\cdot p^{k_2}\in N\cap\pi_p^{-1}(n_0)$
(so $k_2\in K_{N,p,n_0}$) and 
$l=l_0\cdot p^k\in L$
\begin{eqnarray}\label{7.31}
\delta(m,n,l)&=& \delta(m_0,n_0,l_0)\cdot 
\delta(p^{k_1},p^{k_2},p^k)
\end{eqnarray}
and
\begin{eqnarray}\label{7.32}
\delta(p^{k_1},p^{k_2},p^k)= 
\left\{\begin{array}{ll}
0&\textup{ if }k_1=k_2<k\\
 &\textup{ or if }k_1\neq k_2,k\neq \max(k_1,k_2),\\
\varphi(p^{k_1})&\textup{ if }k_1=k_2>k, \\
\varphi(p^{k_1})=1&\textup{ if }k_1=k_2=k=0,\\
\varphi(p^{k_1})\cdot\frac{p-2}{p-1}&
\textup{ if }k_1=k_2=k\geq 1,\\
\varphi(p^{k_1})&\textup{ if }k_1<k_2,k=k_2,\\
\varphi(p^{k_2})&\textup{ if }k_1>k_2,k=k_1.
\end{array}\right. 
\end{eqnarray}
If $\delta(m_0,n_0,l_0)=0$, then 
$\chi_{p,m_0,n_0}(l_a)=\chi_{p,m_0,n_0}(l_b)=0$.
This proves part (b). And \eqref{7.29} and \eqref{7.30}
hold trivially in this case.

Suppose now $\delta(m_0,n_0,l_0)\neq 0$. Then for
$l=l_0\cdot p^k\in L$ as above
\begin{eqnarray}
&&\frac{\chi_{p,m_0,n_0}(l)}{\delta(m_0,n_0,l_0)}=
\sum_{(k_1,k_2)\in K_{M,p,m_0}\times K_{N,p,n_0}} 
\delta(p^{k_1},p^{k_2},p^k) \nonumber \\
&=& \delta_{(k\in K_{M,p,m_0})}\cdot 
\sum_{k_2\in K_{N,p,n_0}:\, k_2<k}\varphi(p^{k_2})
+ \delta_{(k\in K_{N,p,n_0})}\cdot 
\sum_{k_1\in K_{M,p,m_0}:\, k_1<k}\varphi(p^{k_1})\nonumber \\
&+& \delta_{(k\in K_{M,p,m_0}\cap K_{N,p,n_0}\cap\N)}\cdot 
\varphi(p^k)\cdot\frac{p-2}{p-1} 
+\delta_{(k\in K_{M,p,m_0}\cap K_{N,p,n_0}\cap\{0\})}
\nonumber \\
&+& \sum_{k_1\in K_{M,p,m_0}\cap K_{N,p,n_0}:\, k_1>k}
\varphi(p^{k_1}).\label{7.33}
\end{eqnarray}
The following notations will be used to rewrite the difference
$(\chi_{p,m_0,n_0}(l_a)-\chi_{p,m_0,n_0}(l_b))/
\delta(m_0,n_0,l_0)$
in \eqref{7.39}. 
\begin{eqnarray}
A_{M,1}&:=& \sum_{k_1\in K_{M,p,m_0}:\, k_1<k_b}
\varphi(p^{k_1}),\label{7.34}\\
A_{N,1}&:=& \sum_{k_2\in K_{N,p,n_0}:\, k_2<k_b}
\varphi(p^{k_2}),\label{7.35}\\
A_{M,2}&:=& \sum_{k_1\in K_{M,p,m_0}:\, k_b<k_1<k_a}
\varphi(p^{k_1}),\label{7.36}\\
A_{N,2}&:=& \sum_{k_2\in K_{N,p,n_0}:\, k_b<k_2<k_a}
\varphi(p^{k_2}),\label{7.37}\\
A_3&:=& \sum_{k_1\in K_{M,p,m_0}\cap K_{N,p,n_0}:\, k_b<k_1<k_a}
\varphi(p^{k_1}).\label{7.38}
\end{eqnarray}
\begin{eqnarray}
&&\frac{\chi_{p,m_0,n_0}(l_a)-\chi_{p,m_0,n_0}(l_b)}
{\delta(m_0,n_0,l_0)}\nonumber\\
&=& (\delta_{(k_a\in K_{M,p,m_0})}-\delta_{(k_b\in K_{M,p,m_0})})
\cdot A_{N,1} \label{7.39}\\
&+& (\delta_{(k_a\in K_{N,p,n_0})}-\delta_{(k_b\in K_{N,p,n_0})})
\cdot A_{M,1} \nonumber \\
&+& (\delta_{(k_a\in K_{M,p,m_0}, k_b\in K_{N,p,n_0})}
+\delta_{(k_a\in K_{N,p,n_0}, k_b\in K_{M,p,m_0})})\cdot
\varphi(p^{k_b}) \nonumber\\
&-& \delta_{(k_b\in K_{M,p,m_0}\cap K_{N,p,n_0}\cap\N)}
\cdot (p-2)p^{k_b-1}   
-\delta_{(k_b\in K_{M,p,m_0}\cap K_{N,p,n_0}\cap\{0\})}\nonumber\\
&+& \delta_{(k_a\in K_{M,p,m_0})}\cdot A_{N,2}
+ \delta_{(k_a\in K_{N,p,n_0})}\cdot A_{M,2}
- A_3 \nonumber\\
&-& \delta_{(k_a\in K_{M,p,m_0}\cap K_{N,p,n_0})}
\cdot p^{k_a-1}.\nonumber
\end{eqnarray}
Observe
\begin{eqnarray}\label{7.40}
\left\{ 
\begin{array}{l}A_{M,1}\\ A_{N,1}\end{array}
\right\} 
\left\{ \begin{array}{ll} \leq p^{k_b-1}&\textup{if }k_b>0\\
=0&\textup{if }k_b=0\end{array}\right. ,\\ 
A_3\leq \left\{
\begin{array}{l}A_{M,2}\\ A_{N,2}\end{array}
\right\}\leq p^{k_a-1}-p^{k_b}.\label{7.41}
\end{eqnarray}

Now we prove \eqref{7.29}.
Suppose $k_a\in K_{M,p,m_0}-K_{N,p,n_0}$. Then \eqref{7.39} is
\begin{eqnarray}\nonumber
&&(1-\delta_{(k_b\in K_{M,p,m_0})})\cdot A_{N,1}
+\delta_{(k_b\in K_{N,p,n_0}\cap \N)}
\cdot (-A_{M,1}+(p-1)p^{k_b-1})\\
&-&\delta_{(k_b\in K_{M,p,m_0}\cap K_{N,p,n_0}\cap \N)}\cdot
(p-2)p^{k_b-1} \nonumber\\
&+&(\delta_{(k_b\in K_{N,p,n_0}\cap\{0\})}-
\delta_{(k_b\in K_{M,p,m_0}\cap K_{N,p,n_0}\cap\{0\})}) 
+ (A_{N,2}-A_3)\label{7.42}\\
&\geq & 0\cdot A_{N,1}+ 
(\delta_{(k_b\in K_{N,p,n_0}\cap \N)}
-\delta_{(k_b\in K_{M,p,m_0}\cap K_{N,p,n_0}\cap \N)})\cdot
(p-2)p^{k_b-1}\nonumber \\
&& + (0)+(0)\nonumber\\
&\geq& 0. \nonumber
\end{eqnarray}
The case $k_a\in K_{N,p,n_0}-K_{M,p,m_0}$ is treated
analogously. This shows \eqref{7.29}.

Now we prove \eqref{7.30}. 
If $k_a\notin K_{M,p,m_0}\cup K_{N,p,n_0}$, then \eqref{7.39} is
\begin{eqnarray}
&-&\delta_{(k_b\in K_{M,p,m_0})}\cdot A_{N,1}
-\delta_{(k_b\in K_{N,p,n_0})}\cdot A_{M,1} \nonumber \\
&-& \delta_{(k_b\in K_{M,p,m_0}\cap K_{N,p,n_0}\cap \N)}\cdot (p-2)p^{k_b-1} \nonumber\\
&-&\delta_{(k_b\in K_{M,p,m_0}\cap K_{N,p,n_0}\cap\{0\})}
-A_3 \label{7.43}\\
&\leq& 0.\nonumber
\end{eqnarray} 
Suppose $k_a\in K_{M,p,m_0}\cap K_{N,p,n_0}$. 
Then \eqref{7.39} is
\begin{eqnarray}\nonumber
&&(1-\delta_{(k_b\in K_{M,p,m_0})})\cdot A_{N,1}
+(1-\delta_{(k_b\in K_{N,p,n_0})})\cdot A_{M,1}\\
&+&\delta_{(k_b\in (K_{M,p,m_0}\cup K_{N,p,n_0})\cap \N)}
\cdot (p-1)p^{k_b-1} 
+\delta_{(k_b\in K_{M,p,m_0}\cap K_{N,p,n_0}\cap \N)}\cdot
p^{k_b-1} \nonumber\\
&+&\delta_{(k_b\in (K_{M,p,m_0}\cup K_{N,p,n_0})\cap \{0\})} 
+ (A_{N,2}+A_{M,2}-A_3)-p^{k_a-1}\label{7.44}\\
&\leq & p^{k_b} + 
\Bigl(\sum_{k_1\in K_{M,p,m_0}\cup K_{N,p,n_0}:\, 
k_b<k_1<k_a}\varphi(p^{k_1})
\Bigr)  -p^{k_a-1} \nonumber\\
&\leq& p^{k_b}+(p^{k_a-1}-p^{k_b})-p^{k_a-1}\nonumber\\
&\leq& 0.\nonumber
\end{eqnarray}
Part (a) is proved.
\hfill$\Box$

\begin{examples}\label{t7.7}
(i) The Milnor lattice with monodromy $(H_{Mil},h_{Mil})$
of an $A_\mu$ singularity $x_1^{\mu+1}$ in one variable 
is a single Orlik block with set
\begin{eqnarray}\label{7.45}
M&=&\{m\in\N\,|\, m|(\mu+1)\}-\{1\}
\end{eqnarray}
of orders of eigenvalues of the monodromy.
This is well known. It also follows from Theorem \ref{t1.3}
(a) and from the fact that all eigenvalues have
multiplicity 1 and the set of their orders is $M$.
For any prime number and any $m_0\in \pi_p(M)$
\begin{eqnarray}\label{7.46}
K_{M,p,m_0}&=&\left\{\begin{array}{ll}
\Z_{[0,v_p(\mu+1)]} & \textup{ if }m_0\neq 1,\\
\Z_{[1,v_p(\mu+1)]} & \textup{ if }m_0=1
\end{array}\right. 
\end{eqnarray}
(where $[1,0]=\emptyset$ and $\Z_{[1,0]}=\emptyset$).

\medskip
(ii) We consider an $A_\mu$-singularity $x_1^{\mu+1}$
as in (i) and an $A_\nu$-singularity $x_2^{\nu+1}$ 
with set $N$ of orders of eigenvalues of its monodromy.
$N$ and $K_{N,p,n_0}$ for $n_0\in\pi_p(N)$ are
as in \eqref{7.45} and \eqref{7.46}, with 
$\mu$ replaced by $\nu$.
We will show with Theorem \ref{t7.4}
that the Thom-Sebastiani sum 
$A_{\mu}\otimes A_{\nu}$, i.e. $x_1^{\mu+1}+x_2^{\nu+1}$, 
satisfies Orlik's conjecture,
i.e. its Milnor lattice with monodromy 
admits a standard decomposition into Orlik blocks.

It is the tensor product of the Milnor lattices
with monodromies of the two $A$-type singularities
\cite{ST71}, so the tensor product of Orlik blocks
with sets $M$ and $N$. In order to apply Theorem \ref{t7.4},
we have to show that \eqref{7.23} or \eqref{7.24} holds
for any prime number $p$ and any $p$-edge 
$(l_a,l_b)\in E_p(L)$. 

Then $k_a:=v_p(l_a)>k_b:=v_p(l_b)\geq 0$. 
The shape \eqref{7.46} of $K_{M,p,m_0}$
and analogously for $K_{N,p,n_0}$ shows that
the properties $(k_a\in K_{M,p,m_0}\textup{ or not})$
and $(k_a\in K_{N,p,n_0}\textup{ or not}$)
are independent of the choice of $m_0\in \pi_p(M)$
and $n_0\in\pi_p(N)$. Therefore the hypotheses
in  \eqref{7.29} and \eqref{7.30}
are independent of the choice of $m_0\in \pi_p(M)$
and $n_0\in\pi_p(N)$.
This shows that \eqref{7.23} or \eqref{7.24}
holds for any prime number $p$ and any 
$p$-edge $(l_a,l_b)\in E_p(L)$. 
Therefore Theorem \ref{t7.4} applies. 
The Thom-Sebastiani sum $x_1^{\mu+1}+x_2^{\nu+1}$
satisfies Orlik's conjecture.

This example is a very special case of Theorem \ref{t1.3}
(d). But we find it instructive to see the sdiOb-sufficiency
condition at work in a simple case. 

\medskip
(iii) The following is a small abstract example 
of a pair $(M,N)$ which is not sdiOb-sufficient. 
Consider $M:=\{3\}$ and $N:=\{2,3\}$. 
We list some relevant data:
\begin{eqnarray*}
L= \{1,3,6\},\ E_2(L)=\{(6,3)\},\ E_3(L)=\{(3,1)\},\\
\pi_2(M)=\{3\},\ \pi_3(M)=\{1\},\ 
\pi_2(N)=\{1,3\}, \pi_3(N)=\{2,1\},\\
\begin{array}{l|l|l|l|l|l|l}
l\in L & \delta(3,2,l) & \delta(3,3,l) 
& \chi_{2,3,1}(l) & \chi_{2,3,3}(l) 
& \chi_{3,1,2}(l) & \chi_{3,1,1}(l) \\ \hline
1 & 0 & 2 & 0 & 2 & 0 & 2\\
3 & 0 & 1 & 0 & 1 & 0 & 1\\
6 & 1 & 0 & 1 & 0 & 1 & 0
\end{array}
\end{eqnarray*}
The $3$-edge $(3,1)$ satisfies \eqref{7.24}.
But the $2$-edge $(6,3)$ satisfies neither
\eqref{7.23} nor \eqref{7.24}.
The pair $(M,N)$ is not sdiOb-sufficient.
Theorem \ref{t7.4} does not apply.
Remark \ref{t7.5} (i) even claims that
the tensor product $H^{[\Phi_3]}\otimes H^{[\Phi_2\Phi_3]}$ 
of the two Orlik blocks $H^{[\Phi_3]}$ and $H^{[\Phi_2\Phi_3]}$
does not admit a standard decomposition into
Orlik blocks. This is true. It can be proved with
Theorem \ref{t5.1} and Example \ref{t7.7} (ii) as follows. 
\begin{eqnarray*}
H^{[\Phi_3]}\otimes H^{[\Phi_2\Phi_3]}
&\cong& H^{[\Phi_3]}\otimes (H^{[\Phi_2]}\oplus H^{[\Phi_3]})\\
&& (\textup{apply Theorem \ref{t5.1} (a) to } 
H^{[\Phi_2\Phi_3]})\\
&\cong& H^{[\Phi_3]}\otimes H^{[\Phi_2]}
\oplus H^{[\Phi_3]}\otimes H^{[\Phi_3]})\\
&\cong & H^{[\Phi_6]} \oplus (H^{[\Phi_3\Phi_1]}\oplus
H^{[\Phi_1]})\\
&& (\textup{apply Example \ref{t7.7} (ii) to both summands})\\
&\cong& H^{[\Phi_6\Phi_1]}\oplus H^{[\Phi_3\Phi_1]}\\
&& (\textup{apply Theorem \ref{t5.1} (a) to }
H^{[\Phi_6\Phi_1]})\\
&\not\cong& H^{[\Phi_6\Phi_3\Phi_1]}\oplus H^{[\Phi_1]}\\
&&(\textup{apply Theorem \ref{t5.1} (b) with}\\
&& f_1=\Phi_6, f_2=1, f_3=\Phi_1, f_4=\Phi_3).
\end{eqnarray*}
So, $H^{[\Phi_3]}\otimes H^{[\Phi_2\Phi_3]}$
is isomorphic to the non-standard decomposition
$H^{[\Phi_6\Phi_1]}\oplus H^{[\Phi_3\Phi_1]}$ into 
Orlik blocks, but not to the standard decomposition
$H^{[\Phi_6\Phi_3\Phi_1]}\oplus H^{[\Phi_1]}$
into Orlik blocks.

\medskip
(iv) Here we present a small extension of the example
in (iii) which is an sdiOb-sufficient pair 
$(M,N)$, namely $M=\{3\}$ and $N=\{1,2,3\}$. 
So Theorem \ref{t7.4} applies and shows that
$H^{[\Phi_3]}\otimes H^{[\Phi_1\Phi_2\Phi_3]}$ 
admits a standard decomposition into Orlik blocks.
We list some relevant data:
\begin{eqnarray*}
L= \{1,3,6\},\ E_2(L)=\{(6,3)\},\ E_3(L)=\{(3,1)\},\\
\pi_2(M)=\{3\},\ \pi_3(M)=\{1\},\ 
\pi_2(N)=\{1,3\}, \pi_3(N)=\{2,1\},\\
\begin{array}{l|l|l|l|l|l|l|l}
l & \delta(3,1,l) & \delta(3,2,l) & \delta(3,3,l) 
& \chi_{2,3,1}(l) & \chi_{2,3,3}(l) 
& \chi_{3,1,2}(l) & \chi_{3,1,1}(l) \\ \hline
1 & 0 & 0 & 2 & 0 & 2 & 0 & 2\\
3 & 1 & 0 & 1 & 1 & 1 & 0 & 2\\
6 & 0 & 1 & 0 & 1 & 0 & 1 & 0
\end{array}
\end{eqnarray*}
The 3-edge $(3,1)$ satisfies \eqref{7.24} and \eqref{7.23}.
The 2-edge $(6,3)$ satisfies \eqref{7.24}.
\end{examples}

\section{Proof of Theorem \ref{t7.4}}\label{c8}
\setcounter{equation}{0}

\noindent
This section is devoted to the proof of Theorem \ref{t7.4}.
The proof is similar to the one of Theorem \ref{t6.2},
especially the beginning. But the later part is 
much more involved.
We will use all the notations of section \ref{c7}.

Consider the pair  $(G\otimes H, g\otimes h)$. 
Formula \eqref{7.26} for $\divv(g\otimes h)$ implies 
that the characteristic polynomial of $g\otimes h$ is 
\begin{eqnarray}\label{8.1}
p_{G\otimes H,g\otimes h}&=& 
\prod_{l\in L}\Phi_l^{\chi(l)}.
\end{eqnarray}
Let $p_1, ..., p_{\chi_0}\in\C[t]$ be the unique
unitary polynomials with 
$p_{G\otimes H,g\otimes h}=p_1\cdot ...\cdot p_{\chi_0}$
and $p_{\chi_0}|p_{\chi_0-1}|...|p_2|p_1$
and $p_1$ the minimal polynomial of $g\otimes h$. 
Here 
\begin{eqnarray}\label{8.2}
\chi_0=\max_{l\in L}\chi(l)\quad\textup{and}\quad  
p_k=\prod_{l\in L:\, \chi(l)\geq k}\Phi_l\in\Z[t].
\end{eqnarray}
If $G\otimes H$ admits a standard decomposition into
Orlik blocks, that is isomorphic to 
$\bigoplus_{k=1}^{\chi_0}H^{[p_k]}$. 

We want to apply Theorem \ref{t3.1} to $(G\otimes H,g\otimes h)$
instead of $(H,h)$. 

Let $a_0$ and $b_0$ be generating elements of the Orlik
blocks $(G,g)$ and $(H,h)$, so 
$G=\bigoplus_{i=0}^{\rk G-1}\Z\cdot g^i(a_0)$ and 
$H=\bigoplus_{j=0}^{\rk H-1}\Z\cdot h^j(b_0)$. Define
\begin{eqnarray}\label{8.3}
a_i&:=& g^i(a_0) \quad\textup{for }i\geq 0,\\
b_j&:=& h^j(b_0) \quad\textup{for }j\geq 0,\label{8.4}\\
{\bf a}&:=& (a_0,a_1,...,a_{\rk G-1})\quad\textup{a }
\Z\textup{-basis of }G,\label{8.5}\\
{\bf b}&:=& (b_0,b_1,...,b_{\rk H-1})\quad\textup{a }
\Z\textup{-basis of }H.\label{8.6}
\end{eqnarray}
Then the tuple 
\begin{eqnarray}\label{8.7}
\CC^{st}&:=& {\bf a}\otimes {\bf b}\nonumber\\
&:=&(a_0\otimes b_0,a_0\otimes b_1,...,
a_0\otimes b_{\rk H-1},\\
&& a_1\otimes b_0,a_1\otimes b_1,...,
a_1\otimes b_{\rk H-1},\nonumber\\
&& \hspace*{2cm} \vdots \nonumber\\
&& a_{\rk G-1}\otimes b_0,a_{\rk G-1}\otimes b_1,...,
a_{\rk G-1}\otimes b_{\rk H-1}),\nonumber
\end{eqnarray} 
is a $\Z$-basis of $G\otimes H$. 

Observe $(g\otimes h)^k(a_i\otimes b_j)=
a_{i+k}\otimes b_{j+k}$. Consider the tuples
\begin{eqnarray}\label{8.8}
\CC^{dec,i}&:=& (a_i\otimes b_0,a_{i+1}\otimes b_1,...,
a_{i+\deg p_{i+1}-1}\otimes b_{\deg p_{i+1}-1})\\
&&
\textup{for }i\in\Z_{[0,\chi_0-1]},\nonumber \\
\CC^{dec}&:=& 
(\CC^{dec,0},\CC^{dec,1},...,\CC^{dec,\chi_0-1}).\label{8.9}
\end{eqnarray} 
We will show that $\CC^{dec}$ 
is a $\Z$-basis of $G\otimes H$ if and only if
$(M,N)$ is sdiOb-sufficient.
This and Theorem \ref{t3.1} show that $(G\otimes H,g\otimes h)$
admits a standard decomposition into Orlik blocks
if $(M,N)$ is sdiOb-sufficient. 

It remains to show that the matrix 
$M(\CC^{st},\CC^{dec})$ with
$\CC^{dec}=\CC^{st}\cdot M(\CC^{st},\CC^{dec})$ 
has determinant $\pm 1$ if and only if 
the pair $(M,N)$ is sdiOb-sufficient.
We will calculate this determinant up to the sign.
We will show
\begin{eqnarray}\label{8.10}
\det M(\CC^{st},\CC^{dec})
&=& (\pm 1)\cdot \prod_{p\in\PP}
\prod_{(l_a,l_b)\in E_p(L)}p^{\varphi(l_b)\cdot 
\Xi(l_a,l_b)},
\end{eqnarray}
where for $(l_a,l_b)\in E_p(L)$ 
\begin{eqnarray}
\Xi_1(l_a,l_b)&:=& \min(\chi(l_a),\chi(l_b))\in\N_0,
\nonumber \\
\Xi_{2,p}(l_a,l_b)&:=& \sum_{(m_0,n_0)\in\pi_p(M)\times \pi_p(N)}
\min(\chi_{p,m_0,n_0}(l_a),\chi_{p,m_0,n_0}(l_b))
\in\N_0, \nonumber \\
\Xi(l_a,l_b)&:=& \Xi_1(l_a,l_b)-\Xi_{2,p}(l_a,l_b)\in\N_0.
\label{8.11} 
\end{eqnarray}
Obviously $\Xi(l_a,l_b)=0$ if and only if 
\eqref{7.23} or \eqref{7.24} holds. Therefore
$\det M(\CC^{st},\CC^{dec})=\pm 1$ if and only if 
the pair $(M,N)$ is sdiOb-sufficient.

As in the proof of theorem \ref{t6.2},
for the calculation of the determinant,
we will consider also certain tuples of eigenvectors
of $g$, $h$ and $g\otimes h$. 
Let $\{\kappa_1,\kappa_2,...,\kappa_{\rk G}\}$
be the set of eigenvalues of $g$, and let 
$\{\lambda_1,\lambda_2,...,\lambda_{\rk H}\}$
be the set of eigenvalues of $h$, in both sets 
the indices are chosen such that
\begin{eqnarray}\label{8.12}
\ord(\kappa_\alpha)\leq \ord(\kappa_\beta)
\textup{ and }\ord(\lambda_\alpha)\leq \ord(\lambda_\beta)
\textup{ if }\alpha<\beta.
\end{eqnarray}
Then $a_0$ determines a basis of eigenvectors
${\bf u}=(u_1,...,u_{\rk G})$ of $G_\C$,
and $b_0$ determines a basis of eigenvectors 
${\bf v}=(v_1,...,v_{\rk H})$ of $H_\C$, with 
\begin{eqnarray}\label{8.13}
a_0&=& \sum_{\alpha=1}^{\rk G}u_\alpha,\quad
g(u_\alpha)= \kappa_\alpha\cdot u_\alpha\\
b_0&=&\sum_{\beta=1}^{\rk H}v_\beta,\quad
h(v_\beta)=\lambda_\beta\cdot v_\beta.\label{8.14}
\end{eqnarray}
The base change matrices 
\begin{eqnarray}\label{8.15}
M({\bf u},{\bf a})&=& \begin{pmatrix}
1&\kappa_1^1&\cdots & \kappa_1^{\rk G-1}\\
\vdots & \vdots & & \vdots \\
1&\kappa_{\rk G}^1&\cdots & \kappa_{\rk G}^{\rk G-1}
\end{pmatrix}\\
\textup{and }
M({\bf v},{\bf b})&=& \begin{pmatrix}
1&\lambda_1^1&\cdots & \lambda_1^{\rk H-1}\\
\vdots & \vdots & & \vdots \\
1&\lambda_{\rk H}^1&\cdots & \lambda_{\rk H}^{\rk H-1}
\end{pmatrix}
\end{eqnarray}\label{8.16}
are Vandermonde matrices. The tuple
\begin{eqnarray}\label{8.17}
\CC^I&:=& {\bf u}\otimes {\bf v}\nonumber\\
&:=&(u_1\otimes v_1,u_1\otimes v_2,...,
u_1\otimes v_{\rk H},\\
&& u_2\otimes v_1,u_2\otimes v_2,...,
u_2\otimes v_{\rk H},\nonumber\\
&& \hspace*{2cm} \vdots \nonumber\\
&& u_{\rk G}\otimes v_1,u_{\rk G}\otimes v_2,...,
u_{\rk G}\otimes v_{\rk H}),\nonumber
\end{eqnarray} 
is a $\C$-basis of $G_\C\otimes H_\C$.
And the base change matrix with $\CC^{st}$ is 
\begin{eqnarray}\label{8.18}
M(\CC^I,\CC^{st})&=& M({\bf u},{\bf a})
\otimes M({\bf v},{\bf b}).
\end{eqnarray}

Let $\{\mu_1,\mu_2,...,\mu_{\deg p_1}\}$ 
be the set of eigenvalues of $g\otimes h$.
For $\gamma\in\Z_{[1,\deg p_1]}$ define
\begin{eqnarray}\label{8.19}
C(\gamma)&:=& \{(\alpha,\beta)\in\Z_{[1,\rk G]}
\times \Z_{[1,\rk H]}\,|\, 
\kappa_\alpha\lambda_\beta=\mu_\gamma\}\\
&=:& \bigcup_{k=1}^{|C(\gamma)|}
\{(\alpha(\gamma,k),\beta(\gamma,k)\}.\nonumber
\end{eqnarray}
The eigenvalues $\mu_1,...,\mu_{\deg p_1}$ of $g\otimes h$ 
are indexed such that
\begin{eqnarray}\label{8.20}
\chi_0=|C(1)|\geq |C(2)|\geq ...\geq |C(\deg p_1)|\geq 1.
\end{eqnarray}
For any $\gamma\in \Z_{[1,\deg p_1]}$, 
the space $\bigoplus_{k=1}^{|C(\gamma)|}\C\cdot
u_{\alpha(\gamma,k)}\otimes v_{\beta(\gamma,k)}\subset 
G_\C\otimes H_\C$
is the eigenspace with eigenvalue $\mu_\gamma$
of $g\otimes h$. For any $\gamma\in\Z_{[1,\deg p_1]}$
and any $i\geq 0$, the vector
\begin{eqnarray}\label{8.21}
w_{i,\gamma}&:=& \sum_{k=1}^{|C(\gamma)|}
\kappa_{\alpha(\gamma,k)}^i\cdot 
u_{\alpha(\gamma,k)}\otimes v_{\beta(\gamma,k)}
\end{eqnarray}
is an eigenvector with eigenvalue $\mu_\gamma$ of $g\otimes h$.
It is useful as for $j\geq 0$ 
\begin{eqnarray}
a_{i+j}\otimes b_j &=& 
\sum_{\alpha=1}^{\rk G}\sum_{\beta=1}^{\rk H}
\kappa_\alpha^{i+j}\lambda_\beta^j\cdot u_\alpha\otimes
v_\beta \nonumber \\
&=& \sum_{\gamma=1}^{\deg p_1}\mu_\gamma^j\cdot w_{i,\gamma}.
\label{8.22}
\end{eqnarray}
Consider for $\gamma\in\Z_{[1,\deg p_1]}$ and
$i\in\Z_{[0,\chi_0-1]}$ 
the following tuples of eigenvectors of $g\otimes h$, 
\begin{eqnarray}
\CC^{II,\gamma}&:=& 
(u_{\alpha(\gamma,1)}\otimes v_{\beta(\gamma,1)},...,
u_{\alpha(\gamma,|C(\gamma)|)}\otimes 
v_{\beta(\gamma,|C(\gamma)|)}),\label{8.23}\\
\CC^{II}&:=& (\CC^{II,1},...,\CC^{II,\deg p_1}),\label{8.24}\\
\CC^{III,\gamma}&:=& 
(w_{0,\gamma},...,w_{|C(\gamma)|-1,\gamma}),\label{8.25}\\
\CC^{III}&:=& (\CC^{III,1},...,\CC^{III,\deg p_1}),
\label{8.26}\\
\CC^{V,i}&:=& (w_{i,1},w_{i,2},...,w_{i,\deg p_1}),
\label{8.27}\\ 
\CC^{IV,i}&:=& (w_{i,1},w_{i,2},...,w_{i,\deg p_{i+1}}),
\label{8.28}\\ 
\CC^{IV}&:=& (\CC^{IV,1},...,\CC^{IV,\chi_0}).\label{8.29}
\end{eqnarray}
$\CC^{II,\gamma}$ and $\CC^{III,\gamma}$ are 
$\C$-bases of the eigenspace in $G_\C\otimes H_\C$ 
with eigenvalue $\mu_\gamma$ of $g\otimes h$. 
The base change matrix 
$M(\CC^{II,\gamma},\CC^{III,\gamma})$ rewrites the 
relation \eqref{8.21}. It is the Vandermonde matrix
\begin{eqnarray}\label{8.30}
M(\CC^{II,\gamma},\CC^{III,\gamma})&=& 
\begin{pmatrix} 1&\kappa_{\alpha(\gamma,1)}^1 & ... & 
\kappa_{\alpha(\gamma,1)}^{|C(\gamma)|-1}\\
\vdots & \vdots & & \vdots \\
1&\kappa_{\alpha(\gamma,|C(\gamma)|)}^1& \cdots & 
\kappa_{\alpha(\gamma,|C(\gamma)|)}^{|C(\gamma)|-1}
\end{pmatrix}.
\end{eqnarray}
$\CC^I$, $\CC^{II}$, $\CC^{III}$ and $\CC^{IV}$
are $\C$-bases of $G_\C\otimes H_\C$. 
The base change matrices $M(\CC^I,\CC^{II})$
and $M(\CC^{III},\CC^{IV})$ are just permutation matrices.
The entries of $\CC^{V,i}$ are linearly independent.
Therefore the following rectangular 
$(\deg p_1)\times(\deg p_{i+1})$-matrix is well defined.
It rewrites the relations \eqref{8.22}.
\begin{eqnarray}\label{8.31}
M(\CC^{V,i},\CC^{dec,i})&=& 
\begin{pmatrix} 1& \mu_1^1& \cdots& \mu_1^{\deg p_{i+1}-1}\\
\vdots & \vdots & & \vdots \\
1& \mu_{\deg p_1}^1& \cdots & \mu_{\deg p_1}^{\deg p_{i+1}-1}
\end{pmatrix}.
\end{eqnarray}
Now we want to describe the base change matrix
$M(\CC^{IV},\CC^{dec})$. 
Observe that $w_{i,\gamma}$ is in the case
$\gamma>\deg p_{i+1}$ a linear combination of the
entries of $\CC^{III,\gamma}$. 
Therefore the matrix $M(\CC^{IV},\CC^{dec})$
is a block upper triangular matrix whose diagonal blocks
are obtained from the matrices in \eqref{8.31}
by cutting off the lower lines. The diagonal blocks
are the Vandermonde matrices
\begin{eqnarray}\label{8.32}
M^{[i]}:= 
\begin{pmatrix} 1& \mu_1^1& \cdots& \mu_1^{\deg p_{i+1}-1}\\
\vdots & \vdots & & \vdots \\
1& \mu_{\deg p_{i+1}}^1& \cdots & 
\mu_{\deg p_{i+1}}^{\deg p_{i+1}-1}\end{pmatrix}
\quad\textup{for }i\in\Z_{[0,\chi_0-1]}.
\end{eqnarray}
Now the matrix $M(\CC^{st},\CC^{dec})$ can be written 
as the product of matrices
\begin{eqnarray}\label{8.33}
M(\CC^{st},\CC^I)M(\CC^I,\CC^{II})M(\CC^{II},\CC^{III})
M(\CC^{III},\CC^{IV})M(\CC^{IV},\CC^{dec}).
\end{eqnarray}
The absolute value of its determinant is
\begin{eqnarray}\nonumber
&&|\det M(\CC^{st},\CC^{dec})|\\
&=& |\frac{\det M(C^{II},\C^{III})\cdot 
\det M(C^{IV},\CC^{dec})}
{\det M(\CC^I,\CC^{st})}|\nonumber\\
&=& |\frac{\prod_{\gamma=1}^{\deg p_1}
\det M(\CC^{II,\gamma},\CC^{III,\gamma})\cdot
\prod_{i=0}^{\chi_0-1}M^{[i]}}
{(M({\bf u},{\bf a})^{\rk H}\cdot M({\bf v},{\bf b})^{\rk G}}|.
\label{8.34}
\end{eqnarray}
Here we used the equality 
$\det(A\otimes B)=(\det A)^b\cdot (\det B)^a$ for
square matrices $A\in M_{a\times a}(\C)$ and
$B\in M_{b\times b}(\C)$ 
(which becomes obvious if one looks at Jordan normal forms).

The quotient in \eqref{8.34} is a quotient of determinants of
Vandermonde matrices. The determinants are
\begin{eqnarray}\label{8.35}
\det M({\bf u},{\bf a})&=& 
\prod_{1\leq \alpha_1<\alpha_2\leq \rk G}
(\kappa_{\alpha_2}-\kappa_{\alpha_1}),\\
\det M({\bf v},{\bf b})&=& 
\prod_{1\leq \beta_1<\beta_2\leq \rk H}
(\lambda_{\beta_2}-\lambda_{\beta_1}),\label{8.36}\\
\det M(\CC^{II,\gamma},\CC^{III,\gamma})
&=& \prod_{1\leq k_1<k_2\leq |C(\gamma)|}
(\kappa_{\alpha(\gamma,k_2)}-\kappa_{\alpha(\gamma,k_1)}),
\label{8.37}\\
\det M^{[i]}
&=& \prod_{1\leq k_1<k_2\leq \deg p_{i+1}}
(\mu_{k_2}-\mu_{k_1}).\label{8.38}
\end{eqnarray}

Only now we use that $\kappa_\alpha$, $\lambda_\beta$ 
and $\mu_\gamma$ are unit roots and how they are related. 
For $\det M({\bf u},{\bf a})$ we obtain with the formulas
\eqref{4.3} and \eqref{4.9}
\begin{eqnarray}
|\det M({\bf u},{\bf a})|&=& 
\prod_{m_c,m_d\in M:\, m_c>m_d}|\Res(\Phi_{m_c},\Phi_{m_d})|
\nonumber\\
&&\cdot \prod_{m\in M}\sqrt{\discr{\Phi_m}}.\label{8.39}
\end{eqnarray}
Now choose a prime number $p$.
The exponent of $p$ in 
$|\det M({\bf u},{\bf a})|^2 \in\N$
is because of \eqref{4.3} and \eqref{4.9} and 
Theorem \ref{t4.5} (d) and (e)
\begin{eqnarray}\label{8.40}
v_p(|\det M({\bf u},{\bf a})|^2)
&=&2\sum_{(m_c,m_d)\in E_p(M)}\varphi(m_d)\\
&& + \sum_{m\in M:\, v_p(m)\geq 1}
\varphi(m)\cdot (v_p(m)-\frac{1}{p-1}).\nonumber
\end{eqnarray}
Analogously, we obtain for 
$\det M({\bf v},{\bf b})$ and 
$\prod_{i=0}^{\chi_0-1} \det M^{[i]}$
\begin{eqnarray}\label{8.41}
v_p(|\det M({\bf v},{\bf b})|^2)
&=&2\sum_{(n_c,n_d)\in E_p(N)}\varphi(n_d)\\
&&+ \sum_{n\in N:\, v_p(n)\geq 1}
\varphi(n)\cdot (v_p(n)-\frac{1}{p-1}),\nonumber\\
v_p(\left(\prod_{i=0}^{\chi_0-1}|\det M^{[i]}|\right)^2)
&=& 2 \sum_{(l_a,l_b)\in E_p(L)} \varphi(l_b)
\cdot\min(\chi(l_b),\chi(l_a))\hspace*{1cm} \label{8.42}\\
&&+ \sum_{l\in L:\, v_p(l)\geq 1}\varphi(l)\cdot
(v_p(l)-\frac{1}{p-1})\cdot \chi(l). \nonumber
\end{eqnarray}
As the squares of the absolute values of the 
other factors in the quotient in 
\eqref{8.34} are positive integers, 
the square of the absolute value of
the factor $\prod_{\gamma=1}^{\deg p_1}
\det M(\CC^{II,\gamma},\CC^{III,\gamma})$
is a positive rational number, and the value
$v_p(|\prod_{\gamma=1}^{\deg p_1}
\det M(\CC^{II,\gamma},\CC^{III,\gamma})|^2)$ is
a well defined integer. It is the most difficult part.
We will discuss it below.

In order to prove \eqref{8.10}, we have to show for 
any prime number $p$
\begin{eqnarray}\nonumber
&&-rk H\cdot\left(2\sum_{(m_c,m_d)\in E_p(M)}\varphi(m_d)
+ \sum_{m\in M:\, v_p(m)\geq 1}
\varphi(m)\cdot (v_p(m)-\frac{1}{p-1})\right)\\
&&-\rk G\cdot\left(2\sum_{(n_c,n_d)\in E_p(N)}\varphi(n_d)
+ \sum_{n\in N:\, v_p(n)\geq 1}
\varphi(n)\cdot (v_p(n)-\frac{1}{p-1})\right)\nonumber \\
&&+2\sum_{(l_a,l_b)\in E_p(L)}\varphi(l_b)\cdot
\Xi_{2,p}(l_a,l_b)
+ \sum_{l\in L:\, v_p(l)\geq 1}
\varphi(l)\cdot (v_p(l)-\frac{1}{p-1})\cdot\chi(l)\nonumber\\
&&+v_p(\Bigl|\prod_{\gamma=1}^{\deg p_1}
\det M(\CC^{II,\gamma},\CC^{III,\gamma})\Bigr|^2)
=0.\label{8.43}
\end{eqnarray}

Now we discuss the most difficult part.
Consider a difference 
$\kappa_{\alpha(\gamma,k_2)}-\kappa_{\alpha(\gamma,k_1)}$ 
in \eqref{8.37}. Denote for a moment 
\begin{eqnarray*}
m_1:=\ord(\kappa_{\alpha(\gamma,k_1)}),&&
m_2:=\ord(\kappa_{\alpha(\gamma,k_2)}),\\
n_1:=\ord(\lambda_{\beta(\gamma,k_1)}),&& 
n_2:=\ord(\lambda_{\beta(\gamma,k_2)}),\\
l:=\ord(\mu_\gamma),&&
\nu:=\frac{\kappa_{\alpha(\gamma,k_1)}}
{\kappa_{\alpha(\gamma,k_2)}}.
\end{eqnarray*}
Theorem \ref{t4.5} (a) and (b) give 
\begin{eqnarray}
\norm_{\ord\nu} (1-\nu)
= \left\{\begin{array}{ll}
\pm q &\textup{ if }\ord\nu\textup{ is a (positive)}\\ 
&\textup{ power of a prime number }q\\
\pm 1&\textup{ else.}\end{array}\right. \label{8.44}
\end{eqnarray}
Therefore the difference
$\kappa_{\alpha(\gamma,k_2)}-\kappa_{\alpha(\gamma,k_1)}
=\kappa_{\alpha(\gamma,k_2)}\cdot (1-\nu)$ 
makes a contribution to 
$v_p(|\prod_{\gamma=1}^{\deg p_1}
\det M(\CC^{II,\gamma},\CC^{III,\gamma})|^2)$
only if $\ord\nu$ is a (positive) power of $p$. 
By Theorem \ref{t4.5}  (c)(i), this holds only
if $\pi_p(m_1)=\pi_p(m_2)$. 
Because of 
$\kappa_{\alpha(\gamma,k)}\lambda_{\beta(\gamma,k)}
=\mu_\gamma$ for any $k\in\Z_{[1,|C(\gamma)|]}$, we have 
\begin{eqnarray}\label{8.45}
\frac{\kappa_{\alpha(\gamma,k_1)}}
{\kappa_{\alpha(\gamma,k_2)}}
= \frac{\lambda_{\beta(\gamma,k_2)}}
{\lambda_{\beta(\gamma,k_1)}}.
\end{eqnarray}
Therefore also $\pi_p(n_1)=\pi_p(n_2)$ is needed. 
The following lemma makes the sizes
of the possible contributions precise.

\begin{lemma}\label{t8.1}
Fix a prime number $p$. 
Fix $m_1,m_2\in M$, $n_1,n_2\in N$ and $l\in L$
with $\pi_p(m_1)=\pi_p(m_2)$ and $\pi_p(n_1)=\pi_p(n_2)$
and $\max(v_p(m_1),v_p(m_2))>0$. 
The contribution of all differences
$\kappa_{\alpha(\gamma,k_2)}-\kappa_{\alpha(\gamma,k_1)}$ 
in \eqref{8.37} with 
\begin{eqnarray}\label{8.46}
\ord\gamma=l,\quad \ord\kappa_{\alpha(\gamma,k_i)}=m_i,\quad
\ord\lambda_{\beta(\gamma,k_i)}=n_i
\end{eqnarray}
to $v_p(|\prod_{\gamma=1}^{\deg p_1}
\det M(\CC^{II,\gamma},\CC^{III,\gamma})|^2)$
is as follows. We can suppose $m_1\leq m_2$
(by exchanging them if necessary).
If $m_1=m_2$ we can suppose $n_1\leq n_2$.  
We have the following ten cases.
The contribution has always the shape
$\varphi(l)\cdot\delta(m_1,n_1,l)\cdot(\textup{a factor }F)$.
\begin{eqnarray*}
\begin{array}{ll}
\textup{Case} & \textup{the factor }F \\ \hline
(C1)\ v_p(m_1)<v_p(m_2)=v_p(n_2)>v_p(n_1): & 2 \\
(C2a)\ v_p(m_1)<v_p(m_2)<v_p(n_1)=v_p(n_2)=v_p(l): & 
2\\
(C2b)\ v_p(n_1)<v_p(n_2)<v_p(m_1)=v_p(m_2)=v_p(l): &
2\\
(C3)\ v_p(m_1)<v_p(m_2)=v_p(n_1)=v_p(l)>v_p(n_2): &
2\cdot\frac{\varphi(p^{v_p(n_2)})}
{\varphi(p^{v_p(n_1)})} \\
(C4a)\ v_p(m_1)<v_p(m_2)=v_p(n_1)=v_p(n_2)=v_p(l): &
2\cdot\frac{p-2}{p-1} \\
(C4b)\ v_p(n_1)<v_p(n_2)=v_p(m_1)=v_p(m_2)=v_p(l): &
2\cdot\frac{p-2}{p-1} \\
(C5a)\ 1\leq v_p(m_1)=v_p(m_2)<v_p(n_1)=v_p(n_2)=v_p(l): & 
v_p(m_1)-\frac{1}{p-1} \\
(C5b)\ 1\leq v_p(n_1)=v_p(n_2)<v_p(m_1)=v_p(m_2)=v_p(l): & 
v_p(n_1)-\frac{1}{p-1} \\
(C6)\ v_p(m_1)=v_p(m_2)=v_p(n_1)=v_p(n_2)>v_p(l): &
v_p(m_1)-\frac{1}{p-1} \\
(C7)\  v_p(m_1)=v_p(m_2)=v_p(n_1)=v_p(n_2)=v_p(l)>0: &
v_p(m_1)-\frac{2}{p-1}
\end{array}
\end{eqnarray*}
(In the case $(C1)$  
$v_p(l)\leq\max(v_p(m_1),v_p(n_1))$ with equality if 
$v_p(m_1)\neq v_p(n_1)$.)
\end{lemma}

{\bf Proof:}
$\kappa_{\alpha(\gamma,k_i)}\lambda_{\beta(\gamma,k_i)}=
\mu_\gamma$ implies $\delta(m_i,n_i,l)>0$ for $i\in\{1,2\}$.
With Remark \ref{t7.5} (iii), which describes
the set $\{(m,n)\in \N^2\,|\, \delta(m,n,l)>0\}$,
one obtains easily that only the ten cases in the 
lemma are possible (assuming $m_1\leq m_2$, and
assuming $n_1\leq n_2$ in the case $m_1=m_2$).

Write for $m\in\N$
\begin{eqnarray}\label{8.47}
\Z_m:=\Z_{[0,m-1]}\quad\textup{and}\quad
\Z_m^*:=\{a\in\Z_m\,|\, \gcd(a,m)=1\}.
\end{eqnarray}
The contribution of all differences 
$\kappa_{\alpha(\gamma,k_2)}-\kappa_{\alpha(\gamma,k_1)}$ 
in \eqref{8.37} with \eqref{8.46}
to $v_p(|\prod_{\gamma=1}^{\deg p_1}
\det M(\CC^{II,\gamma},\CC^{III,\gamma})|^2)$
is 
\begin{eqnarray}\nonumber
2\sum_{k\geq 1}\frac{1}{\varphi(p^k)}
\Bigl|\{(a_1,b_1,a_2,b_2,c_1,c_2)\hspace*{4cm}\\
\in\Z_{m_1}^*\times
\Z_{n_1}^*\times \Z_{m_2}^*\times \Z_{n_2}^*\times
\Z_l^*\times\Z_{p^k}^*\,| \nonumber \\
\Bigl(a_1<a_2\textup{ if }m_1=m_2\textup{ and }
n_1=n_2\Bigr),\nonumber\\
\frac{a_1}{m_1}+\frac{b_1}{n_1}\equiv\frac{c_1}{l}
\equiv\frac{a_2}{m_2}+\frac{b_2}{n_2}\modd\Z \nonumber\\
\frac{a_1}{m_1}-\frac{a_2}{m_2}\equiv\frac{c_2}{p^k}
\equiv\frac{b_2}{n_2}-\frac{b_1}{n_1}\modd\Z\}\Bigr|.
\label{8.48}
\end{eqnarray}
This follows with the identifications
\begin{eqnarray*}
\kappa_{\alpha(\gamma,k_j)}=e^{2\pi i \frac{a_j}{m_j}},\ 
\lambda_{\beta(\gamma,k_j)}=e^{2\pi i \frac{b_j}{m_j}},\ 
\gamma=e^{2\pi i \frac{c_1}{l}},\ 
\frac{\kappa_{\alpha(\gamma,k_1)}}{\kappa_{\alpha(\gamma,k_2)}}
=e^{2\pi i \frac{c_2}{p^k}}
\end{eqnarray*}
from the fact that the norm in \eqref{8.44}
has $\varphi(\ord\nu)=\varphi(p^k)$ factors
which together give a factor $p$.
In the case $(m_1=m_2\textup{ and }n_1=n_2)$ the condition
$k_1<k_2$ in \eqref{8.37} leads to $a_1<a_2$
(or just as well $a_1>a_2$) in \eqref{8.48}.
In the other cases, the condition $k_1<k_2$
is taken care of by $m_1<m_2$ or 
$(m_1=m_2\textup{ and }n_1<n_2)$. 

It remains to calculate the number in \eqref{8.48}
in each of the ten cases. 
For an arbitrary fixed $c_1\in \Z_l^*$, the set
\begin{eqnarray}\label{8.49}
\{(a_1,b_1)\in\Z_{m_1}^*\times \Z_{n_1}^*\, |\,
\frac{a_1}{m_1}+\frac{b_1}{n_1}\equiv\frac{c_1}{l}\modd\Z\}
\end{eqnarray}
has $\delta(m_1,n_1,l)$ elements. This follows from
\eqref{7.13} $(\Psi_m\cdot \Psi_n=\sum_{\tilde l\geq 1}
\delta(m,n,\tilde l)\Psi_{\tilde l}$).
Therefore the set
\begin{eqnarray}\label{8.50}
\{(a_1,b_1,c_1)\in\Z_{m_1}^*\times \Z_{n_1}^*\times \Z_l^*\, |\,
\frac{a_1}{m_1}+\frac{b_1}{n_1}\equiv\frac{c_1}{l}\modd\Z\}
\end{eqnarray}
has $\varphi(l)\cdot\delta(m_1,n_1,l)$ elements.

Now fix for a moment such numbers $a_1,b_1$ and $c_1$.
If for some $k\geq 1$ and some $c_2\in\Z_{p^k}^*$
elements $a_2\in\Z_{m_2}$ and $b_2\in\Z_{n_2}$
with 
\begin{eqnarray}\label{8.51}
\frac{a_1}{m_1}-\frac{a_2}{m_2}\equiv \frac{c_2}{p^k}
\equiv \frac{b_2}{n_2}-\frac{b_1}{n_1}\modd\Z
\end{eqnarray}
exist, they are uniquely determined by these equations,
respectively by
\begin{eqnarray}\label{8.52}
a_2&\equiv& a_1\frac{m_2}{m_1}-c_2\frac{m_2}{p^k}
\mod m_2\Z,\\
b_2&\equiv& b_1\frac{n_2}{n_1}+c_2\frac{n_2}{p^k}\mod n_2\Z.
\label{8.53}
\end{eqnarray}
The questions are whether they exist, and if yes,
whether $a_2$ is in $\Z_{m_2}^*$ and 
$b_2$ is in $\Z_{n_2}^*$. 

Because of $m_2=m_1\cdot p^{v_p(m_2)-v_p(m_1)}$, 
$a_2$ is in $\Z_{m_2}$ if and only if $k\leq v_p(m_2)$.
In the cases with $m_1<m_2$, $a_2$ is in $\Z_{m_2}^*$
if and only if $k=v_p(m_2)$. 
In the cases with $m_1=m_2$ and $k<v_p(m_2)$,
$a_2$ is in $\Z_{m_2}^*$. In the cases with $m_1=m_2$
and $k=v_p(m_2)$, we need 
$a_1\not\equiv c_2\pi_p(m_1)\modd p\Z$.
The set of triples $(a_1,b_1,c_1)$ in \eqref{8.50}
which satisfy this, has 
$\varphi(l)\cdot\delta(m_1,n_1,l)\cdot\frac{p-2}{p-1}$
elements. 

Similarly, $n_2=n_1\cdot p^{v_p(n_2)-v_p(n_1)}$,
but here also $n_1>n_2$ is possible.  
If $n_1\leq n_2$, then 
$b_2$ is in $\Z_{n_2}$ if and only if $k\leq v_p(n_2)$.
In the cases with $n_1<n_2$, $b_2$ is in $\Z_{n_2}^*$
if and only if $k=v_p(n_2)$. 
In the cases with $n_1=n_2$ and $k<v_p(n_2)$,
$b_2$ is in $\Z_{n_2}^*$. In the cases with $n_1=n_2$
and $k=v_p(n_2)$, we need 
$b_1\not\equiv -c_2\pi_p(n_1)\modd p\Z$.
The set of triples $(a_1,b_1,c_1)$ in \eqref{8.50}
which satisfy this, has 
$\varphi(l)\cdot\delta(m_1,n_1,l)\cdot\frac{p-2}{p-1}$
elements. 

Now consider the cases with $n_1>n_2$. Then $b_2\in\Z_{n_2}^*$
if and only if $k=v_p(n_1)$ and 
$v_p(b_1+c_2\pi_p(n_2))=k-v_p(n_2)$. We claim that the
set of triples $(a_1,b_1,c_1)$ in \eqref{8.50}
which satisfy this, has $\varphi(l)\cdot\delta(m_1,n_1,l)\cdot
\frac{\varphi(p^{v_p(n_2)})}{\varphi(p^k)}$ elements. 
The claim is a consequence of the following facts:
\begin{list}{}{}
\item[(1)] The projection from the set in \eqref{8.50} to
$\Z_{p^k}^*$ which sends a triple $(a_1,b_1,c_1)$ to the
class of $b_1$ in $\Z_{p^k}^*$ is surjective and all fibers
have the same size. 
\item[(2)] For a fixed $\gamma\in\Z_{p^k}^*$, the set
$\{\beta\in\Z_{p^k}^*\, |\, v_p(\beta+\gamma)=k-v_p(n_2)\}$
has $\varphi(p^{v_p(n_2)})$ elements.
\end{list}
Fact (2) follows from 
$\Lambda(p^k,p^k,p,v_p(n_2))=\varphi(p^k)$
in Theorem \ref{4.5} (c)(iii) and from the factor
$\varphi(p^{v_p(n_2)})^{-1}$ in the definition \eqref{4.19}
of $\Lambda(p^k,p^k,p,v_p(n_2))$.

Now the ten cases in Lemma \ref{t8.1}
can be obtained by combining the results of the discussion
above, with some extra arguments for $(C6)$ and $(C7)$. 
We discuss the cases separately.

$(C1)$: $v_p(m_1)<v_p(m_2)$ gives the unique $k=v_p(m_2)$ and
$a_2\in\Z_{m_2}^*$. Then automatically $k=v_p(n_2)>v_p(n_1)$,
and thus also $b_2\in\Z_{n_2}^*$. 
As $c_2\in\Z_{p^k}^*$ was fixed arbitrarily,
the number in \eqref{8.48} is 2 times the number of elements
of the set in \eqref{8.50}, so the factor $F$ is $2$.

$(C2a)$: $v_p(m_1)<v_p(m_2)$ gives the unique $k=v_p(m_2)$ and
$a_2\in\Z_{m_2}^*$. Then automatically $k<v_p(n_1)=v_p(n_2)$,
and thus also $b_2\in\Z_{n_2}^*$. 
As $c_2\in\Z_{p^k}^*$ was fixed arbitrarily,
the number in \eqref{8.48} is 2 times the number of elements
of the set in \eqref{8.50}, so the factor $F$ is $2$.

$(C2b)$: Analogous to $(C2a)$. 

$(C3)$: $v_p(m_1)<v_p(m_2)$ gives the unique $k=v_p(m_2)$ and
$a_2\in\Z_{m_2}^*$. Then automatically $k=v_p(n_1)>v_p(n_2)$.
For $b_2\in\Z_{n_2}^*$ we need the set of triples 
$(a_1,b_1,c_1)$ which satisfy 
$v_p(b_1+c_2\pi_p(n_2))=k-v_p(n_2)$.
Their number is $\varphi(l)\cdot\delta(m_1,n_1,l)\cdot
\frac{\varphi(p^{v_p(n_2)})}{\varphi(p^k)}$. 
As $c_2\in\Z_{p^k}^*$ was fixed arbitrarily,
the number in \eqref{8.48} is 2 times this number, 
so the factor $F$ is 
$2\cdot\frac{\varphi(p^{v_p(n_2)})}{\varphi(p^k)}$.

$(C4a)$: $v_p(m_1)<v_p(m_2)$ gives the unique $k=v_p(m_2)$ and
$a_2\in\Z_{m_2}^*$. Then automatically $k=v_p(n_1)=v_p(n_2)$.
For $b_2\in\Z_{n_2}^*$ we need the set of triples 
$(a_1,b_1,c_1)$ which satisfy 
$b_1\not\equiv -c_2\pi_p(n_1)\modd p\Z$. 
Their number is $\varphi(l)\cdot\delta(m_1,n_1,l)\cdot
\frac{p-2}{p-1}$. 
As $c_2\in\Z_{p^k}^*$ was fixed arbitrarily,
the number in \eqref{8.48} is 2 times this number, 
so the factor $F$ is $2\cdot\frac{p-2}{p-1}$.

$(C4b)$: Analogous to $(C4a)$. 

$(C5a)$: $v_p(m_1)=v_p(m_2)$ gives $k\in\Z_{[1,v_p(m_1)]}$. 
For any $k\in\Z_{[1,v_p(m_1)-1]}$, $a_2\in\Z_{m_2}^*$
and $b_2\in\Z_{n_2}^*$ hold automatically.
For $k=v_p(m_2)<v_p(n_2)$,  $b_2\in\Z_{n_2}^*$ holds
automatically, but for $a_2\in \Z_{m_2}^*$ we need the set
of triples $(a_1,b_1,c_1)$ in \eqref{8.50} with 
$a_1\not\equiv c_2\pi_p(m_1)\modd p\Z$.
Their number is $\varphi(l)\cdot\delta(m_1,n_1,l)\cdot
\frac{p-2}{p-1}$. 
As $m_1=m_2$ and $n_2=n_1$, the condition $a_1<a_2$ 
in \eqref{8.48} cancels the factor 2 in \eqref{8.48}. 
As $c_2\in\Z_{p^k}^*$ was fixed arbitrarily,
the number in \eqref{8.48} is 
$\varphi(l)\cdot\delta(m_1,n_1,l)\cdot
((v_p(m_1)-1)\cdot 1 + \frac{p-2}{p-1})$, so the factor
$F$ is $v_p(m_1)-\frac{1}{p-1}$. 

$(C5b)$: Analogous to $(C5a)$. 

$(C6)$ and $(C7)$: $v_p(m_1)=v_p(m_2)=v_p(n_1)=v_p(n_2)$ gives 
$k\in\Z_{[1,v_p(m_1)]}$. 
For any $k\in\Z_{[1,v_p(m_1)-1]}$, $a_2\in\Z_{m_2}^*$
and $b_2\in\Z_{n_2}^*$ hold automatically.
As in the case $(C5a)$, these $k$ give the contribution
$v_p(m_1)-1$ to the factor $F$. 
Now consider $k=v_p(m_1)$. Then we need for 
$a_2\in \Z_{m_2}^*$ 
and $b_2\in\Z_{n_2}^*$ the set of triples $(a_1,b_1,c_1)$ 
in \eqref{8.50} with 
$a_1\not\equiv c_2\pi_p(m_1)\modd p\Z$
and $b_1\not\equiv -c_2\pi_p(n_1)\modd p\Z$. 
The relation in \eqref{8.50} implies (for $k=v_p(m_2)$)
\begin{eqnarray}\label{8.54}
(a_1\pi_p(n_1)+b_1\pi_p(m_1))\pi_p(l)\equiv
c_1\pi_p(m_1)\pi_p(n_1)\cdot p^{k-v_p(l)}\modd p^k\Z.
\end{eqnarray}
Therefore the condition 
$b_1\not\equiv -c_2\pi_p(n_1)\modd p\Z$
is equivalent to the condition 
\begin{eqnarray}\label{8.55}
a_1\pi_p(l)\not\equiv
c_1\pi_p(m_1)p^{k-v_p(l)}+c_2\pi_p(m_1)\pi_p(l) \modd p\Z.
\end{eqnarray}
Of course, the condition 
$a_1\not\equiv c_2\pi_p(m_1)\modd p\Z$
is equivalent to the condition
\begin{eqnarray}\label{8.56}
a_1\pi_p(l)\not\equiv c_2\pi_p(m_1)\pi_p(l) \modd p\Z.
\end{eqnarray}
In the case $(C6)$ we have $k>v_p(l)$, so then the conditions
\eqref{8.55} and \eqref{8.56} coincide. 
Then $k=v_p(m_1)$ gives the same contribution
$\frac{p-2}{p-1}$ to the factor $F$ as in the case
$(C5a)$. Then the factor $F$ is as in the case $(C5a)$.
In the case $(C7)$ we have $k=v_p(l)$,
so then the conditions \eqref{8.55} and \eqref{8.56}
exclude for the class of $a_1$ in $\Z_p^*$ two different
numbers. Then the set of triples in \eqref{8.50} which
satisfy \eqref{8.55} and \eqref{8.56} is
$\varphi(l)\cdot \delta(m_1,n_1,l)\cdot\frac{p-3}{p-1}$.
Then $k=v_p(m_1)$ gives the contribution 
$\frac{p-3}{p-1}$ to the factor $F$.
Then the factor $F$ is $(v_p(m_1)-1)\cdot 1+\frac{p-3}{p-1}
=v_p(m_1)-\frac{2}{p-1}$.
\hfill $\Box$ 

\bigskip
The following notations will be useful.
\begin{eqnarray}\label{8.57}
T&:=& \{(m,n,l)\in M\times N\times L\,|\, \delta(m,n,l)>0\},\\
\pi_{p,3}&:& \N^3\to \N^3,\ (m,n,l)\mapsto
(\pi_p(m),\pi_p(n),\pi_p(l)),\label{8.58}\\
T_0&:=& \pi_{p,3}(T),\label{8.59}\\
K_{{\bf t}_0}&:=& \{(k_1,k_2,k)\in 
K_{M,p,m_0}\times K_{N,p,n_0}\times K_{L,p,l_0}\,|\, 
\delta(p^{k_1},p^{k_2},p^k)>0\} \nonumber\\
&&\textup{for }{\bf t}_0=(m_0,n_0,l_0)\in T_0.\label{8.60}
\end{eqnarray}
Observe that for any 
${\bf t}=(m,n,l)\in\N^3$
\begin{eqnarray}\label{8.61}
\delta({\bf t})=
\delta(\pi_{p,3}({\bf t}))\cdot
\delta(p^{v_p(m)},p^{v_p(n)},p^{v_p(l)}).
\end{eqnarray}
Thus for ${\bf t}=(m,n,l)\in T$ and 
${\bf t}_0=\pi_{p,3}({\bf t})\in T_0$
we have $\delta({\bf t}_0)>0$ and 
$(v_p(m),v_p(n),v_p(l))\in K_{{\bf t}_0}$.

In the next lemma, we split each summand in
\eqref{8.43} into a sum of pieces over $T$.
Afterwards we will show that the sum of all pieces
over each fiber of the map $\id\times\id\times \pi_p:
T\to M\times N\times \pi_p(L)$ 
is already equal to 0. That implies \eqref{8.43}.

\begin{lemma}\label{t8.2}
The seven summands in \eqref{8.43} can be written
as sums of pieces over $T$ as follows.
\begin{eqnarray*}
-\rk H\cdot 2\sum_{(m_c,m_d)\in E_p(M)}\varphi(m_d)
&=& \sum_{{\bf t}\in T}\varphi(l)\delta({\bf t})\cdot 
A^{(1)}_{{\bf t}},\\
-\rk H\cdot\sum_{m\in M:\, v_p(m)\geq 1}\varphi(m)
(v_p(m)-\frac{1}{p-1})
&=& \sum_{{\bf t}\in T}\varphi(l)\delta({\bf t})\cdot 
A^{(2)}_{{\bf t}},\\
-\rk G\cdot 2\sum_{(n_c,n_d)\in E_p(N)}\varphi(n_d)
&=& \sum_{{\bf t}\in T}\varphi(l)\delta({\bf t})\cdot 
A^{(3)}_{{\bf t}},\\
-\rk G\cdot \sum_{n\in N:\, v_p(n)\geq 1}\varphi(n)
(v_p(n)-\frac{1}{p-1})
&=& \sum_{{\bf t}\in T}\varphi(l)\delta({\bf t})\cdot 
A^{(4)}_{{\bf t}},\\
2\sum_{(l_a,l_b)\in E_p(L)}\varphi(l_b)\cdot
\Xi_{2,p}(l_a,l_b)
&=& \sum_{{\bf t}\in T}\varphi(l)\delta({\bf t})\cdot 
A^{(5)}_{{\bf t}},\\
\sum_{l\in L:\, v_p(l)\geq 1}\varphi(l)
(v_p(l)-\frac{1}{p-1})\cdot \chi(l)
&=& \sum_{{\bf t}\in T}\varphi(l)\delta({\bf t})\cdot 
A^{(6)}_{{\bf t}},\\
v_p(\Bigl|\prod_{\gamma=1}^{\deg p_1}\det 
M(\CC^{II,\gamma},\CC^{III,\gamma})\Bigr|^2)
&=& \sum_{{\bf t}\in T}\varphi(l)\delta({\bf t})\cdot 
A^{(7)}_{{\bf t}},
\end{eqnarray*}
where the summands $A^{(j)}_{{\bf t}}$ are as follows.
Here $(m,n,l):={\bf t}$, 
$(m_0,n_0,l_0):=\pi_{p,3}({\bf t})$ and 
$(k_1,k_2,k):=(v_p(m),v_p(n),v_p(l))$.
\begin{eqnarray}
A^{(1)}_{{\bf t}}&=& -2\cdot |K_{M,p,m_0}\cap\Z_{>k_1}|,
\label{8.62}\\
A^{(2)}_{{\bf t}}&=& -\delta_{(k_1>0)}\cdot(k_1-\frac{1}{p-1}),
\label{8.63}\hspace*{1cm}\\
A^{(3)}_{{\bf t}}&=& -2\cdot |K_{N,p,n_0}\cap\Z_{>k_2}|,
\label{8.64}\\
A^{(4)}_{{\bf t}}&=& -\delta_{(k_2>0)}\cdot(k_2-\frac{1}{p-1}),
\label{8.65}\\
A^{(6)}_{{\bf t}}&=& \delta_{(k>0)}\cdot(k-\frac{1}{p-1}),
\label{8.66}
\end{eqnarray}
and
\begin{eqnarray}\label{8.67}
A^{(5)}_{{\bf t}}&=& A^{(5,1)}_{{\bf t}}+A^{(5,2)}_{{\bf t}}
\quad\textup{with}\\ 
A^{(5,1)}_{{\bf t}}&=& 
2|(K_{M,p,m_0}\cup K_{N,p,n_0}-(K_{M,p,m_0}\cap K_{N,p,n_0}))
\cap \Z_{>\max(k_1,k_2)}|, \hspace*{1cm} \label{8.68}\\
A^{(5,2)}_{{\bf t}}&=& 
\delta_{(k>0,\, \max(k_1,k_2)\in K_{M,p,m_0}\cap K_{N,p,n_0})}
\cdot \frac{2}{p-1},\label{8.69}
\end{eqnarray}
and
\begin{eqnarray}\label{8.70}
A^{(7)}_{{\bf t}}&=& 
\sum_{\tiny \begin{array}{r}
(C)\in\{(C1),(C2a),(C2b),(C3),(C4a),\\
(C4b),(C5a),(C5b),(C6),(C7)\}\end{array}} A^{(C)}_{{\bf t}}
\end{eqnarray}
with $A^{(C)}_{{\bf t}}$ as follows.
\begin{eqnarray}
A^{(C1)}_{{\bf t}}&=& 
2\cdot |K_{M,p,m_0}\cap K_{N,p,n_0}\cap\Z_{>\max(k_1,k_2)}|,
\label{8.71} \\
A^{(C2a)}_{{\bf t}}&=& 
2\cdot |K_{M,p,m_0}\cap\Z_{[k_1+1,k_2-1]}|,\label{8.72} \\
A^{(C2b)}_{{\bf t}}&=& 
2\cdot |K_{N,p,n_0}\cap\Z_{[k_2+1,k_1-1]}|,\label{8.73} \\
A^{(C3)}_{{\bf t}}&=& A^{(C1)}_{{\bf t}},\label{8.74}
\end{eqnarray}
\begin{eqnarray}
A^{(C4a)}_{{\bf t}}&=& 
\delta_{(k_2>k_1,\,k_2\in K_{M,p,m_0})}
\cdot 2\cdot\frac{p-2}{p-1},\label{8.75}\\
A^{(C4b)}_{{\bf t}}&=& 
\delta_{(k_1>k_2,\,k_1\in K_{N,p,n_0})}
\cdot 2\cdot\frac{p-2}{p-1},
\label{8.76}
\end{eqnarray}
\begin{eqnarray}
A^{(C5a)}_{{\bf t}}&=& 
\delta_{(1\leq k_1<k_2)}\cdot(k_1-\frac{1}{p-1}),\label{8.77}\\
A^{(C5b)}_{{\bf t}}&=& 
\delta_{(1\leq k_2<k_1)}\cdot(k_2-\frac{1}{p-1}),\label{8.78}\\
A^{(C6)}_{{\bf t}}&=& 
\delta_{(k_1=k_2>k)}\cdot(k_1-\frac{1}{p-1}),\label{8.79}\\
A^{(C7)}_{{\bf t}}&=& 
\delta_{(k_1=k_2=k\geq 1)}\cdot(k_1-\frac{2}{p-1}).\label{8.80}
\end{eqnarray}
\end{lemma}

{\bf Proof:}
First \eqref{8.62}--\eqref{8.66} are proved.
Observe 
\begin{eqnarray}\label{8.81}
\rk G&=& \sum_{m\in M}\varphi(m),\quad
\rk H = \sum_{n\in N}\varphi(n),\\
\varphi(m)\varphi(n)&=& \sum_{l\in L}\varphi(l)\delta(m,n,l)
\quad\textup{for }(m,n)\in M\times N,
\label{8.82}
\end{eqnarray}
where \eqref{8.82} follows from \eqref{7.13}. 
Now one finds
\begin{eqnarray*}
&&\rk H\cdot \sum_{(m\cdot p^a,m)\in E_p(M)}
\varphi(m)\\
&=& \Bigl(\sum_{n\in N}\varphi(n)\Bigr)
\Bigl(\sum_{m\in M}\varphi(m)\cdot 
|K_{M,p,\pi_p(m)}\cap \Z_{>v_p(m)}|\Bigr)\\
&=& \sum_{(m,n,l)\in T}\varphi(l)\delta(m,n,l)
\cdot |K_{M,p,\pi_p(m)}\cap \Z_{>v_p(m)}|.
\end{eqnarray*}
This shows \eqref{8.62}.
Similar calculations show \eqref{8.63}--\eqref{8.65}.
With \eqref{7.19},
\begin{eqnarray*}
\chi(l)&=& \sum_{(m,n)\in M\times N}\delta(m,n,l),
\end{eqnarray*}
one sees also \eqref{8.66} immediately.

Now we will prove \eqref{8.67}--\eqref{8.69}. 
It is more difficult.
The second equality sign below uses Lemma \ref{t7.6} (a).
There we sum over ${\bf t}_0=(m_0,n_0,l_0)\in T_0$. 
After the third equality sign below, we split on purpose
$\Z_{>k}=\Z_{[k+1,\max(k_1,k_2)]}\cup\Z_{>\max(k_1,k_2)}$,
as that will allow some simplification later.
Note that for $(k_1,k_2,k)\in K_{{\bf t}_0}$
the inequality $\max(k_1,k_2)>k$ implies 
$k_1=k_2$ and $\delta(p^{k_1},p^{k_2},p^k)=\varphi(p^{k_1})$.
\begin{eqnarray*}
&&\sum_{(l_a,l_b)\in E_p(L)}\varphi(l_b)\cdot
\Xi_{2,p}(l_a,l_b)\\
&=& \sum_{{\bf t}_0\in T_0}
\sum_{k,k_3\in K_{L,p,l_0}:\, k_3>k} \varphi(l_0p^k)
\min(\chi_{p,m_0,n_0}(l_0p^{k_3}),
\chi_{p,m_0,n_0}(l_0p^{k}))\\
&=& \sum_{{\bf t}_0\in T_0}
\sum_{k\in K_{L,p,l_0}} \varphi(l_0p^k)\cdot\Bigl(
\chi_{p,m_0,n_0}(l_0p^k)\\
&& \hspace*{3cm}\cdot |((K_{M,p,m_0}\cup K_{N,p,n_0})
-(K_{M,p,m_0}\cap K_{N,p,n_0}))\cap \Z_{>k}|\\
&&\hspace*{2cm} + 
\sum_{k_3\in K_{L,p,l_0}-(K_{M,p,m_0}\cup K_{N,p,n_0}):\, k_3>k}\chi_{p,m_0,n_0}(l_0p^{k_3})\\
&&\hspace*{2cm} + 
\sum_{k_3\in K_{M,p,m_0}\cap K_{N,p,n_0}:\, k_3>k}\chi_{p,m_0,n_0}(l_0p^{k_3}) \Bigr)
\end{eqnarray*}
\begin{eqnarray*}
&=& \sum_{{\bf t}_0\in T_0}\varphi(l_0)\delta({\bf t}_0)
\cdot \sum_{k\in K_{L,p,l_0}} \varphi(p^k)\cdot\Bigl(
\nonumber\\
&&\sum_{(k_1,k_2)\in K_{M,p,m_0}\times K_{N,p,n_0}}\Bigl[ \delta(p^{k_1},p^{k_2},p^k)  \\
&&\hspace*{2cm} \cdot |((K_{M,p,m_0}\cup K_{N,p,n_0})
-(K_{M,p,m_0}\cap K_{N,p,n_0}))\cap \Z_{>\max(k_1,k_2)}|
\nonumber\\
&& \hspace*{1cm} +\ \delta_{(k_1=k_2>k)}\cdot
 \varphi(p^{k_1}) \nonumber\\
&& \hspace*{2cm}\cdot |((K_{M,p,m_0}\cup K_{N,p,n_0})
-(K_{M,p,m_0}\cap K_{N,p,n_0}))\cap \Z_{[k+1,k_1]}|
\Bigr] 
\end{eqnarray*}
\begin{eqnarray}
&+&   
\sum_{k_3\in K_{L,p,l_0}-(K_{M,p,m_0}\cup K_{N,p,n_0}):\, k_3>k}
\sum_{k_4\in K_{M,p,m_0}\cap K_{N,p,n_0}:\, k_4>k_3}
\varphi(p^{k_4}) \nonumber\\
&+&
\sum_{k_3\in K_{M,p,m_0}\cap K_{N,p,n_0}:\, k_3>k}
\Bigl[ \sum_{k_4\in K_{M,p,m_0}\cap K_{N,p,n_0}:\, k_4>k_3}
\varphi(p^{k_4}) +\varphi(p^{k_3})\frac{p-2}{p-1} \nonumber\\
&&+\sum_{k_1\in K_{M,p,m_0}:\, k_1<k_3}\varphi(p^{k_1})
+\sum_{k_2\in K_{N,p,n_0}:\, k_2<k_3}\varphi(p^{k_2}) 
\Bigr] \Bigr)\label{8.83}
\end{eqnarray}
The first three lines in this formula \eqref{8.83} give 
$\sum_{{\bf t}\in T}\varphi(l)\delta({\bf t})\cdot
\frac{1}{2}A^{(5,1)}_{\bf t}$.
We will show that the last five lines 
give $\sum_{{\bf t}\in T}\varphi(l)\delta({\bf t})\cdot
\frac{1}{2}A^{(5,2)}_{\bf t}$.

First we take care of line eight of \eqref{8.83}.
It is
\begin{eqnarray}
&&\sum_{{\bf t}_0\in T_0}\varphi(l_0)\delta({\bf t}_0)
\cdot \sum_{k\in K_{L,p,l_0}} \varphi(p^k)\cdot\Bigl( \nonumber\\
&&\sum_{k_3\in K_{M,p,m_0}\cap K_{N,p,n_0}:\, k_3>k}\Bigl[
\sum_{k_1\in K_{M,p,m_0}:\, k_1<k_3}\varphi(p^{k_1})
+\sum_{k_2\in K_{N,p,n_0}:\, k_2<k_3}\varphi(p^{k_2}) 
\Bigr] \Bigr)\nonumber\\
&=&\sum_{{\bf t}_0\in T_0}\varphi(l_0)\delta({\bf t}_0)
\cdot\sum_{k_3\in K_{M,p,m_0}\cap K_{N,p,n_0}:\, k_3>0}
\sum_{k=0}^{k_3-1}\varphi(p^k) \nonumber\\
&&\hspace*{0.5cm}\cdot \ \Bigl[
\sum_{k_1\in K_{M,p,m_0}:\, k_1<k_3}\varphi(p^{k_1})
+\sum_{k_2\in K_{N,p,n_0}:\, k_2<k_3}\varphi(p^{k_2}) 
\Bigr] \nonumber\\
&=& \sum_{{\bf t}\in T}\varphi(l)\delta({\bf t})
\cdot\delta_{(k_1\neq k_2,\, 
\max(k_1,k_2)\in K_{M,p,m_0}\cap K_{N,p,n_0})}\cdot
\frac{1}{p-1}.\label{8.84}
\end{eqnarray}
Here we use $K_{L,p,l_0}\supset\Z_{[0,k_3]}$ if
$k_3\in K_{M,p,m_0}\cap K_{N,p,n_0}$, 
we identify $\max(k_1,k_2)$ with $k_3$, 
we identify $l$ with $l_0p^{k_3}$, 
and we use
\begin{eqnarray}\label{8.85}
\sum_{k=0}^{k_3-1}\varphi(p^k)=p^{k_3-1}
=\varphi(p^{k_3})\cdot \frac{1}{p-1}.
\end{eqnarray}
Now we simplify the lines four to seven in 
\eqref{8.83}.
We will rename $k_1=k_2$ in the lines four and five
as $k_4$, and write in the lines six and seven   
the sum over $k_3$
as $|\Z_{[k+1,k_4-1]}-(K_{M,p,m_0}\cup K_{N,p,n_0})|$
respectively as
$|\Z_{[k+1,k_4-1]}\cap K_{M,p,m_0}\cap K_{N,p,n_0}|$. 
Here we use $K_{L,p,l_0}\supset \Z_{[0,k_4]}$ for
$k_4\in K_{M,p,m_0}\cap K_{N,p,n_0}$. 
Then we obtain for the lines four to seven of formula
\eqref{8.83}
\begin{eqnarray}\nonumber
&&\sum_{{\bf t}_0\in T_0}\varphi(l_0)\delta({\bf t}_0)
\cdot \Bigl(\\
&&\sum_{k\in K_{L,p,l_0}} \varphi(p^k)
\sum_{k_4\in K_{M,p,m_0}\cap K_{N,p,n_0}:\, k_4>k}
 \varphi(p^{k_4})\cdot \bigl[|\Z_{[k+1,k_4-1]}|
 +\frac{p-2}{p-1}\bigr] \Bigr)\nonumber\\
&=& \sum_{{\bf t}_0\in T_0}\varphi(l_0)\delta({\bf t}_0)
\cdot \Bigl( \nonumber\\
&&\sum_{k_4\in K_{M,p,m_0}\cap K_{N,p,n_0}:\, k_4>0}
\varphi(p^{k_4}) \cdot
\sum_{k=0}^{k_4-1} \varphi(p^k)
\bigl[ k_4-\frac{1}{p-1}-k\bigr]\Bigr) \nonumber\\
&=& \sum_{{\bf t}_0\in T_0}\varphi(l_0)\delta({\bf t}_0)
\sum_{k_4\in K_{M,p,m_0}\cap K_{N,p,n_0}:\, k_4>0}
\varphi(p^{k_4}) 
\cdot\bigl[ p^{k_4-1}-\frac{1}{p-1}\bigr] .\label{8.86}
\end{eqnarray}
Here we used
\begin{eqnarray}\label{8.87}
\sum_{k=0}^{k_4-1}\varphi(p^k)=p^{k_4-1},\quad 
\sum_{k=1}^{k_4-1}p^{k-1}k
=\frac{(k_4-1)p^{k_4}-k_4p^{k_4-1}+1}{(p-1)^2}.
\end{eqnarray}
On the other hand,
\begin{eqnarray*}
&&\sum_{{\bf t}\in T}\varphi(l)\delta({\bf t})\cdot
\delta_{(k_1=k_2,\, k>0)}\cdot\frac{1}{p-1}\\
&=& \sum_{{\bf t}_0\in T_0}\varphi(l_0)\delta({\bf t}_0)
\frac{1}{p-1}\sum_{k_4\in K_{M,p,m_0}\cap K_{N,p,n_0}:\, k_4>0}
\sum_{k=1}^{k_4}\varphi(p^k)\delta(p^{k_4},p^{k_4},p^k)\\
&=& \sum_{{\bf t}_0\in T_0}\varphi(l_0)\delta({\bf t}_0)
\sum_{k_4\in K_{M,p,m_0}\cap K_{N,p,n_0}:\, k_4>0}
\varphi(p^{k_4})\cdot\bigl[ p^{k_4-1}-\frac{1}{p-1}\bigr]\\
&=& \textup{the last line of \eqref{8.86}.}
\end{eqnarray*}
Therefore the lines four to eight of formula \eqref{8.83}
are equal to 
\begin{eqnarray*}
&&\sum_{{\bf t}\in T}\varphi(l)\delta({\bf t})
\cdot\delta_{(k_1\neq k_2,\, 
\max(k_1,k_2)\in K_{M,p,m_0}\cap K_{N,p,n_0})}\cdot
\frac{1}{p-1}\\
&+& \sum_{{\bf t}\in T}\varphi(l)\delta({\bf t})\cdot
\delta_{(k_1=k_2,\, k>0)}\cdot\frac{1}{p-1}
= \frac{1}{2}\sum_{{\bf t}\in T}\varphi(l)\delta({\bf t})\cdot
A^{(5,2)}_{\bf t}.
\end{eqnarray*}
The formulas \eqref{8.67}--\eqref{8.69} are proved.

Now we prove the formulas for the ten parts
$A_{\bf t}^{(C)}$ of $A_{\bf t}^{(7)}$.
Lemma \ref{t8.1} will be applied.
All formulas except the formula for $A_{\bf t}^{(C3)}$
are immediate consequences of it. 
In all ten cases in Lemma \ref{t8.1} except the case
$(C3)$, we put $(m,n,l)_{\textup{Lemma \ref{t8.2}}}
:=(m_1,n_1,l)_{\textup{Lemma \ref{t8.1}}}$, 
so that $(k_1,k_2,k)_{\textup{Lemma \ref{t8.2}}}
=(v_p(m_1),v_p(n_1),v_p(l))_{\textup{Lemma \ref{t8.1}}}$. 
We sum over the possible pairs 
$(m_2,n_2)_{\textup{Lemma \ref{t8.1}}}$.

It rests to prove formula \eqref{8.74} for 
$A_{\bf t}^{(C3)}$. Here we put
$(m,n)_{\textup{Lemma \ref{t8.2}}}:=
(m_1,n_2)_{\textup{Lemma \ref{t8.1}}}$, so that 
$(k_1,k_2)_{\textup{Lemma \ref{t8.2}}}
=(v_p(m_1),v_p(n_2))_{\textup{Lemma \ref{t8.1}}}$. 
We sum over the possible values
$k_3:=v_p(m_2)_{\textup{Lemma \ref{t8.1}}}$.
They form the set $K_{M,p,m_0}\cap K_{N,p,n_0}\cap 
\Z_{>\max(k_1,k_2)}$.
We need a case discussion.

The case $k_1<k_2$: Then $k=k_2$ and 
\begin{eqnarray*}
&&\Bigl(\varphi(p^{v_p(l)})
\delta(p^{v_p(m_1)},p^{v_p(n_1)},p^{v_p(l)})
\cdot \frac{\varphi(p^{v_p(n_2)})}{\varphi(p^{v_p(n_1)})}
\Bigr)_{\textup{Lemma \ref{t8.1}}}\\
&=& \varphi(p^{k_3})
\delta(p^{k_1},p^{k_3},p^{k_3})
\cdot \frac{\varphi(p^{k_2})}{\varphi(p^{k_3})}\\
&=&\varphi(p^{k_2})\varphi(p^{k_1})
= \varphi(p^k)\delta(p^{k_1},p^{k_2},p^k).
\end{eqnarray*}
Summing over $k_3\in K_{M,p,m_0}\cap K_{N,p,n_0}\cap 
\Z_{>\max(k_1,k_2)}$ gives formula \eqref{8.74}
for $A_{\bf t}^{(C3)}$. 

The  case $k_1>k_2$: Analogous to the case $k_1<k_2$.

The case $k_1=k_2=0$: Analogous to the case $k_1<k_2$.

The case $k_1=k_2>0$: Then $k\in\Z_{[0,k_1]}$ and 
\begin{eqnarray*}
&&\Bigl(\varphi(p^{v_p(l)})
\delta(p^{v_p(m_1)},p^{v_p(n_1)},p^{v_p(l)})
\cdot \frac{\varphi(p^{v_p(n_2)})}{\varphi(p^{v_p(n_1)})}
\Bigr)_{\textup{Lemma \ref{t8.1}}}\\
&=& \varphi(p^{k_3})
\delta(p^{k_1},p^{k_3},p^{k_3})
\cdot \frac{\varphi(p^{k_2})}{\varphi(p^{k_3})}
=\varphi(p^{k_1})\varphi(p^{k_1})\\
&=& \varphi(p^{k_1})\cdot \Bigl(
\sum_{k=0}^{k_1-1}\varphi(p^k) 
+ \varphi(p^{k_1})\frac{p-2}{p-1}\Bigr)\\
&=& \sum_{k=0}^{k_1} \varphi(p^k)\cdot\delta(p^{k_1},p^{k_2},p^k).
\end{eqnarray*}
Summing over $k_3\in K_{M,p,m_0}\cap K_{N,p,n_0}\cap 
\Z_{>\max(k_1,k_2)}$ gives formula \eqref{8.74}
for $A_{\bf t}^{(C3)}$. 
This finishes the proof of Lemma \ref{t8.2}.
\hfill $\Box$

\bigskip
Lemma \ref{t8.2} makes it now easy to prove 
\eqref{8.43}, i.e. that its left hand side vanishes.

\begin{eqnarray*}
&&A_{\bf t}^{(1)} + A_{\bf t}^{(C2a)} + A_{\bf t}^{(C4a)}\\
&=& -2\cdot|K_{M,p,m_0}\cap \Z_{>\max(k_1,k_2)}|
+ \delta_{(k_2>k_1,\, k_2\in K_{M,p,m_0})}\cdot
\frac{-2}{p-1},\\
&&A_{\bf t}^{(3)} + A_{\bf t}^{(C2b)} + A_{\bf t}^{(C4b)}\\
&=& -2\cdot|K_{N,p,n_0}\cap \Z_{>\max(k_1,k_2)}|
+ \delta_{(k_1>k_2,\, k_1\in K_{N,p,n_0})}\cdot
\frac{-2}{p-1},\\
&&A_{\bf t}^{(5,1)} + A_{\bf t}^{(C1)} + A_{\bf t}^{(C3)}\\
&=& 2\cdot|K_{M,p,m_0}\cap \Z_{>\max(k_1,k_2)}|
+ 2\cdot|K_{N,p,n_0}\cap \Z_{>\max(k_1,k_2)}|.
\end{eqnarray*}
Therefore 
\begin{eqnarray}\label{8.88}
A_{\bf t}^{(\textup{part 1})}
&:=&\sum_{\tiny \begin{array}{l}(j)\in\{(1),(3),(5,1),(C1),\\
(C2a),(C2b),(C3),(C4a),(C4b)\}\end{array}} A_{\bf t}^{(j)}\\
&=& 
(\delta_{(k_2>k_1,\, k_2\in K_{M,p,m_0})} 
+ \delta_{(k_1>k_2,\, k_1\in K_{N,p,n_0})})
\cdot \frac{-2}{p-1}.\nonumber
\end{eqnarray}

For the other terms, we make a case discussion.

The case $k_1>k_2$: Then $k=k_1$ and 
\begin{eqnarray*}
A_{\bf t}^{(2)}+A_{\bf t}^{(6)}&=&0,\\
A_{\bf t}^{(4)}+A_{\bf t}^{(C5b)}&=&0,\\
A_{\bf t}^{(C5a)}=A_{\bf t}^{(C6)}&=& A_{\bf t}^{(C7)}=0,\\
A_{\bf t}^{(5,2)}+A_{\bf t}^{(\textup{part 1})}&=& 0.
\end{eqnarray*}
Therefore then 
$\sum_{\textup{all possible }(j)}A_{\bf t}^{(j)}=0,$
and \eqref{8.43} is true.

The case $k_1<k_2$: Analogous to the case $k_1>k_2$.

The case $k_1=k_2=0$: Then $k=k_1$ and 
\begin{eqnarray*}
A_{\bf t}^{(j)}=0\quad\textup{for }
(j)\in\{(2),(4),(6),(5,2),(C5a),\\(C5b),(C6),(C7),
(\textup{part 1})\}.
\end{eqnarray*}
Also then 
$\sum_{\textup{all possible }(j)}A_{\bf t}^{(j)}=0,$
and \eqref{8.43} is true.

The case $k_1=k_2>0$: Then $k\in\Z_{[0,k_1]}$ and 
\begin{eqnarray*}
0&=&A_{\bf t}^{(C5a)}=A_{\bf t}^{(C5b)}
=A_{\bf t}^{(\textup{part 1})},\\
A_{\bf t}^{(\textup{part 2})}&:=&
A_{\bf t}^{(2)}+ A_{\bf t}^{(4)} +A_{\bf t}^{(6)}
+ A_{\bf t}^{(5,2)}+ A_{\bf t}^{(C6)} +A_{\bf t}^{(C7)}\\
&=& -k_1+\delta_{(k_1>k)þ}\cdot\frac{1}{p-1}
+\delta_{(k>0)}\cdot (k+\frac{1}{p-1}).
\end{eqnarray*}
This sum does not vanish for a single $k$.
But summing over $k\in\Z_{[0,k_1]}$, we obtain 0:
\begin{eqnarray*}
&&\sum_{k=0}^{k_1}\varphi(p^k)\delta(p^{k_1},p^{k_2},p^k)
\cdot A_{\bf t}^{(\textup{part 2})}\\
&=& \varphi(p^{k_1})\sum_{k=0}^{k_1}\varphi(p^k)
(1-\delta_{(k_1=k)}\frac{1}{p-1})\cdot 
A_{\bf t}^{(\textup{part 2})}\\
&=&\varphi(p^{k_1})\Bigl((-k_1+\frac{1}{p-1})
+ \sum_{k=1}^{k_1-1}p^{k-1}((-k_1+k)(p-1)+2)\\
&& +\ p^{k_1-1}(p-2)\frac{1}{p-1}\Bigr)\\
&=& \varphi(p^{k_1})\Bigl((-k_1+\frac{1}{p-1})
+ (p^{k_1-1}-1)(-k_1+\frac{2}{p-1})+\\
&&+\ \frac{(k_1-1)p^{k_1}-k_1p^{k_1-1}+1}{p-1}
+\ p^{k_1-1}\frac{p-2}{p-1}\Bigr)\\
&=& \varphi(p^{k_1})\cdot 0 = 0.
\end{eqnarray*}
Here we used \eqref{8.87}.
Therefore then for any fixed 
$(m,n,l_0)\in M\times N\times \pi_p(L)$ the sum
\begin{eqnarray*}
\sum_{k\in K_{L,p,l_0}}\varphi(l_0p^k)\delta(m,n,l_0p^k)
\cdot \sum_{\textup{all possible }(j)}A_{(m,n,l_0p^k)}^{(j)}
\end{eqnarray*}
vanishes. Thus \eqref{8.43} is true also in the case
$k_1=k_2>0$. 
This finishes the proof of Theorem \ref{t7.4}.
\hfill $\Box$

\begin{remark}\label{t8.3}
At the heart of the proof of Theorem \ref{t7.4} in this
section is the proof that the set $\CC^{dec}$ in 
\eqref{8.9} is a $\Z$-basis of $G\otimes H$ in the special 
case $M,N\subset\{p^k\,|\, k\in\N_0\}$ 
for some prime number $p$ 
if the pair $(M,N)$ is sdiOb-sufficient.
Our proof for this case is long and tedious. 
It would be desirable to have an elegant or short proof
for this case. 
\end{remark}

\section{A compatibility condition for sets of orders of 
eigenvalues of Orlik blocks}\label{c9}
\setcounter{equation}{0}

\noindent
This section proposes and discusses a condition for finite sets 
$M\subset \N$ of orders of eigenvalues of Orlik blocks.
It is given in Definition \ref{t9.4} (c).
It has a number of good properties, which will be 
given below in Lemma \ref{t9.5}, 
Theorem \ref{t9.6}, Lemma \ref{t9.8}, 
Theorem \ref{t9.9}, Theorem \ref{t9.10}
and Lemma \ref{t9.12}. 
Theorem \ref{t9.9} contains the statement that
$\Or(M)\otimes\Or(N)$ admits a standard decomposition
into Orlik blocks if the finite sets $M$ and $N\subset \N$
satisfy this condition.
This is prepared by several definitions.

\begin{definition}\label{t9.1}
(a) An {\it excellent order} $\succ$ on a set
$\Z_{[0,s(\succ)]}$ for some bound $s(\succ)\in\N_0$
is a strict order (so transitive and for all 
$a,b\in \Z_{[0,s(\succ)]}$ either $a=b$ or
$a\succ b$ or $b\succ a$) which is determined by the
set 
\begin{eqnarray}\label{9.1}
S(\succ)&:=& \{k\in\Z_{[0,s(\succ)]}\, |\, k\succ 0\}
\end{eqnarray}
in the following way:
\begin{eqnarray}\label{9.2}
\left. \begin{array}{l}
\succ\textup{ equals }>\textup{ on }S(\succ)\cup\{0\},\\
\succ\textup{ equals }<\textup{ on }
\Z_{[0,s(\succ)]}-S(\succ).\end{array}\right\}
\end{eqnarray}
($S(\succ)=\emptyset$ is allowed.) 
The maximal element of $\Z_{[0,s(\succ)]}$ with  respect
to $\succ$ is called $s^+(\succ)$, so $s^+(\succ)\succ k$
for any other element $k\in\Z_{[0,s(\succ)]}$.

(b) The trivial excellent order is $\succ_0$
with $s(\succ_0):=0$, so it is the empty order on 
$\Z_{[0,s(\succ_0)]}=\{0\}$ 
(and, of course $S(\succ_0)=\emptyset$).

\medskip
(c) The tensor product 
of two excellent orders $\succ_1$ and $\succ_2$
is the excellent order $\succ_1\otimes\succ_2$ with 
\begin{eqnarray}\label{9.3}
s(\succ_1\otimes \succ_2)&:=&\max(s(\succ_1),s(\succ_2))
\quad\textup{and}\\
S(\succ_1\otimes \succ_2)&:=&
(S(\succ_1)\cup S(\succ_2))-(S(\succ_1)\cap S(\succ_2)).
\label{9.4}
\end{eqnarray}
\end{definition}

\begin{examples}\label{t9.2}
(i) The excellent order $\succ_1$ with
$s(\succ_1)=7$ and $S(\succ_1)=\{6,4,1\}$ is given by
\begin{eqnarray*}
s^+(\succ_1)=6\succ_1 4\succ_1 1\succ_1 0 \succ_1 2
\succ_1 3\succ_1  5\succ_1 7.
\end{eqnarray*}
(ii)  The excellent order $\succ_2$ with
$s(\succ_2)=6$ and $S(\succ_2)=\{6,5,2,1\}$ is given by
\begin{eqnarray*}
s^+(\succ_2)=6\succ_2 5\succ_2 2\succ_2 1 
\succ_2 0\succ_2 3\succ_2 4.
\end{eqnarray*}
(iii) The excellent order $\succ_3:=(\succ_1\otimes \succ_2)$
for $\succ_1$ and $\succ_2$ in (i) and (ii) satisfies
$s(\succ_3)=7$, $S(\succ_3)=\{5,4,2\}$ and is given by
\begin{eqnarray*}
5\succ_3 4\succ_3 2\succ_3 0 \succ_3 1\succ_3 3\succ_3 6
\succ 7.
\end{eqnarray*}
(iv) For any excellent order $\succ$, the tensor product
with the trivial excellent order is $\succ$ itself, 
$\succ\otimes \succ_0=\succ$.
\end{examples}

\begin{definition}\label{t9.3}
(a) A {\it path} in a finite directed graph $(V,E)$ 
(so $V$ is a finite non-empty set and $E\subset V\times V$)
is a  tuple $(v_1,...,v_l)$ for some $l\in\Z_{\geq 2}$ 
with $v_j\in V$ and $(v_j,v_{j+1})\in E$ for 
$j\in\Z_{[1,l-1]}$. It is a path from $v_1$ to $v_l$,
so with source $v_1$ and target $v_l$.

\medskip
(b) A finite directed graph $(V,E)$ {\it has a center}
$v_V\in V$ if it has no path from
any vertex to itself and if it has at least one path from
$v_V$ to any other vertex $v\in V$.
(The center is unique, which justifies the notation
$v_V$.)

\medskip
(c) Consider a tuple $(\succ_p)_{p\in P}$ of excellent 
orders for a finite set $P\subset \PP$ of prime numbers.
It defines a finite directed graph $(V,E_V)$ with center $v_V$
as follows. Its set $V=V((\succ_p)_{p\in P})$ of vertices
is the {\it quadrant} in $\N$ 
\begin{eqnarray}\label{9.5}
V&:=& \{\prod_{p\in P}p^{k_p}\,|\, 
k_p\in\Z_{[0,s(\succ_p)]}\}.
\end{eqnarray}
Its set of edges $E_V=E((\succ_p)_{p\in P})$ is the set
\begin{eqnarray}
E_V&:=& \bigcup_{p\in P}E_{V,p}\textup{ with }\label{9.6}\\
E_{V,p}&:=& \{(m_a,m_b)\in V\times V\,|\, 
\pi_p(m_a)=\pi_p(m_b), v_p(m_a)\succ_p v_p(m_b)\}.
\nonumber
\end{eqnarray}
The edges in $E_{V,p}$ are called $p$-edges.
So, the underlying undirected graph coincides with the
undirected graph which underlies $(V,E(V))$
(defined in Definition \ref{t6.1} (c)).
But the directions of edges might have changed.
The graph $(V,E_V)$ is obviously centered with center
\begin{eqnarray}\label{9.7}
v_V&=&\prod_{p\in P}p^{s^+(\succ_p)}.
\end{eqnarray}
\end{definition}

\begin{definition}\label{t9.4}
(a) Let $\succ$ be an excellent order on the set
$\Z_{[0,s(\succ)]}$. A set $K\subset \N_0$ is 
{\it subset compatible} with $\succ$ if a bound
$k_K\in\Z_{[0,s(\succ)]}$ with 
\begin{eqnarray}\label{9.8}
K&=&\{k\in \Z_{[0,s(\succ)]}\, |\, k\succ k_K\}
\end{eqnarray}
exists or if $K=\Z_{[0,s(\succ)]}$. 
$(k_K=s^+(\succ)$ gives $K=\emptyset$, which is allowed.)

\medskip
(b) For a finite non-empty set $M\subset\N$, let
\begin{eqnarray}\label{9.9}
\PP(M)&:=& \{p\in\PP\,|\, M\neq \pi_p(M)\}\\
&=& \{p\in \PP\,|\, \exists\ m\in M\textup{ with }v_p(m)>0\}
\nonumber
\end{eqnarray}
be the set of prime numbers which turn up as factors of some
numbers in $M$.

\medskip
(c) A finite non-empty set $M\subset\N$ is
{\it compatible} with a tuple $(\succ_p)_{p\in P}$
of excellent orders for a finite set $P\supset \PP(M)$ 
of prime numbers if 
\begin{eqnarray}\label{9.10}
M&\subset& V((\succ_p)_{p\in P})
\end{eqnarray}
and if for any prime number $p\in\PP(M)$
and any $m_0\in\pi_p(M)$ the set $K_{M,p,m_0}$
(with $\pi_p^{-1}(m_0)\cap M
=\{m_0\cdot p^k\,|\, k\in K_{M,p,m_0}\}$, see \eqref{7.28})
is subset compatible with $\succ_p$.
(So, here the $\succ_p$ for $p\in P-\PP(M)$
are irrelevant. But $P\supset \PP(M)$ instead
of $P=\PP(M)$ will be useful.)

\medskip
(d) A map $\chi:\N\to\N_0$ with finite support
$\supp(\chi):=\{m\in\N\,|\, \chi(m)\neq 0\}$ is 
{\it compatible} with a tuple $(\succ_p)_{p\in P}$
of excellent orders for a finite set 
$P\supset \PP(\supp(\chi))$ of prime numbers if 
\begin{eqnarray}\label{9.11}
\supp(\chi)&\subset& V((\succ_p)_{p\in P})
\end{eqnarray}
and if for any edge $(m_a,m_b)\in E_V$
\begin{eqnarray}\label{9.12}
\chi(m_a)&\geq & \chi(m_b).
\end{eqnarray}

\medskip
(e) A {\it covering} of a map $\chi:\N\to \N_0$
with finite support is a tuple 
$(M_1,...,M_l)$ ($l\in\N_0$) 
of finite non-empty sets $M_j\subset\N$ with
\begin{eqnarray}\label{9.13}
\chi(m)&=& |\{j\in\{1,...,l\}\,|\, m\in M_j\}|
\textup{ for any }m\in\N.
\end{eqnarray} 
Here obviously $l\geq\max(\chi(m)\,|\, m\in\N_0)=:l_\chi$.
In the case $\supp(\chi)=\emptyset$ we have 
$l=0$ and an empty tuple.
The {\it standard covering} of $\chi$ is the tuple
$(M_1^{(st)},...,M_{l_\chi}^{(st)})$ with 
\begin{eqnarray}\label{9.14}
M_j^{(st)}=\{m\in \supp(\chi)\,|\, \chi(m)\geq j\}
\textup{ for }j\in\{1,...,l_\chi\}.
\end{eqnarray}
It is the unique covering with $M_1\supset ...\supset M_l$,
and it satisfies $M_1^{(st)}=\supp(\chi)$. 

\medskip
(f) Let $\chi:\N\to \N_0$ have finite support, 
let $P\supset \PP(\supp(\chi))$ be a finite set of prime
numbers, and let
$(\succ_p)_{p\in P}$ be a tuple of excellent orders
with \eqref{9.11}.
A covering $(M_1,...,M_l)$ of $\chi$ is called
{\it compatible} with $(\succ_p)_{p\in P}$
if each set $M_j$ is compatible with  $(\succ_p)_{p\in P}$. 
\end{definition}

The following lemma expresses the compatibility conditions
in Definition \ref{t9.4} (c) and (d) in a different way,
and it shows their relationship.

\begin{lemma}\label{t9.5}
(a) Let $M\subset\N$ be a finite non-empty set,
let $P\supset \PP(M)$ be a finite set of prime
numbers, 
and let $(\succ)_{p\in P}$ be a tuple of excellent orders
with \eqref{9.10}. (Recall the definition of $(V,E_V,v_V)$ 
in Definition \ref{t9.3} (c).) 
The following three conditions are equivalent:
\begin{list}{}{}
\item[(i)]
$M$ is compatible with $(\succ)_{p\in P}$.
\item[(ii)]
$(M,E_V\cap M\times M)$ is a directed graph with center $v_V$
(so $v_V\in M)$, and if $M$ contains the target of a path
in $V$, it contains all vertices in this path.
\item[(iii)]
If $m_b\in M$ and $(m_a,m_b)\in E_V$ then 
$m_a\in M$. 
\end{list}

\medskip
(b) Let $\chi:\N\to\N_0$ be a map with finite support,
let $P\supset \PP(\supp(\chi))$ be a finite set of prime
numbers, and let $(\succ_p)_{p\in P}$ be a tuple of 
excellent orders with \eqref{9.11}. 
The following three conditions are equivalent.
\begin{list}{}{}
\item[(i)] 
$\chi$ is compatible with $(\succ_p)_{p\in P}$.
\item[(ii)] 
$\chi$ has a covering $(M_1,...,M_l)$ which is compatible
with $(\succ_p)_{p\in P}$.
\item[(iii)]
The standard covering of $\chi$ is compatible with
$(\succ_p)_{p\in P}$.
\end{list}
\end{lemma}

{\bf Proof:}
(a) (i)$\Rightarrow$(iii): 
Suppose $m_b\in M$ and $(m_a,m_b)\in E_{V,p}$
for some prime number $p\in P$. 
Then $\pi_p(m_a)=\pi_p(m_b)$ and 
\begin{eqnarray*}
v_p(m_a)\succ_p v_p(m_b)\in K_{M,p,\pi_p(m_b)}
=K_{M,p,\pi_p(m_a)}.
\end{eqnarray*}
As $M$  is compatible with $(\succ_q)_{q\in P}$,
also $v_p(m_a)\in K_{M,p,\pi_p(m_a)}$.
Therefore $m_a\in M$. 

(iii)$\Rightarrow$(ii): Apply (iii) several times. 

(ii)$\Rightarrow$(i): 
Consider a prime number $p\in \PP(M)$,
an element $m_0\in \pi_p(M)$, any number
$k_1\in K_{M,p,m_0}$ and any number 
$k_2\in \Z_{[0,s(\succ_p)]}$ with $k_2\succ_p k_1$.
We have to show $k_2\in K_{M,p,m_0}$. 

Define $m_i:=m_0\cdot p^{k_i}$ for $i\in\{1,2\}$.
Then $m_1\in M$, and we have to show $m_2\in M$.
But obviously $(m_2,m_1)$ is a path in $(V,E_V)$.
(ii) applies and yields $m_2\in M$. 

\medskip
(b) The case $\supp(\chi)=\emptyset$ is trivial.
We suppose $\supp(\chi)\neq\emptyset$. 

(i)$\Rightarrow$(iii): 
Let $M_j^{(st)}$ be one of the sets in the standard covering
of $\chi$. It is sufficient to prove property (iii) in (a)
for $M_j^{(st)}$. Suppose $m_b\in M_j^{(st)}$
and $(m_a,m_b)\in E_V$. 
Because of \eqref{9.12} and $(m_a,m_b)\in E_V$,
$\chi(m_a)\geq \chi(m_b)$. By definition of $M_j^{(st)}$,
then $m_a\in M_j^{(st)}$. 

(iii)$\Rightarrow$(ii): Trivial.

(ii)$\Rightarrow$(i): 
Let $(M_1,...,M_l)$ be a covering of $\chi$ which is 
compatible with $(\succ_p)_{p\in P}$. 
Let $(m_a,m_b)\in E_{V}$. We have to show 
$\chi(m_a)\geq \chi(m_b)$. 
Consider one of the sets $M_j$ with $m_b\in M_j$. 
By hypothesis it is compatible with $(\succ_p)_{p\in P}$.
Part (a) and $m_b\in M_j$ give $m_a\in M_j$.
Therefore 
\begin{eqnarray*}
\chi(m_a)&=&|\{i\in\{1,...,l\}\,|\, m_a\in M_i\}|\\
&\geq & |\{i\in\{1,...,l\}\,|\, m_b\in M_i\}|=\chi(m_b).
\end{eqnarray*}
\hfill$\Box$ 

\bigskip

Theorem \ref{t9.6} compares different decompositions
into Orlik blocks.

\begin{theorem}\label{t9.6}
Let $(\succ_p)_{p\in P}$ be at tuple of excellent orders
for a finite set $P$ of prime numbers, and let 
$\chi:\N\to \N_0$ be a map with finite support which
is compatible with $(\succ_p)_{p\in P}$. 
Let $(M_1^{(1)},...,M_{l_1}^{(1)}))$ and $(M_1^{(2)},...,
M_{l_2}^{(2)})$ be two coverings of $\chi$ which
are both compatible with $(\succ_p)_{p\in P}$.
Then the corresponding sums of Orlik blocks
are isomorphic,
\begin{eqnarray}\label{9.15}
\bigoplus_{i=1}^{l_1}\Or(M_i^{(1)})&\cong &
\bigoplus_{j=1}^{l_2}\Or(M_j^{(2)}),
\end{eqnarray}
and $l_1=l_2=l_\chi(:=\max(\chi(m)\,|\, m\in\N))$.
\end{theorem}

{\bf Proof:} At the end of the proof, we will apply 
Theorem \ref{t5.1} (b). But before that, 
the main work is the discussion how coverings of $\chi$ 
which are compatible with $(\succ_p)_{p\in P}$ 
can look like. There is not so much freedom. 

Let $(M_1^{(st)},...,M_{l_\chi}^{(st)})$ be the standard
covering of $\chi$. For each $j\in\Z_{[1,l_\chi]}$
consider the set 
\begin{eqnarray*}
M_j^{(st)}-M_{j+1}^{(st)}=\{m\in \supp(\chi)\,|\, 
\chi(m)=j\}.
\end{eqnarray*}
as a subgraph of $(V,E_V)$ and denote its components by
$M_{j,i}^{(st)}$ for $i\in\Z_{[1,c_j]}$ for some 
$c_j\in\N$. We define a directed graph 
$(V^{(comp)},E^{(comp)})$ with set of vertices
\begin{eqnarray}\label{9.16}
V^{(comp)}&:=& \{M_{j,i}^{(st)}\,|\, 
j\in\Z_{[1,l_\chi]},i\in \Z_{[1,c_j]}\}
\end{eqnarray}
and set of edges
\begin{eqnarray}\label{9.17}
E^{(comp)}&:=&\{(v_1,v_2)\in V^{(comp)}\times V^{(comp)}\,|\, 
v_1\neq v_2, \textup{ an edge in }E_V \nonumber \\
&&\textup{from a vertex in }v_1
\textup{ to a vertex in }v_2\textup{ exists}\}.
\end{eqnarray}
Furthermore, we define the map
\begin{eqnarray}\label{9.18}
\chi^{(comp)}:V^{(comp)}\to\N_0,\quad
M_{j,i}^{(st)}\mapsto j.
\end{eqnarray}
By Lemma \ref{t9.5} (b) and by hypothesis,
the standard covering is compatible with
$(\succ_p)_{p\in P}$. Therefore
for any edge $(v_1,v_2)\in E^{(comp)}$
$\chi^{(comp)}(v_1)\geq \chi^{(comp)}(v_2)$. 
In fact, for any edge $(v_1,v_2)\in E^{(comp)}$ even
\begin{eqnarray}\label{9.19}
\chi^{(comp)}(v_1)> \chi^{(comp)}(v_2)
\end{eqnarray}
holds. Because if on the contrary 
$\chi^{(comp)}(v_1)=\chi^{(comp)}(v_2)$,
then $v_1$ and $v_2$ would be subsets of 
the subgraph $M_j^{(st)}-M_{j+1}^{(st)}$ 
where $j:=\chi^{(comp)}(v_1)$, 
and because of $(v_1,v_2)\in E^{(comp)}$
an edge from a vertex in $v_1$ to a vertex in $v_2$
would exist, so they 
would not be different components. 

$M_1^{(st)}=\supp(\chi)$ is a directed graph with
center $v_V$. Therefore and because of \eqref{9.19},
also $(V^{(comp)},E^{(comp)})$ is a directed graph
with center. Its center is the component 
$M_{l_\chi,1}^{(st)}$, which contains $v_V$, and
which is the only component of $M_{l_\chi}^{(st)}$. 

The following claim makes the shape of any
covering $(M_1,...,M_{l_1})$ of $\chi$ which is 
compatible with $(\succ_p)_{p\in P}$ explicit.

\medskip
{\bf Claim 1:} {\it A covering $(M_1,...,M_{l_1})$ of $\chi$
is compatible with $(\succ_p)_{p\in P}$ if and only
if the following holds.

$l_1=l_\chi$, and a tuple $(V_1,...,V_{l_\chi})$
of subsets of $V^{(comp)}$ exists such that
\begin{eqnarray}\label{9.20}
M_j&=& \bigcup_{v\in V_j}v \subset \supp(\chi),\\
\chi^{(comp)}(v)&=& |\{j\in \Z_{[1,l_\chi]}\,|\, v\in V_j\}|,
\label{9.21}
\end{eqnarray}
and such that each subset $V_j$ is as a subgraph
of $V^{(comp)}$ a directed graph with center
$M_{l_\chi}^{(st)}$ and contains all vertices of a path
in $V^{(comp)}$ whose target it contains.}

\medskip
(Examples are given in Examples \ref{t9.7}.)

{\bf Proof of Claim 1:}
$\Rightarrow$: Suppose that $(M_1,...,M_{l_1})$ is a
covering of $\chi$ which is compatible with
$(\succ_p)_{p\in P}$. First we will prove \eqref{9.20}
and \eqref{9.21} (with $l_1$ instead of $l_\chi$).

Consider an edge $(m_a,m_b)\in E_V$ with
$m_a$ and $m_b$ in the same component $M_{j,i}^{(st)}$. 
Then 
\begin{eqnarray*}
|\{k\in\Z_{[1,l_1]}\,|\, m_a\in M_k\}|&=& \chi(m_a)
=\chi^{(comp)}(M_{j,i}^{(st)})=j\\
&=&\chi(m_b)=|\{k\in\Z_{[1,l_1]}\,|\, m_b\in M_k\}|.
\end{eqnarray*}
Let $M_k$ be one subset which contains $m_b$.
Because $M_k$ is compatible with $(\succ_p)_{p\in P}$ 
and because of Lemma \ref{t9.5} (a), also $m_a\in M_k$. 
Therefore $m_a$ and $m_b$ are elements of the
same sets $M_k$. 

Because $M_{j,i}^{(st)}$ is a connected subgraph
of $(V,E_V)$, this implies that these sets $M_k$
contain all of $M_{j,i}^{(st)}$
and that all other sets $M_{\www k}$ do not intersect
$M_{j,i}^{(st)}$. 

Therefore a tuple $(V_1,...,V_{l_1})$ of subsets of 
$V^{(comp)}$ with \eqref{9.20} and \eqref{9.21}
(with $l_1$ instead of $l_\chi$) exists. 

Suppose that $v_2\in V_j$ for some $j\in\Z_{[1,l_1]}$
and that $(v_1,v_2)\in E^{(comp)}$ is an edge.
Then an edge $(m_1,m_2)\in E_V$ with 
$m_1\in v_1$ and $m_2\in v_2\subset M_j$ exists. 
$M_j$ is compatible with $(\succ_p)_{p\in P}$.
Lemma \ref{t9.5} (a) implies $m_1\in M_j$. 
Therefore $v_1\in V_j$.

This shows that $V_j$ contains any vertex $v\in V^{(comp)}$
such that a path to a vertex $v_2\in V_j$ exists.
And especially, therefore $V_j$ is a directed subgraph
of $V^{(comp)}$ with center $M_{l_\chi}^{(st)}$. 

Any set $V_j$ contains the vertex $M_{l_\chi}^{(st)}$,
and $\chi^{(comp)}(M_{l_\chi}^{(st)})=l_\chi$.
With \eqref{9.21} this gives $l_1=l_\chi$. 

\medskip
$\Leftarrow$: \eqref{9.20} and \eqref{9.21}
show that $(M_1,...,M_{l_\chi})$ is a covering of $\chi$. 
It remains to show that each set $M_j$ is compatible
with $(\succ_p)_{p\in P}$. 

Let $(m_a,m_b)\in E_V$ be an edge and $m_b\in M_j$
for some $j\in\Z_{[1,l_\chi]}$. It is sufficient to 
show $m_a\in M_j$, because then one can apply
Lemma \ref{t9.5} (a) (iii)$\Rightarrow$(i). 

Let $v_a\in V^{(comp)}$ respectively $v_b\in V_j$ 
be the vertex of $V^{(comp)}$ which contains $m_a$ 
respectively $m_b$.

$\chi(m_a)<\chi(m_b)$ is impossible because $\chi$
is compatible with $(\succ_p)_{p\in P}$. 

If $\chi(m_a)=\chi(m_b)$, then $v_a=v_b$, and 
$m_a\in M_j$ because of \eqref{9.20}.

If $\chi(m_a)>\chi(m_b)$, then $(v_a,v_b)\in E^{(comp)}$.
Then the hypothesis on $V_j$ implies $v_a\in V_j$.
And $m_a\in M_j$ because of \eqref{9.20}.
This finishes the proof of Claim 1.
\hfill ($\Box$) 

\medskip
{\bf Claim 2:} {\it 
(a) Let $(M_1,...,M_{l_\chi})$ be a covering of $\chi$
which is compatible with $(\succ_p)_{p\in P}$,
and which satisfies $M_1\not\subset M_2$ and 
$M_2\not\subset M_1$. Then also
$(M_1\cup M_2,M_1\cap M_2,M_3,...,M_{l_\chi})$ is a
covering of $\chi$ which is compatible with
$(\succ_p)_{p\in P}$. 

\medskip
(b) Iterating the procedure in (a), one can go from any
covering of $\chi$ which is compatible with
$(\succ_p)_{p\in P}$ to the standard covering and
thus also to any other covering of $\chi$ which is 
compatible with $(\succ_p)_{p\in P}$. }

\medskip
{\bf Proof of Claim 2:}
(a) Because of \eqref{9.20} and \eqref{9.21} it is clear
that $(M_1\cup M_2,M_1\cap M_2,M_3,...,M_{l_\chi})$
is a covering of $\chi$ and that it satisfies
\eqref{9.20} and \eqref{9.21}. 
The corresponding tuple of subsets of $V^{(comp)}$ is
$(V_1\cup V_2,V_1\cap V_2,V_3,...,V_{l_\chi})$. 
It is also clear that $V_1\cup V_2$ and $V_1\cap V_2$
are directed subgraphs of $V^{(comp)}$ with center
$M_{l_\chi}^{(st)}$ and that $V_1\cup V_2$
respectively $V_1\cap V_2$ contains all vertices of 
a path in $V^{(comp)}$ whose target it contains.

(b) Iterating the procedure in (a), one can increase the
set $M_1$ until it is $M_1^{(st)}=\supp(\chi)$.
Then again iterating the procedure in (a), but leaving
$M_1$ as it is, one can increase the set $M_2$
until it is $M_2^{(st)}$. One can continue this until
one reaches the standard covering of $\chi$. 
\hfill $(\Box$)

\medskip
{\bf Claim 3:} In the situation of part (a) of Claim 2,
\begin{eqnarray}\label{9.22}
\bigoplus_{j=1}^{l_\chi}\Or(M_j)
&\cong& \Or(M_1\cup M_2)\oplus \Or(M_1\cap M_2)\oplus
\bigoplus_{j=3}^{l_\chi}\Or(M_j).\hspace*{1cm}
\end{eqnarray}

\medskip
{\bf Proof of Claim 3:} We will apply Theorem \ref{t5.1} (b).
Define
\begin{eqnarray*}
f_1:=1,\, f_2:=\prod_{m\in M_2-M_1}\Phi_m,\ 
f_3:=\prod_{m\in M_1\cap M_2}\Phi_m,
\ f_4:=\prod_{m\in M_1-M_2}\Phi_m.
\end{eqnarray*}
Then 
\begin{eqnarray*}\begin{array}{rclrcl}
H^{[f_1f_3f_4]}&\cong& \Or(M_1),
& H^{[f_2f_3]}&\cong&\Or(M_2),\\
H^{[f_1f_3]}&\cong& \Or(M_1\cap M_2), 
& H^{[f_2f_3f_4]}&\cong& \Or(M_1\cup M_2).\end{array}
\end{eqnarray*}
$|\Res(f_1,f_4)|=\Res(1,f_4)=1$ is trivially true. 
We want to show $|\Res(f_2,f_4)|=1$.
Suppose $|\Res(f_2,f_4)|>1$.  Because of
\eqref{4.24} and \eqref{4.8} there is an edge
in $E(\N)$ which connects a vertex $m_1\in M_1-M_2$
with a vertex $m_2\in M_2-M_1$. 
Then either $(m_1,m_2)$ or $(m_2,m_1)$ 
is an edge in $E_V$. For $j\in\{1,2\}$, 
let $V_j\subset V^{(comp)}$ be the subset with
$M_j=\bigcup_{v\in V_j}v_j$, and let $v_j\in V_j$
be the vertex with $m_j\in v_j$. 
Then $(v_1,v_2)$ or $(v_2,v_1)$ is an edge in 
$E^{(comp)}$, but $v_1\in V_1-V_2$ and $v_2\in V_2-V_1$,
as $m_1\in M_1-M_2$ and $m_2\in M_2-M_1$.
This is a contradiction to the property of $V_1$ and $V_2$
that they contain with the endpoint of a path in $V^{(comp)}$
also the source of the path. 
Because of this contradiction $|\Res(f_2,f_4)|=1$.

Theorem \ref{t5.1} (b) can be applied and gives
\begin{eqnarray}\label{9.23}
\Or(M_1)\oplus \Or(M_2)&\cong& \Or(M_1\cup M_2)\oplus
\Or(M_1\cap M_2).
\end{eqnarray}
This implies \eqref{9.22}
and finishes the proof of Claim 3.\hfill ($\Box$) 

\medskip
Claim 3 and part (b) of Claim 2 show that any
covering $(M_1,...,M_{l_\chi})$ of $\chi$ which is
compatible with $(\succ_p)_{p\in P}$ satisfies
\begin{eqnarray}\label{9.24}
\bigoplus_{j=1}^{l_\chi}\Or(M_j) &\cong &
\bigoplus_{j=1}^{l_\chi}\Or(M_j^{(st)}).
\end{eqnarray}
This implies immediately also \eqref{9.15}
(and there $l_1=l_2=l_\chi$).
This finishes the proof of Theorem \ref{t9.6}.
\hfill $\Box$

\begin{examples}\label{t9.7}
(i) Fix two prime numbers $p_1$ and $p_2$ and four numbers
$k_1,k_2,k_3,k_4\in\N$ with $k_1<k_3$, $k_2<k_4$, and 
define the following three {\it rectangles} of numbers,
\begin{eqnarray*}
N_0&:=& \{p_1^{l_1}p_2^{l_2}\,|\, 
(l_1,l_2)\in\Z_{[0,k_1]}\times \Z_{[0,k_2]}\},\\
N_1&:=& \{p_1^{l_1}p_2^{l_2}\,|\, 
(l_1,l_2)\in\Z_{[k_1+1,k_3]}\times \Z_{[0,k_2]}\},\\
N_2&:=& \{p_1^{l_1}p_2^{l_2}\,|\, 
(l_1,l_2)\in\Z_{[0,k_1]}\times \Z_{[k_2+1,k_4]}\}.
\end{eqnarray*}
Define a map $\chi:\N\to\N_0$ with finite support by
\begin{eqnarray*}
\chi(m):=\left\{\begin{array}{ll}
2&\textup{if }m\in N_0,\\
1&\textup{if }m\in N_1\cup N_2,\\
0&\textup{if }m\in \N-(N_0\cup N_1\cup N_2).
\end{array}\right. 
\end{eqnarray*}
Then the standard covering of $\chi$ is given by 
$l_\chi=2$ and
\begin{eqnarray*}
M_2^{(st)}=N_0\subset M_1^{(st)}=N_0\cup N_1\cup N_2.
\end{eqnarray*}

Consider the excellent orders $\succ_{p_1}$ and
$\succ_{p_2}$ with
\begin{eqnarray*}
s(\succ_{p_1})=k_3,\ S(\succ_{p_1})=\emptyset, 
\textup{ so }\succ_{p_1}\textup{ is }
<\textup{ on }\Z_{[0,k_3]},\\
s(\succ_{p_2})=k_4,\ S(\succ_{p_2})=\emptyset, 
\textup{ so }\succ_{p_2}\textup{ is }
<\textup{ on }\Z_{[0,k_4]}.
\end{eqnarray*}
$\chi$ and its standard covering are compatible with
the tuple $(\succ_{p_1},\succ_{p_2})$ of excellent orders.
The only other covering of $\chi$ which is compatible
with $(\succ_{p_1},\succ_{p_2})$, consists of 
\begin{eqnarray*}
M_1=N_0\cup N_1,\quad M_2=N_0\cup N_2,
\end{eqnarray*}
because
\begin{eqnarray*}
V^{(comp)}=\{M_2^{(st)},M_{1,1}^{(st)},M_{1,2}^{(st)}\}\quad
\textup{with }M_{1,i}^{(st)}=N_i.
\end{eqnarray*}

(ii) Keep the data from (i). Define a new map
$\www\chi:\N\to \N_0$ with finite support by
\begin{eqnarray*}
\www\chi(m):=\chi(m)+\delta_{(m=p_1^{k_1+1}p_2^{k_2+1})}.
\end{eqnarray*}
Then its standard covering is given by $l_{\www\chi}=2$ and
\begin{eqnarray*}
\www{M}_2^{(st)}=N_0,\ 
\www{M}_1^{(st)}=N_0\cup N_1\cup N_2 \cup 
\{p_1^{k_1+1}p_2^{k_2+1}\}.
\end{eqnarray*}
$\www\chi$ and its standard covering are compatible with
the tuple $(\succ_{p_1},\succ_{p_2})$ of excellent orders.
No other covering of $\www\chi$ is compatible with
$(\succ_{p_1},\succ_{p_2})$, because the difference set
$\www{M}_1^{(st)}-\www{M}_2^{(st)}$ has only one 
component, so $V^{(comp)}=\{\www{M}_2^{(st)},
\www{M}_1^{(st)}-\www{M}_2^{(st)}\}$.
\end{examples}

Now we turn to a situation which is motivated
by Theorem \ref{t7.4}.
Two finite non-empty sets $M\subset \N$ and $N\subset\N$ 
give rise to the map 
$\chi:\N\to\N_0$ defined in \eqref{7.19} with finite
support $L$ (as in \eqref{7.20}).

First we consider a special case.

\begin{lemma}\label{t9.8}
Let $\succ_1$ and $\succ_2$ be two excellent orders.
Let $K_1\subset\N_0$ and $K_2\subset \N_0$ be
two finite non-empty sets such that $K_j$
is subset compatible with $\succ_j$ (Definition \ref{t9.4} (a))
for $j\in\{1,2\}$. Let $p$ be a prime number.
Define $M:=\{p^k\,|\, k\in K_1\}$ and 
$N:=\{p^k\,|\, k\in K_2\}$. Consider $\chi:\N\to\N_0$
as in \eqref{7.19}.

Then $\chi$ is compatible with $\succ_1\otimes\succ_2$
(considered as a tuple of excellent orders over 
$P=\{p\}$).
\end{lemma}

{\bf Proof:}
Define for $j\in\{1,2\}$ 
\begin{eqnarray}\label{9.25}
f(K_j)&:=& \left\{\begin{array}{ll}
\max(k\in\N_0\,|\, \Z_{[0,k]}\subset K_j) & \textup{ if }
0\in K_j,\\
-1 & \textup{ if }0\notin K_j,\end{array}\right.\\
e(K_j)&:=& \left\{\begin{array}{ll}
\max(k\in\N_0\,|\, \Z_{[0,k]}\cap K_j=\emptyset) & \textup{ if }
0\notin K_j,\\
-1 & \textup{ if }0\in K_j,\end{array}\right.
\label{9.26}\\
k^{(j)}&:=& \max(f(K_j),e(K_j)),\label{9.27}\\
k^{max}&:=& \max(k^{(1)},k^{(2)}).\label{9.28}
\end{eqnarray}
Then either ($f(K_j)\geq 0$ and $e(K_j)=-1$) or
($f(K_j)=-1$ and $e(K_j)\geq 0$).
In both cases $k^{(1)},k^{(2)},k^{max}\in\N_0$.

\medskip
{\bf Claim 1:} {\it $\chi$ has constant values on
$\{p^k\,|\, k\in \Z_{[0,k^{max}]}\}$. }

\medskip
{\bf Proof of Claim 1:}
Here $\chi(p^k)$ for any $k\in\N_0$ 
is given by a formula similar to \eqref{7.33}, namely 
\begin{eqnarray}\label{9.29}
&&\chi(p^k)=\sum_{(k_1,k_2)\in K_1\times K_2}
\delta(p^{k_1},p^{k_2},p^k)\\
&=& \delta_{(k\in K_1)}\cdot \sum_{k_2\in K_2:\, k_2<k}
\varphi(p^{k_2}) + 
\delta_{(k\in K_2)}\cdot \sum_{k_1\in K_1:\, k_1<k}
\varphi(p^{k_1}) \nonumber\\
&+& \delta_{(k\in K_1\cap K_2\cap\N)}\cdot (p-2)p^{k-1}
+\delta_{(k\in K_1\cap K_2\cap\{0\})} 
+ \sum_{k_1\in K_1\cap K_2:\, k_1>k} \varphi(p^{k_1}).
\nonumber
\end{eqnarray}

We can restrict to the two cases $k^{max}=f(K_1)$
and $k^{max}=e(K_1)$.

The case $k^{max}=f(K_1)\geq 0$: 
We consider $k\in \Z_{[0,k^{max}]}$ and claim 
\begin{eqnarray}\label{9.30}
\chi(p^k)&=& \sum_{k_2\in  K_2:\, k_2\leq k^{max}}
\varphi(p^{k_2})
+ \sum_{k_1\in K_1\cap K_2:\, k_1>k^{max}}\varphi(p^{k_1}).
\end{eqnarray}
We prove this in the two subcases $k\notin K_2$
and $k\in K_2$. In both cases $\Z_{[0,k^{max}]}\subset K_1$,
and $\delta_{(k\in K_1)}=1$, 
and $\sum_{k_1\in K_1:\, k_1<k}\varphi(p^{k_1})
=\delta_{(k>0)}\cdot p^{k-1}$.

The subcase $k\notin K_2$: Then only the first and the last
summand in \eqref{9.29} do not vanish. The last summand 
can be split into
\begin{eqnarray*}
\sum_{k_2\in K_2:\, k<k_2\leq k^{max}}\varphi(p^{k_2})
+ \sum_{k_1\in K_1\cap K_2:\, k_1>k^{max}} \varphi(p^{k_1}).
\end{eqnarray*}
The first summand in \eqref{9.29} and this sum 
give \eqref{9.30}.

The subcase $k\in K_2$: Then \eqref{9.29} takes the 
following shape, where we split again the last summand
into the two pieces above,
\begin{eqnarray*}
\chi(k)&=& \sum_{k_2\in K_2:\, k_2<k}\varphi(p^{k_2})
+ \delta_{(k>0)}\cdot p^{k-1} 
+ \delta_{(k>0)}\cdot (p-2)p^{k-1}\\ 
&+& \delta_{(k=0)}
+ \sum_{k_2\in K_2:\, k<k_2\leq k^{max}}\varphi(p^{k_2}) 
+ \sum_{k_1\in K_1\cap K_2:\, k_1>k^{max}} \varphi(p^{k_1}).
\end{eqnarray*}
In the case $k>0$, the second and third summand give
$\varphi(p^k)$. In the case $k=0$, the fourth summand
$\delta_{(k=0)}$ gives $\varphi(p^k)=1$. In both cases
\eqref{9.30} is true.

The case $k^{max}=e(K_1)$: 
We consider $k\in \Z_{[0,k^{max}]}$ and claim 
\begin{eqnarray}\label{9.31}
\chi(p^k)&=&  \sum_{k_1\in K_1\cap K_2:\, k_1>k^{max}}\varphi(p^{k_1}).
\end{eqnarray}
Here $\Z_{[0,k^{max}]}\cap K_1=\emptyset$ and especially 
$k\notin K_1$. Only the last summand in \eqref{9.29} 
does not vanish. It gives \eqref{9.31}.
This finishes the proof of Claim 1. \hfill ($\Box$)

\medskip
{\bf Claim 2:} For $j\in\{1,2\}$
\begin{eqnarray}\label{9.32}
K_j\cap \Z_{>k^{(j)}} = S(\succ_j)\cap \Z_{>k^{(j)}}.
\end{eqnarray}

\medskip
{\bf Proof of Claim 2:}
If $k^{(j)}=e(K_j)$ then $K_j=K_j\cap\Z_{>k^{(j)}}$ and
especially $0\notin K_j$. 
Then the subset compatibility of
$K_j$ with $\succ_j$ shows 
$K_j= S(\succ_j)\cap \Z_{>k^{(j)}}$. 
Then \eqref{9.32} is clear.

If $k^{(j)}=f(K_j)$ then $K_j\supset \Z_{[0,k^{(j)}]}$
and $k^{(j)}+1\notin K_j$ and $0\in K_j$.
Then the subset compatibility of $K_j$
with $\succ_j$ shows
$K_j=(S(\succ_j)\cap \Z_{>k^{(j)}})\cup \Z_{[0,k^{(j)}]}$.
Also then \eqref{9.32} is clear.
This finishes the proof of Claim 2.\hfill ($\Box$)

\medskip
Write $\succ_3:=(\succ_1\otimes \succ_2)$.
First we will show \eqref{9.11}, then \eqref{9.12}.
By subset compatibility $K_j\subset \Z_{[0,s(\succ_j)]}$.
By definition of $\chi$,
\begin{eqnarray*}
\supp(\chi)&\subset& \{p^k\,|\, k\leq \max(K_1\cup K_2)\}\\
&\subset& \{p^k\,|\, k\leq s(\succ_3)\} =V(\succ_3).
\end{eqnarray*}
This is \eqref{9.11}. For the proof of \eqref{9.12},
we have to show the following:
For $k_1,k_2\in\Z_{[0,s(\succ_3)]}$ with $k_1>k_2$ 
\begin{eqnarray}\label{9.33}
\left. \begin{array}{lll}
k_1\succ_3 k_2&\Rightarrow& \chi(p^{k_1})\geq \chi(p^{k_2}),\\
k_2\succ_3 k_1&\Rightarrow& \chi(p^{k_2})\geq \chi(p^{k_1}).
\end{array}\right\}
\end{eqnarray}
We consider several cases.

The case $k_1\leq k^{max}$: Then Claim 1 implies
$\chi(p^{k_1})=\chi(p^{k_2})$. Then it does not matter
whether $k_1\succ_3 k_2$ or $k_2\succ_3 k_1$. 

The case $k_1>k^{max}$ and $k_1\in S(\succ_3)$:
By the relation between $S(\succ_3)$ and $\succ_3$
then $k_1\succ_3 k_2$. 
By definition of $S(\succ_3)$ in \eqref{9.4}
and by \eqref{9.32} $k_1\in (K_1\cup K_2)-(K_1\cap K_2)$. 
Lemma \ref{t7.6} (a) implies 
$\chi(p^{k_1})\geq \chi(p^{k_2})$. This shows \eqref{9.33}.

The case $k_1>k^{max}$ and $k_1\notin S(\succ_3)$:
By the relation between $S(\succ_3)$ and $\succ_3$
then $k_2\succ_3 k_1$. 
By definition of $S(\succ_3)$ in \eqref{9.4}
and by \eqref{9.32} $k_1\notin (K_1\cup K_2)-(K_1\cap K_2)$. 
Lemma \ref{t7.6} (a) implies 
$\chi(p^{k_2})\geq \chi(p^{k_1})$. This shows \eqref{9.33}.
This finishes the proof of Lemma \ref{t9.8}.
\hfill$\Box$

\begin{theorem}\label{t9.9}
Let $M\subset\N$ and $N\subset \N$ be two finite non-empty sets,
let $P\supset \PP(M)\cup\PP(N)$ be a finite set of prime 
numbers, and let $(\succ_p^M)_{p\in P}$ and
$(\succ_p^N)_{p\in P}$ be two tuples of excellent orders
such that $M$ and $(\succ_p^M)_{p\in P}$
are compatible and $N$ and $(\succ_p^N)_{p\in P}$
are compatible. 

Let $\chi:\N\to\N_0$ be as in \eqref{7.19}
with finite support $L$.
For any $p\in P$ and any $(m_0,n_0)\in\pi_p(M)\times
\pi_p(N)$ define $\chi_{p,m_0,n_0}:L\to\N_0$
as in \eqref{7.21}.
Extend it to $\N-L$ with values 0.
Write $\succ_p^L:=\succ_p^M\otimes\succ_p^N$ for any
$p\in P$. 

\medskip
(a) Any $\chi_{p,m_0,n_0}$ is compatible with
the excellent order $\succ_p^L$ in the following
sense (which was not considered in Definition
\ref{t9.4} (d)):
\begin{eqnarray}\label{9.34}
(l_a,l_b)\in E_{V,p}&\Rightarrow& 
\chi_{p,m_0,n_0}(l_a)\geq \chi_{p,m_0,n_0}(l_b).
\end{eqnarray}
Here $V=V((\succ_q^L)_{q\in P})$, and
$E_{V,p}$ in \eqref{9.6} is determined by $\succ_p^L$.

\medskip
(b) $\chi$ is compatible with the tuple
$(\succ_p^L)_{p\in P}$ of excellent orders.

\medskip
(c) The pair $(M,N)$ is sdiOb-sufficient.
\end{theorem}

{\bf Proof:}
(a) Define $l_0:=\pi_p(l_a)=\pi_p(l_b)$.
If $\delta(m_0,n_0,l_0)=0$ then by Lemma \ref{t7.6} (b)
$\chi_{p,m_0,n_0}(l_a)=\chi_{p,m_0,n_0}(l_b)=0$,
and \eqref{9.34} holds.
Consider the case $\delta(m_0,n_0,l_0)>0$. 
Define $K_1:=K_{M,p,m_0}$ and $K_2:=K_{N,p,n_0}$. 
Then for $k\in\N_0$
\begin{eqnarray}\label{9.35}
(\chi\textup{ in Lemma \ref{t9.8}})(p^k)&=&
\frac{\chi_{p,m_0,n_0}(l_0\cdot p^k)}{\delta(m_0,n_0,l_0)}.
\end{eqnarray}
Therefore Lemma \ref{t9.8} implies \eqref{9.34}.

(b) We have to show for any edge $(l_a,l_b)\in E_V$
$\chi(l_a)\geq \chi(l_b)$. There is a unique prime number
$p$ with $(l_a,l_b)\in E_{V,p}$. 
By \eqref{7.22}
\begin{eqnarray*}
\chi &=& \sum_{(m_0,n_0)\in\pi_p(M)\times \pi_p(N)}
\chi_{p,m_0,n_0}.
\end{eqnarray*}
This and part (a) show $\chi(l_a)\geq \chi(l_b)$.

(c) We have to show for any prime number $p$ and
any $p$-edge $(l_a,l_b)\in E_p(L)$ that \eqref{7.23}
or \eqref{7.24} holds.

Either $(l_a,l_b)\in E_{V,p}$ or $(l_b,l_a)\in E_{V,p}$.
In the case $(l_a,l_b)\in E_{V,p}$, \eqref{9.34} gives
\eqref{7.23}.
In the case $(l_b,l_a)\in E_{V,p}$, \eqref{9.34} gives
\eqref{7.24}. \hfill $\Box$

\begin{theorem}\label{t9.10}
Let $\chi_1:\N\to\N_0$ and $\chi_2:\N\to\N_0$
be two maps with finite supports $M=\supp(\chi_1)$
and $N=\supp(\chi_2)$. 
Define a map $\chi_3:\N\to\N_0$ by
\begin{eqnarray}\label{9.36}
\Bigl(\prod_{m\in M}\Phi_m^{\chi_1(m)}\Bigr)
\otimes \Bigl(\prod_{n\in N}\Phi_n^{\chi_2(n)}\Bigr)
= \prod_{l\in \N_0}\Phi_l^{\chi_3(l)}.
\end{eqnarray}
$\chi_3$ has finite support $\supp(\chi_3)=:L$. 

Let $P\supset\PP(M)\cup\PP(N)$ be a finite set of
prime numbers. Let $(\succ_p^M)_{p\in P}$
and $(\succ_p^N)_{p\in P}$ be two tuples of excellent
orders such that $\chi_1$ is compatible with
$(\succ_p^M)_{p\in P}$ and $\chi_2$ is compatible
with $(\succ_p^N)_{p\in P}$. 
Write $\succ_p^L:=(\succ_p^M\otimes \succ_p^N)$
for any $p\in P$. 

\medskip
(a) $\chi_3$ is compatible with the tuple
$(\succ_p^L)_{p\in P}$ of excellent orders.

\medskip
(b) Let $(M_1^{(st)},...,M_{l_{\chi_1}}^{(st)})$,
$(N_1^{(st)},...,N_{l_{\chi_2}}^{(st)})$ and
$(L_1^{(st)},...,L_{l_{\chi_3}}^{(st)})$ be
the standard coverings of $\chi_1$, $\chi_2$
and $\chi_3$. Then
\begin{eqnarray}\label{9.37}
\Bigl(\bigoplus_{i=1}^{l_{\chi_1}}\Or(M_i^{(st)})\Bigr)
\otimes
\Bigl(\bigoplus_{j=1}^{l_{\chi_2}}\Or(N_j^{(st)})\Bigr)
\cong \bigoplus_{k=1}^{l_{\chi_3}}\Or(L_k^{(st)}),
\end{eqnarray}
so the tensor product of sums of Orlik blocks on the left
hand side admits a standard decomposition into Orlik blocks.
\end{theorem}

{\bf Proof:}
(a) For $(i,j)\in\Z_{[0,l_{\chi_1}]}\times \Z_{[0,l_{\chi_2}]}$
define a map $\chi_{3,i,j}:\N\to \N_0$ by
\begin{eqnarray}\label{9.38}
\Bigl(\prod_{m\in M_i^{(st)}}\Phi_m^{\chi_1(m)}\Bigr)
\otimes \Bigl(\prod_{n\in N_j^{(st)}}\Phi_n^{\chi_2(n)}\Bigr)
= \prod_{l\in \N_0}\Phi_l^{\chi_{3,i,j}(l)}.
\end{eqnarray}
$\chi_{3,i,j}$ has finite support.

By Lemma \ref{t9.5} (b) and the hypotheses on
$\chi_1$ and $\chi_2$, the set $M_i^{(st)}$
is compatible with $(\succ_p)^M_{p\in P}$,
and the set $N_j^{(st)}$ is compatible with 
$(\succ_p)^N_{p\in P}$. 

By Theorem \ref{t9.9} (b), $\chi_{3,i,j}$ is 
compatible with $(\succ_p)^L_{p\in P}$.

Because of \eqref{9.38}, 
$\chi_3=\sum_{(i,j)}\chi_{3,i,j}$.
Therefore also $\chi_3$ is compatible with
$(\succ_p)^L_{p\in P}$.

(b) Let $(L_1^{(i,j,st)},...,L_{l_{\chi_{3,i,j}}}^{(i,j,st)})$
be the standard covering of $\chi_{3,i,j}$. 
By Lemma \ref{t9.5} (b) each set 
$L_k^{(i,j,st)}$ is compatible with $(\succ_p)^L_{p\in P}$.

By Theorem \ref{t9.9} (c), each pair
$(M_i^{(st)},N_j^{(st)})$ is sdiOb-sufficient.

By Theorem \ref{t7.4}, the tensor product 
$\Or(M_i^{(st)})\otimes \Or(N_j^{(st)})$ admits
a standard decomposition into Orlik blocks.
In other words,
\begin{eqnarray}\label{9.39}
\Or(M_i^{(st)})\otimes \Or(N_j^{(st)})
&\cong& \bigoplus_{k=1}^{l_{\chi_{3,i,j}}}\Or(L_k^{(i,j,st)}).
\end{eqnarray}

Because of $\chi_3=\sum_{(i,j)}\chi_{3,i,j}$,
the tuple 
\begin{eqnarray*}
(L_k^{(i,j,st)}\,|\, (i,j)\in 
\Z_{[0,l_{\chi_1}]}\times \Z_{[0,l_{\chi_2}]},\ 
k\in \Z_{[0,l_{\chi_{3,i,j}}]})
\end{eqnarray*}
is a covering of $\chi_3$. Because all sets in it are
compatible with $(\succ_p)_{p\in P}^L$, it is also
compatible with $(\succ_p)_{p\in P}^L$. 

\eqref{9.39} gives the (possibly) non standard 
decomposition on the right hand side of \eqref{9.40}
for the tensor product of sums of Orlik blocks on
the left hand side of \eqref{9.40},
\begin{eqnarray}\label{9.40}
\Bigl(\bigoplus_{i=1}^{l_{\chi_1}}\Or(M_i^{(st)})\Bigr)
\otimes
\Bigl(\bigoplus_{j=1}^{l_{\chi_2}}\Or(N_j^{(st)})\Bigr)
= \bigoplus_{(i,j,k)}\Or(L_k^{(i,j,st)}).
\end{eqnarray}
We can apply Theorem \ref{t9.6}.
It says that the sum of Orlik blocks on the right hand
side of \eqref{9.40} is isomorphic to the standard
decomposition into Orlik blocks on the right hand
side of \eqref{9.37}.
\hfill $\Box$

\begin{remarks}\label{t9.11}
(i) In Theorem \ref{t9.9}, the sdiOb-sufficiency of
the pair $(M,N)$ in part (c) is weaker than part (a). 
The sdiOb-sufficiency demands only that for any
fixed $p$-edge $(l_a,l_b)\in E_p(M)$ one has
\eqref{7.23} or \eqref{7.24}.
Part (a) gives for fixed $k_a$ and $k_b$
and any $(l_a,l_b)\in E_p(M)$ with $v_p(l_a)=k_a$
and $v_p(l_b)=k_b$ the same alternative
\eqref{7.23} or \eqref{7.24}.

\medskip
(ii) Therefore one might ask whether a weaker condition
than the compatibility with tuples of excellent orders
might also have good properties. 
Lemma \ref{t9.12} says in a precise sense that this is 
not the case, but that the compatibility with tuples of
excellent orders is needed.
Lemma \ref{t9.12} was the way how we found the
excellent orders and the compatibilities with them.

\medskip
(iii) Condition (ii) in Lemma \ref{t9.12} is 
via Theorem \ref{t7.4} and the (conjectural) Remark
\ref{t7.5} (i) the condition
that the tensor product of an Orlik block $\Or(M)$
with the Milnor lattice of an arbitrary $A_\mu$-singularity
admits a standard decomposition into Orlik blocks.
So, it is a quite natural condition.
\end{remarks}

\begin{lemma}\label{t9.12}
Let $M\subset\N$ be a finite non-empty set.
The following conditions are equivalent.
\begin{list}{}{}
\item[(i)]
For any $n_N\in\N$, the pair
$(M,\{n\in\N\,|\, n|n_N\})$ is sdiOb-sufficient.
\item[(ii)]
For any $n_N\in\N$, the pair
$(M,\{n\in\N\,|\, n|n_N\}-\{1\})$ is sdiOb-sufficient.
\item[(iii)]
A tuple $(\succ_p)_{p\in \PP(M)}$ of excellent orders
exists such that $M$ is compatible with it.
\end{list}
\end{lemma}

{\bf Proof:}
For any $n_N\in\N$, the sets $\{n\in\N\,|\, n|n_N\}$
and $\{n\in\N\,|\, n|n_N\}-\{1\})$ are compatible
with suitable tuples of excellent orders.
Therefore and by Theorem \ref{t9.9} (c),
(iii) implies (i) and (ii). 

(i)$\Rightarrow$(iii): 
Suppose that (i) holds. We will define for any prime number
$p\in\PP(M)$ an excellent order $\succ_p$ and then show
that the set $M$ is compatible with the tuple 
$(\succ_p)_{p\in\PP(M)}$ of excellent orders
(Definition \ref{t9.4} (c)).

Fix a prime number $p\in\PP(M)$. 
We will define an excellent order $\succ_p$
by fixing the number $s(\succ_p)\in\N_0$ and
the set $S(\succ_p)\subset\Z_{[0,s(\succ_p)]}$.  
Then we have to show that for any $m_0\in \pi_p(M)$
the set $K_{M,p,m_0}$ is subset compatible
with $\succ_p$ (Definition \ref{t9.4} (a)). 

Define the map $g_p:\pi_p(M)\to\N_0$ by
\begin{eqnarray}\label{9.41}
g_p(m_0):=\left\{\begin{array}{ll}
\max(k\in\N_0\,|\, \Z_{[0,k]}\subset K_{M,p,m_0})
&\textup{ if }0\in K_{M,p,m_0},\\
\max(k\in\N_0\,|\, \Z_{[0,k]}\cap K_{M,p,m_0}=\emptyset)
&\textup{ if }0\notin K_{M,p,m_0}.
\end{array}\right.
\end{eqnarray}
Define $g_{p,min}:=\min(g_p(m_0)\,|\, m_0\in\pi_p(M)).$
If the set 
\begin{eqnarray}\label{9.42}
\{m_0\in\pi_p(M)\,|\, g_p(m_0)=g_{p,min},
0\notin K_{M,p,m_0}\}
\end{eqnarray}
is not empty, choose an element of it and denote it
by $\www{m_0}$. If the set in \eqref{9.42} is empty,
the set
\begin{eqnarray}\label{9.43}
\{m_0\in\pi_p(M)\,|\, g_p(m_0)=g_{p,min},
0\notin K_{M,p,m_0}\}
\end{eqnarray}
is not empty. Then choose an element of if and denote it
by $\www{m_0}$. Define an excellent order $\succ_p$ by
\begin{eqnarray}
s(\succ_p)&:=& \max(v_p(m)\,|\,m\in M),\nonumber\\
S(\succ_p)&:=&\left\{\begin{array}{ll}
K_{M,p,\www{m_0}}&\textup{ if }0\notin K_{M,p,\www{m_0}},\\
K_{M,p,\www{m_0}}-\{0\}&\textup{ if }0\in K_{M,p,\www{m_0}}.
\end{array}\right. \label{9.44}
\end{eqnarray}
The inclusion  $S(\succ_p)\subset\Z_{[0,s(\succ_p)]}$
is obvious from the definition of $s(\succ_p)$. 

We will show for each $m_0\in \pi_p(M)$ 
\begin{eqnarray}\label{9.45} 
K_{M,p,m_0}\cap\Z_{>g_p(m_0)}
= K_{M,p,\www{m_0}}\cap \Z_{>g_p(m_0)}.
\end{eqnarray}
Because of $g_p(m_0)\geq g_{p,min}=g_p(\www{m_0})$,
this is equivalent to $K_{M,p,m_0}$ being subset compatible
with $\succ_p$. It is sufficient to show \eqref{9.45}
for $m_0\in\pi_p(M)-\{\www{m_0}\}$. 

We will need for this the assumption that (i) holds. 
We take that into account in the following way.
Choose any $k_N\in\N_0$ and define
\begin{eqnarray}\label{9.46}
l_0:=n_0&:=&\prod_{q\in\PP(M)-\{p\}}
q^{1+\max(v_q(m)\,|\, m\in M)},\\
n_N&:=& n_0\cdot p^{k_N},\nonumber\\
N&:=& \{n\in\N\,|\, n|n_N\}.\nonumber
\end{eqnarray}
Then $v_p(n_N)=k_N$, $\pi_p(n_N)=n_0$, and 
\begin{eqnarray}\label{9.47}
\delta(m_0,n_0,l_0)>0\quad\textup{for any }
m_0\in\pi_p(M).
\end{eqnarray}
A priori $\chi_{p,m_0,n_0}/\delta(m_0,n_0,l_0)$
for any $m_0\in \pi_p(M)$ is as in \eqref{7.33}.
But as in the proof of Lemma \ref{t9.8}, 
it boils down to the following much simpler form, 
because of $K_{N,p,n_0}=\Z_{[0,k_N]}$:
\begin{eqnarray}\label{9.48}
\frac{\chi_{p,m_0,n_0}(l_0\cdot p^k)}
{\delta(m_0,n_0,l_0)}
= \left\{\begin{array}{l}
p^{k_N}\qquad\textup{ if }k>k_N\textup{ and }
k\in K_{M,p,m_0},\\
0\ \qquad\textup{ if }k>k_N\textup{ and }
k\notin K_{M,p,m_0},\\
\sum_{k_1\in K_{M,p,m_0}:\, k_1\leq k_N}
\varphi(p^{k_1})\ \textup{ if }k\leq k_N.
\end{array}\right. 
\end{eqnarray} 
The last sum is in $\Z_{[0,p^{k_N}]}$. 
Thus the quotient in \eqref{9.48} has only 2 or 3 values,
for fixed $m_0$ and varying $k$.  

Now consider any $m_0\in\pi_p(M)-\{\www{m_0}\}$.
We want to prove \eqref{9.45}. 

{\bf 1st case:} 
Suppose $k\in K_{M,p,\www{m_0}}\cap\Z_{>g_p(m_0)}$
and $k\notin K_{M,p,m_0}$. We want to arrive at a 
contradiction. 
We choose $k_N:=k-1$. 
We have $K_{M,p,\www{m_0}}\not\supset\Z_{[0,k-1]}$,
because $g_p(\www{m_0})\leq k-1$. \eqref{9.48}
for $\www{m_0}$ gives 
\begin{eqnarray}\label{9.49}
0\leq\frac{\chi_{p,\www{m_0},n_0}(l_0\cdot p^{k-1})}
{\delta(\www{m_0},n_0,l_0)}
< \frac{\chi_{p,\www{m_0},n_0}(l_0\cdot p^k)}
{\delta(\www{m_0},n_0,l_0)} = p^{k-1}.
\end{eqnarray}
We have $K_{M,p,m_0}\cap\Z_{[0,k-1]}\neq\emptyset$,
because $g_p(m_0)\leq k-1$. 
\eqref{9.48} for $m_0$ gives 
\begin{eqnarray}\label{9.50}
p^{k-1}\geq \frac{\chi_{p,m_0,n_0}(l_0\cdot p^{k-1})}
{\delta(m_0,n_0,l_0)}> 
\frac{\chi_{p,m_0,n_0}(l_0\cdot p^k)}
{\delta(m_0,n_0,l_0)}=0.
\end{eqnarray}
The strict inequalities in \eqref{9.49} and \eqref{9.50}
contradict the sdiOb-sufficiency of the pair
$(M,N)$ for the $p$-edge $(l_a,l_b)=(k,k-1)$
(Definition \ref{t7.3} (d)). 

{\bf 2nd case:} Suppose $k\in K_{M,p,m_0}\cap\Z_{>g_p(m_0)}$
and $k\notin K_{M,p,\www{m_0}}$. 
We exchange the roles of $\www{m_0}$ and $m_0$ in the
1st case, and we arrive in exactly the same way
at a contradiction.

We have proved \eqref{9.45}. 

(ii)$\Rightarrow$(iii): This is similar to the proof
of (i)$\Rightarrow$(iii).
\hfill $\Box$

\section{Chain type singularities}\label{c10}
\setcounter{equation}{0}

\noindent
A chain type singularity is a quasihomogeneous singularity of the
special shape
\begin{eqnarray}\label{10.1}
f=f(x_1,...,x_n)= x_1^{a_1+1}+\sum_{i=2}^nx_{i-1}x_i^{a_i}
\end{eqnarray}
for some $n\in\N$ and some $a_1,...,a_n\in\N$.
This quasihomogeneous polynomial has an isolated singularity. Define 
\begin{eqnarray}\label{10.2}
b_k:= (a_1+1)\cdot a_2\cdot ...\cdot a_k\quad\textup{for }
k\in\{1,...,n\},\quad b_0:=1.
\end{eqnarray}
The Milnor number $\mu$, the weights and the characteristic polynomial
are calculated for example in \cite[Corollary 4.3]{HZ19}
(applying formulas in \cite{MO70}).
Here we need the following quite surprising result of Orlik and Randell.
The function $\chi:\N\to \{0,1,...,n+1\}$ with 
\begin{eqnarray}\label{10.3}
\chi(m)&:=& \left\{\begin{array}{ll}
n+1 & \textup{if }m\not| b_n,\\
\min(i\in\{0,1,...,n\}\, |\, m|b_i)&\textup{ if }m|b_n.
\end{array}\right. 
\end{eqnarray}
will be useful.

\begin{theorem}\label{t10.1} \cite[Theorem (2.11)]{OR77}
For any chain type singularity $f$ as in \eqref{10.1},
an automorphism $h:H_{Mil}\to H_{Mil}$ exists such that $(H_{Mil},h)$ is
an Orlik block, $h_{Mil}=h^\mu$, and the set $M$ of orders of the eigenvalues
of $h$ is as follows,
\begin{eqnarray}\label{10.4}
M= \{m\in \N\, |\, \chi(m)\equiv n\modd 2\}
\subset\{m\in\N\,|\, m|b_n\}.
\end{eqnarray}
\end{theorem}

Theorem \ref{t1.3} (a) says that the pair $(H_{Mil},h_{Mil})$ for each
chain type singularity admits a standard Orlik decomposition. 
Here we give the proof. It is an easy application of Theorem \ref{t10.1}
and Theorem \ref{t6.2}.

\bigskip
{\bf Proof of Theorem \ref{t1.3} (a):} 
We will show that the pair $(M,\mu)$ is sdiOb-sufficient.
This and Theorem \ref{t6.2} imply that $(H_{Mil},h_{Mil})$ 
admits a standard Orlik decomposition.


Consider a prime number $p$ and a $p$-edge $(n_a,n_b)\in E_p(\gamma_\mu(M))$.
Then because of \eqref{6.5}
\begin{eqnarray}\label{10.5}
\gamma_\mu^{-1}(n_a)&=& 
\{m_a^0\cdot c\, |\, c\textup{ divides }\prod_{q\in\PP:v_q(n_a)=0}q^{v_q(\mu)}\},\\
\textup{where } m_a^0&:=& n\cdot \prod_{q\in\PP:\, v_q(n_a)>0}q^{v_q(\mu)}.
\label{10.6}
\end{eqnarray}
We want to prove that the $p$-edge $(n_a,n_b)$ satisfies \eqref{6.3}
if $m_a^0\notin M$ and that it satisfies \eqref{6.4} if $m_a^0\in M$.

But before we make an observation which is valid in both cases.

\medskip
{\it Observation 11.2:
If $(m_c,m_d)\in\gamma_\mu^{-1}(n_a)\times \gamma_\mu^{-1}(n_b)$
is a $p$-edge and $\chi(m_c)>\chi(m_a^0)$ then $\chi(m_d)=\chi(m_c)$.}

\medskip
Proof of Observation 11.2: As $(m_c,m_d)$ is a $p$-edge,
$v_q(m_c)=v_q(m_d)$ for any prime number $q\neq p$ and $v_p(m_c)>v_p(m_d)$.
The number $\chi(m_c)$ is characterized by $v_q(m_c)\leq v_q(b_{\chi(m_c)})$
for any prime number $q$ and $v_r(m_c)>v_r(b_{\chi(m_c)-1})$ for some
prime number $r$. Here $r\neq p$ follows from 
$v_r(b_{\chi(m_c)-1})\geq v_r(b_{\chi(m_a^0)})\geq v_r(m_a^0)$ 
and $v_p(m_c)=v_p(m_a^0)$ (which follows from \eqref{10.5}).
This shows $v_q(m_d)\leq v_q(b_{\chi(m_c)})$ for any prime number $q$
and $v_r(m_d)=v_r(m_c)> v_r(b_{\chi(m_c)-1})$ for the given prime number $r$.
This implies $\chi(m_d)=\chi(m_c)$. \hfill$(\Box)$

\medskip
First suppose $m_a^0\notin M$. Then any $m_c\in\gamma_\mu^{-1}(n_a)\cap M$
satisfies $\chi(m_c)>\chi(m_a^0)$ and $\chi(m_c)\equiv n\modd 2$. 
By the Observation 9.2, any $m_d\in \gamma_\mu^{-1}(n_b)$
with $(m_c,m_d)$ a $p$-edge satisfies $\chi(m_d)=\chi(m_c)\equiv n\modd 2$,
thus $m_d\in M$. This shows \eqref{6.3} for the $p$-edge $(n_a,n_b)$.

Second suppose $m_a^0\in M$. Consider a $p$-edge
$(m_c,m_d)\in \gamma^{-1}(n_a)\times (\gamma^{-1}(n_b)\cap M)$.
If $\chi(m_c)>\chi(m_a^0)$ then $\chi(m_c)=\chi(m_d)$ by the Observation 9.2,
and thus $\chi(m_c)=\chi(m_d)\equiv n\modd 2$, so $m_c\in M$. 
If $\chi(m_c)=\chi(m_a^0)$ then $\chi(m_c)=\chi(m_a^0)\equiv n\modd 2$,
and again $m_c\in M$. So in both cases $m_c\in M$. This
shows \eqref{6.4} for the $p$-edge $(n_a,n_b)$.
\hfill$\Box$ 

\setcounter{lemma}{2}

\begin{remark}\label{t10.2}
The proof gives that $(M,\www\mu)$ is sdiOB-sufficient 
for any $\www\mu\in \N$.
The proof does not use the following formula \eqref{10.7} 
for the Milnor number:
\begin{eqnarray}\label{10.7}
\mu=\sum_{i=0}^n (-1)^i\cdot b_{n-i}=b_n-b_{n-1}+...+(-1)^{n-1}b_1+(-1)^n.
\end{eqnarray}

\end{remark}

\section{Cycle type singularities}\label{c11}
\setcounter{equation}{0}

\noindent
A cycle type singularity is a quasihomogeneous singularity
of the special shape
\begin{eqnarray}\label{11.1}
f=f(x_1,...,x_n)=\sum_{i=1}^{n-1}x_i^{a_i}x_{i+1}+x_n^{a_n}x_1
\end{eqnarray}
for some $n\in\Z_{\geq 2}$ and some $a_1,...,a_n\in\N$ which satisfy 
\begin{eqnarray}\label{11.2}
\textup{for even }n\textup{ neither }a_j=1\textup{ for all even }j
\textup{ nor }a_j=1\textup{ for all odd }j.
\end{eqnarray}
The following well known facts are proved for example in
\cite[Lemma 4.1]{HZ19}.
This polynomial is quasihomogeneous and has an isolated 
singularity. The Milnor number is $\mu=\prod_{i=1}^n a_i$.
The weights have the shape $(w_1,...,w_n)=(\frac{v_1}{d},...,
\frac{v_n}{d})$ with $d=\mu-(-1)^n$ and $v_1,...,v_n\in\N$
(for even $n$ $v_1,...,v_n>0$ requires \eqref{11.2}).
These numbers satisfy $\gcd(v_1,d)=...=\gcd(v_n,d)$.
Define $b:=d/\gcd(v_1,d)\in\N$. Then
\begin{eqnarray}\label{11.3}
p_{H_{Mil},h_{Mil}}=(t^b-1)^{\gcd(v_1,d)}\cdot (t-1)^{(-1)^n}.
\end{eqnarray}
Therefore Orlik's conjecture says here the following.
\begin{eqnarray}\label{11.4}
(H_{Mil},h_{Mil})\cong
\left\{\begin{array}{ll}
(\gcd(v_1,d)-1) H^{[t^b-1)]}& \\
\hspace*{1cm}\oplus H^{[(t^b-1)/(t-1)]} 
& \textup{if }n\textup{ is odd,}\\
\gcd(v_1,d)H^{[t^b-1)]}\oplus H^{[t-1]} 
& \textup{if }n\textup{ is even.}
\end{array}\right.
\end{eqnarray}
It is true by Theorem \ref{t1.3} (b). 
We proved Theorem \ref{t1.3} (b) in \cite[Theorem 1.3]{HM20-1},
using algebraic topology and a spectral sequence and 
building on \cite{Co82}. Cooper's paper had the same aim.
But it contains two serious mistakes. The second one leads
for even $n$ to the (wrong) claim in \cite{Co82} 
that $(H_{Mil},h_{Mil})$ has a decomposition into Orlik blocks, 
but not a standard decomposition into Orlik blocks. 
See the introduction and the Remarks 5.1 in \cite{HM20-1} 
for the relation of \cite{HM20-1} to \cite{Co82}.

\section{Quasihomogeneous singularities and their 
Thom-Sebastiani sums}\label{c12}
\setcounter{equation}{0}

\noindent
Theorem \ref{t12.1} says that the exponent map
$\chi_{Mil}^f$ of the characteristic polynomial
$p_{Mil}^f$ of an isolated quasihomogeneous singularity $f$ 
comes equipped with a canonical compatible tuple
of excellent orders. 
Together with Theorem \ref{t9.10}, it implies Theorem 
\ref{t1.3} (c), namely that the Thom-Sebastiani
sum of two singularities satisfies Orlik's conjecture
if the two singularities satisfy Orlik's conjecture.

\begin{theorem}\label{t12.1}
Consider an isolated quasihomogeneous singularity
$f(x_1,...,x_n)$ with weight system 
$(w_1,...,w_n)\in (\Q\cap(0,1))^n$ with
$w_j=\frac{s_j}{t_j}$ and $s_j,t_j\in\N$ with $\gcd(s_j,t_j)=1$. 
The characteristic polynomial of the monodromy $h_{Mil}^f$
on its Milnor lattice $H_{Mil}^f$ is called here
$p_{Mil}^f$. 
It gives rise to an exponent map $\chi_{Mil}^f:\N\to\N_0$
with finite support $M_f:=\supp(\chi_{Mil}^f)$ by 
$p_{Mil}^f=\prod_{m\in\N}\Phi_m^{\chi_{Mil}^f(m)}$. 
The map $\chi_{Mil}^f$ is compatible with the tuple
$(\succ_p^f)_{p\in P(M_f)}$ of excellent orders which is 
defined as follows,
\begin{eqnarray}\label{12.1}
s(\succ_p^f)&:=& \max(v_p(m)\,|\, m\in M_f),\\
S(\succ_p^f)&:=& \{k\in \Z_{[0,s(\succ_p)]}\,|\, 
|\{j\in\{1,...,n\}\,|\, p^k|t_j\}|\textup{ is odd}\}.
\label{12.2}
\end{eqnarray}
\end{theorem}

Before the proof, we make some remarks, show how Theorem
\ref{t12.1} implies Theorem \ref{t1.3} (c), and
cite two classical results Theorem \ref{t12.4} and 
Theorem \ref{t12.5}, which will be needed in the proof
of Theorem \ref{t12.1}.

\begin{remark}\label{t12.2}
Consider two isolatd quasihomogeneous singularities
$f(x_1,...,x_{n_f})$ and $g(x_{n_f+1},...,x_{n_f+n_g})$.
Denote the tuples of excellent orders of $f$, $g$ and $f+g$
by $(\succ_p^f)_{p\in M_f}$, $(\succ_p^g)_{p\in M_g}$
and $(\succ_p^{f+g})_{p\in M_{f+g}}$. 
Then $M_{f+g}\subset M_f\cup M_g$. Extend them to tuples
of excellent orders 
$(\succ_p^f)_{p\in M_f\cup M_g}$, $(\succ_p^g)_{p\in M_f\cup M_g}$
and $(\succ_p^{f+g})_{p\in M_f\cup M_g}$ by copies of the 
trivial excellent order $\succ_0$ (Definition \ref{t9.1} (b)).
Then \eqref{12.2} and the definition of the tensor product 
of two excellent orders (Definition \ref{t9.1} (c)) 
show immediately
\begin{eqnarray}\label{12.3}
\succ_p^{f+g}&=& \succ_p^f\otimes \succ_p^g\quad\textup{for }
p\in P(M_f\cup M_g).
\end{eqnarray}
\end{remark}

{\bf Proof of Theorem \ref{t1.3} (c):}
Consider the data in Remark \ref{t12.2}.
Identify them as follows with the data in Theorem \ref{t9.10}:
\begin{eqnarray*}
\chi_{Mil}^f=\chi_1,\quad \succ_p^f=\succ_p^M,\quad 
\chi_{Mil}^g=\chi_2,\quad \succ_p^g=\succ_p^N.
\end{eqnarray*}
The basic result
\begin{eqnarray}\label{12.4}
(H_{Mil}^{f+g},h_{Mil}^{f+g})\cong (H_{Mil}^f,h_{Mil}^f)
\otimes (H_{Mil}^g,h_{Mil}^g)
\end{eqnarray}
of Sebastiani and Thom \cite{ST71} implies 
$\chi_{Mil}^{f+g}=\chi_3$. And Remark \ref{t12.2} implies
$\succ_p^{f+g}=\succ_3$.

Let $(M_1^{(st)},...,M_{l_{\chi_1}}^{(st)})$,
$(N_1^{(st)},...,N_{l_{\chi_2}}^{(st)})$ and
$(L_1^{(st)},...,L_{l_{\chi_3}}^{(st)})$ be
the standard coverings of $\chi_1$, $\chi_2$
and $\chi_3$. The assumption that $f$ and $g$ satisfy 
Orlik's conjecture says
\begin{eqnarray}\label{12.5}
(H_{Mil}^f,h_{Mil}^f)\cong 
\bigoplus_{i=1}^{l_{\chi_1}}\Or(M_i^{(st)}),\quad  
(H_{Mil}^g,h_{Mil}^g)\cong 
\bigoplus_{j=1}^{l_{\chi_2}}\Or(N_j^{(st)}).
\end{eqnarray}
Theorem \ref{t9.10} applies because of Theorem \ref{t12.1}.
Together \eqref{12.4}, \eqref{12.5} and \eqref{9.37} in
Theorem \ref{t9.10} (b) give
\begin{eqnarray}\label{12.6}
(H_{Mil}^{f+g},h_{Mil}^{f+g})\cong 
\bigoplus_{k=1}^{l_{\chi_3}}\Or(L_k^{(st)}),
\end{eqnarray}
which is Orlik's conjecture for $f+g$.\hfill$\Box$

\begin{remark}\label{t12.3}
Consider an isolated quasihomogeneous singularity $f$
as in Theorem \ref{t12.1}. Write $\chi:=\chi_{Mil}^f$.
Let $(M_1^{(st)},...,M_{l_\chi}^{(st)}$ be the standard
covering of $\chi$. We can show with some extra work 
which we will carry out in \cite{HM20-2} the following:
The compatibility of $\chi$ with the tuple
$(\succ_p^f)_{p\in P(M_f)}$ implies that each 
set $M_k^{(st)}$ satisfies condition (I) in Theorem 1.2
in \cite{He20}, which is also Theorem 6.2 in \cite{HZ19}.
Therefore Theorem \ref{t12.1} and this implication
solve the problems 6 and 7 in \cite{HZ19}. And therefore 
\begin{eqnarray}\label{12.7}
\Aut_{S^1}(\Or(M_k^{(st)}))=\{\pm h_{[p_k]}^j\,|\, j\in\Z\},
\end{eqnarray}
where $p_k:=\prod_{m\in M_k^{(st)}}\Phi_m$, 
$\Or(M_k^{(st)})=(H^{[p_k]},h_{[p_k]})$ (see Definition  
\ref{t2.6}) and 
$\Aut_{S^1}(\Or(M_k^{(st)}))$ denotes the automorphisms of the
Orlik block $\Or(M_k^{(st)})$, whose eigenvalues are in $S^1$.
We will discuss this in \cite{HM20-2}. 
\end{remark}

Milnor and Orlik \cite{MO70} 
proved a formula for the characteristic
polynomial of an isolated quasihomogeneous singularity.
Recall the notations in Definition \ref{t7.1}.

\begin{theorem}\label{t12.4}\cite{MO70}
Consider an isolated quasihomogeneous singularity
$f(x_1,...,x_n)$ with weight system 
$(w_1,...,w_n)\in (\Q\cap(0,1))^n$ with
$w_j=\frac{s_j}{t_j}$ and $s_j,t_j\in\N$ with $\gcd(s_j,t_j)=1$,
and with characteristic polynomial  
$p_{Mil}^f=\prod_{m\in\N}\Phi_m^{\chi_{Mil}^f(m)}$, where
$\chi_{Mil}^f:\N\to\N_0$ has finite support 
$M_f:=\supp{\chi_{Mil}^f}$. The divisor
$\divv p_{Mil}^f=\sum_{m\in M_f}\chi_{Mil}^f(m)\cdot\Psi_m$
is determined by the weights via the following formula,
\begin{eqnarray}\label{12.8}
\divv p_{Mil}^f = \prod_{j=1}^n\bigl(\frac{1}{s_j}\Lambda_{t_j}
-\Lambda_1\bigr).
\end{eqnarray} 
\end{theorem}

Kouchnirenko \cite{Ko76} gave a characterization of the weight
systems which allow quasihomogeneous polynomials with an
isolated singularity at 0. Roughly these are the weight systems
which allow sufficiently many monomials of weighted degree 1.
His result was rediscovered and generalized. See \cite{HK12}
and \cite{HM20-2} for references.

\begin{theorem}\label{t12.5}\cite[Remarque 1.13 (ii)]{Ko76}
Let a weight system $(w_1,...,w_n)\in (\Q\cap(0,1))^n$
be given. For $J\subset \{1,...,n\}$ and $q\in\Q_{\geq 0}$,
denote
\begin{eqnarray}\label{12.9} 
(\N_0^J)_q:= \{(\alpha_1,...,\alpha_n)\in\N_0^n\,|\, 
\alpha_j=0\textup{ for }j\notin J,\ \sum_{j\in J}w_j\alpha_j=q\}.
\end{eqnarray}
The set of quasihomogeneous polynomials of weighted degree 1
with an isolated singularity at 0 is not empty if and only if 
the weight system satisfies the following condition (C2).
\begin{eqnarray}\label{12.10}
(C2)&& \forall\ J\subset\{1,...,n\}\textup{ with }J\neq\emptyset
\quad\exists K\subset\{1,...,n\} \nonumber\\
&& \textup{with }|K|=|J|\textup{ and }\forall\ k\in K\ 
(\N_0^J)_{1-w_k}\neq\emptyset.
\end{eqnarray}
\end{theorem}

{\bf Proof of Theorem \ref{t12.1}:}
Consider the data in Theorem \ref{t12.1}.  
Write $\chi:=\chi_{Mil}^f$. 
Fix a prime number $p\in \PP(M_f)$ and two numbers
$m_a,m_b\in M_f$ with $(m_a,m_b)\in E_p(M_f)$. 
Because of Definition \ref{t9.4} (c), we have to show
\begin{eqnarray}\label{12.11}
\chi(m_a)\geq \chi(m_b)&&\textup{if }v_p(m_a)\succ_p^f v_p(m_b),\\
\chi(m_a)\leq \chi(m_b)&&\textup{if }v_p(m_b)\succ_p^f v_p(m_a).
\nonumber
\end{eqnarray}
Because of the definition of the excellent order $\succ_p^f$, 
this is equivalent to the claim (which has to be proved):
\begin{eqnarray}\label{12.12}
\chi(m_a)\geq \chi(m_b)&&\textup{if }|\{j\in\{1,...,n\}\,|\, 
p^{v_p(m_a)}|t_j\}|\textup{ is odd,}\\
\chi(m_a)\leq \chi(m_b)&&\textup{if }
|\{j\in\{1,...,n\}\,|\, p^{v_p(m_a)}|t_j\}|\textup{ is even.}
\nonumber
\end{eqnarray}

Define the map $\nu:\N\to\Z$ by
\begin{eqnarray}\label{12.13}
\nu(k):=\sum_{m:\, k|m}\chi(m)\cdot\mu_{Moeb}(\frac{m}{k}),
\end{eqnarray}
where $\mu_{Moeb}:\N\to\{0,1,-1\}$ is the Moebius function with
\begin{eqnarray}\label{12.14}
\mu_{Moeb}(m):=\left\{\begin{array}{ll}
(-1)^r&\textup{if }m=p_1\cdot ...\cdot p_r\textup{ with}\\
& p_1,...,p_r\textup{ different prime numbers,}\\
0 & \textup{else.}\end{array}\right.
\end{eqnarray}
It has finite support and satisfies 
\begin{eqnarray}\label{12.15}
\divv p_{Mil}^f=\sum_{k\in\N}\nu(k)\cdot\Lambda_k,\quad
\chi(m)=\sum_{k:\, m|k}\nu(k).
\end{eqnarray}
Thus
\begin{eqnarray}\label{12.16}
\chi(m_b)-\chi(m_a)=\sum_{k:\, m_b|k,\, m_a\nmid k}\nu(k)
=\sum_{k:\, m_b|k,\, p^{v_p(m_a)}\nmid k}\nu(k).
\end{eqnarray}
Formula \eqref{12.8} allows a good control on the map $\nu$. 
Suppose that the weights $(w_1,...,w_n)$ are numbered such that
\begin{eqnarray}\label{12.17}
\{j\in\{1,...,n\}\,|\, p^{v_p(m_a)}\nmid t_j\}
=\{1,...,\www n\}\quad\textup{for some }\www n\leq n.
\end{eqnarray}
Formula \eqref{12.8} and formula \eqref{7.11} 
($\Lambda_a\cdot \Lambda_b=\gcd(a,b)\Lambda_{\lcm(a,b)}$)
tell 
\begin{eqnarray}\label{12.18}
\sum_{k\in\N}\nu(k)\cdot\Lambda_k &=&\prod_{j=1}^n
\bigl(\frac{1}{s_j}\Lambda_{t_j}-\Lambda_1\bigr),\\
\sum_{k:\, p^{v_p(m_a)}\nmid k}\nu(k)\cdot \Lambda_k
&=&(-1)^{n-\www n}\cdot \prod_{j=1}^{\www n} 
\bigl(\frac{1}{s_j}\Lambda_{t_j}-\Lambda_1\bigr).
\label{12.19}
\end{eqnarray}

Now we claim that the shorter weight system 
$(w_1,...,w_{\www n})$ satisfies Kouchnirenko's condition (C2)
in Theorem \ref{t12.5}. To prove this, start with a 
subset $J\subset\{1,...,\www n\}$ with $J\neq\emptyset$.
The weight system $(w_1,....,w_n)$ satisfies (C2).
Therefore a set $K\subset\{1,...,n\}$ with $|K|=|J|$
and $\forall\ k\in K$ $(\N_0^J)_{1-w_k}\neq\emptyset$ exists.
For any $k\in K$ choose a tuple
$(\alpha_1,...,\alpha_n)\in (\N_0^J)_{1-w_k}$,
so $\sum_{j\in J}w_j\alpha_j=1-w_k$.
As $J\subset\{1,...,\www n\}$, $p^{v_p(m_a)}\nmid t_j$
for any $j\in J$. Therefore $p^{v_p(m_a)}\nmid t_k$.
This shows $k\in\{1,...,\www n\}$. Therefore
$K\subset\{1,...,\www n\}$. Thus the weight system
$(w_1,...,w_{\www n})$ satisfies Kouchnirenko's condition (C2).

Theorem \ref{t12.5} gives us the existence of a 
quasihomogeneous polynomial $\www f$ with this weight
system and an isolated singularity at 0. 
Write $\www\chi:=\chi_{Mil}^{\www f}$ for the exponential
map of its characteristic polynomial $p_{Mil}^{\www f}$. 
Define the map $\www \nu:\N\to\Z$ by \eqref{12.13} with
$\chi$ replaced by $\www\chi$. Then 
\begin{eqnarray}\label{12.20}
\divv p_{Mil}^{\www f}=\sum_{m\in\N}\www\chi(m)\cdot \Psi_m
=\sum_{k\in\N}\www\nu(k)\cdot\Lambda_k.
\end{eqnarray}
By Theorem \ref{t12.4}, this is equal to
$\prod_{j=1}^{\www n}\bigl(
\frac{1}{s_j}\Lambda_{t_j}-\Lambda_1\bigr)$,
which is up to the sign $(-1)^{n-\www n}$ the right hand side
of \eqref{12.19}. Therefore
\begin{eqnarray}\label{12.21}
\supp(\www\nu)&\subset&\{k\in\N\,|\, p^{v_p(m_a)}\nmid k\},\\
\www\nu(k)&=&(-1)^{n-\www n}\cdot \nu(k)\quad 
\textup{for }k\in\supp(\www\nu)).\label{12.22}
\end{eqnarray}
This and \eqref{12.16} show
\begin{eqnarray}
\chi(m_b)-\chi(m_a)
&=&(-1)^{n-\www n}\cdot\sum_{k:\, m_b|k}\www\nu(k)\nonumber\\
&=&(-1)^{n-\www n}\cdot \www\chi(m_b)
\in (-1)^{n-\www n}\cdot \N_0.\label{12.23}
\end{eqnarray}
This implies the claim \eqref{12.12} which we had to prove.
\hfill$\Box$

\section{Looking backward and forward}\label{c13}
\setcounter{equation}{0}

\begin{remarks}\label{t13.1}
(i) By Theorem 1.3, Orlik's Conjecture \ref{t1.2}
holds for any
iterated Thom-Sebastiani sum of chain type singularities
and cycle type singularities (and any quasihomogeneous
singularity with the same weights as such a sum).
Such sums are also called {\it invertible polynomials}.
The Brieskorn-Pham singularities,
$f=f(x_1,...,x_n)=\sum_{i=1}^n x_i^{a_i}$ for some $n\in\N$
and some $a_1,...,a_n\in\Z_{\geq 2}$ are special cases,
as the $A$-type singularities $x_i^{a_i}$ are special
chain type singularities.

\medskip
(ii) For each weight system in $n=2$ variables which allows 
a quasihomogeneous singularity, a singularity of at least
one of the following 3 types exists: 
Brieskorn-Pham singularity (type I),
chain type singularity (type II), 
cycle type singularity (type III). 
This observation and these types are due to Arnold  
\cite{AGV85}. Therefore and because of
Theorem \ref{t1.3} (d), each quasihomogeneous curve 
singularity satisfies Orlik's conjecture. 
Michel and Weber claimed in the introduction of \cite{MW86}
that they have a proof of this. In view of their techniques,
this claim can be trusted. But the proof was not written.

\medskip
(iii) In the case of $n=3$ variables, Arnold distinguishes
7 types I, II, III, IV, V, VI and VII 
of weight systems which allow quasihomogeneous 
singularities (some weight systems are simultaneously of
several types) \cite{AGV85} (see also e.g. \cite{HK12}). 
5 of the 7 types arise from iterated
Thom-Sebastiani sums of chain type singularities and
cycle type singularities, the 2 types III and VI not. 
For each quasihomogeneous singularity with a weight
system of one of the 5 types, Orlik's conjecture holds
by Theorem \ref{t1.3} (d). 
For each quasihomogeneous singularity with a weight
system only of type III or VI, Orlik's conjecture is open.

\medskip
(iv) The isolated hypersurface singularities with modality
$\leq 2$ were classified by Arnold 1972, see \cite{AGV85}
for the results. Each of the families with modality
$\leq 2$ which contains
a quasihomogeneous singularity, contains especially
an iterated Thom-Sebastiani sum of chain type singularities
and cycle type singularities. Therefore for these families
Orlik's conjecture is true.
\end{remarks}

\begin{remarks}\label{t13.2}
(i) Orlik's paper \cite{Or72} contains a second conjecture,
which may not be confused with Conjecture \ref{t1.2}.
It is also often called {\it Orlik's conjecture}.
It is a consequence of Conjecture \ref{t1.2}, 
and it is weaker than Conjecture \ref{t1.2}.
In the case of a quasihomogeneous singularity 
$f\in\C[x_1,...,x_n]$ in $n\geq 3$ variables, 
it predicts the homology of the link 
$K:=f^{-1}\cap S^{2n-1}$ of the singularity. 
The Wang sequence (see below \eqref{13.7}) tells how this
homology looks like if Conjecture \ref{t1.2} is true. 
A first version was fixed by Orlik as Conjecture 3.2 in
\cite{Or72}. With some additional arguments and
Theorem \ref{t12.4}, he 
gave a second more explicit version Conjecture 3.3 
in \cite{Or72} which allows to determine the homology 
of $K$ solely in terms of the weights of $f$. 
Lemma \ref{t13.3} below recalls the first version. 
With Theorem \ref{t12.4} and Theorem \ref{t4.5} (a),
it is not so difficult to derive the second more 
explicit version.

\medskip
(ii) The links of some quasihomogeneous singularities
give interesting examples of Sasakian structures.
This was explored by Boyer, Galicki and others,
see e.g. \cite{Bo08}. There the explicit version
Conjecture 3.3 in \cite{Or72} is cited as Conjecture 19. 
In Theorem 27 in \cite{Bo08}, the links of 12 Brieskorn-Pham
singularities and 2 chain type singularities in $n=5$
variables are considered. Then the link has real dimension
7.

\medskip
(iii) In \cite{Bo08} in Proposition 20, 
Conjecture 19 is claimed to
be true for $n\in\{3,4\}$, for Brieskorn-Pham
polynomials and for chain type singularities. 

Our Theorem \ref{t1.3} (d) and Lemma \ref{t13.3} below
give this Conjecture 19 now for all Thom-Sebastiani
sums of chain type singularities and cycle type 
singularities. These contain the Brieskorn-Pham
singularities and the chain type singularities,
but not all singularities in $n=3$ or $n=4$ variables.

For $n=4$, Boyer cites work of Galicki. 
For $n=3$, Boyer cites \cite{OW71}. See (iv) for that paper.
But Arnold's two types III and VI are open in the case $n=3$. 
For Brieskorn-Pham singularities, 
\cite{Ra75} gives the result.
For chain type singularities, Boyer cites \cite{OR77}.
But that paper gives only our Theorem \ref{t10.1}.
One needs additionally our algebraic Theorem \ref{t6.2},
see section \ref{c10} above. 

\medskip
(iv) The paper \cite{OW71} and also \cite{Or72} 
miss in the classification of quasihomogeneous singularities
with $n=3$ variables 
Arnold's types III and VI. Proposition (3.1.2) and
Theorem (3.1.4) in \cite{OW71} list only the other 5 types. 
Proposition 2.2 in \cite{Or72} even claims that
each weight system which allows a quasihomogeneous
singularity allows also a Thom-Sebastiani sum 
of chain type singularities and cycle type singularities.
This is wrong for $n\geq 3$. 
Proposition 3.4 in \cite{Or72} claims that Conjecture 3.3
in \cite{Or72} (= Conjecture 19 in \cite{Bo08}) is true
for $n=3$. 
It refers to \cite{OW71}. Because there Arnold's types
III and VI are missed, for these types Conjecture 3.3
in \cite{Or72} is open. 
\end{remarks}

Lemma \ref{t13.3} shows how Conjecture 3.2 in \cite{Or72}
follows from Conjecture \ref{t1.2} (=Conjecture 3.1 in
\cite{Or72}). 

\begin{lemma}\label{t13.3}\cite{Or72}
Let $f\in\C[x_1,...,x_n]$ be a quasihomogeneous singularity
with weight system $(w_1,...,w_n)$ with
$w_i=\frac{s_i}{t_i}\in\Q\cap (0,1)$ and $s_i,t_i\in\N$,
$\gcd(s_i,t_i)=1$.
Let $l:=\chi_{Mil}^f(1)\in\N_0$ be the 
multiplicity of $\Phi_1$ as a factor of the
characteristic polynomial $p_{H_{Mil},h_{Mil}}=p_{Mil}^f$.
Then
\begin{eqnarray}\label{13.1}
l&=& \sum_{s=0}^n\sum_{\{i_1,...,i_s\}\subset\{1,...,n\}}
\frac{(-1)^s}{w_{i_1}\cdot ...\cdot w_{i_s}
\cdot\lcm(t_{i_1},...,t_{i_s})}.
\end{eqnarray}
Let $p_1,...,p_k$
be the elementary divisors of the characteristic polynomial 
$p_{Mil}^f$,
so they are products of cyclotomic polynomials with
multiplicities 1 and with
\begin{eqnarray}\label{13.2}
p_{Mil}^f=p_1\cdot ...\cdot p_k,\quad
p_k|p_{k-1}|...|p_2|p_1. 
\end{eqnarray}
Then $k\geq l$ and 
\begin{eqnarray}
p_j(1)\left\{\begin{array}{ll}
=0&\textup{for }j\in\{1,...,l\},\\
\in\N&\textup{for }j\in\{l+1,...,k\},
\end{array}\right. \label{13.3}\\
\textup{and}\quad p_k(1)|p_{k-1}(1)|...|p_{l+1}(1).
\nonumber
\end{eqnarray}

Now suppose that $f$ satisfies Conjecture \ref{t1.2}. 
Then the homology of the link $K:=f^{-1}(1)\cap S^{2n-1}$
is given by 
\begin{eqnarray}\label{13.4}
H_{n-1}(K,\Z)\cong\Z^l,\quad
H_{n-2}(K,\Z)\cong\Z^l\oplus\bigoplus_{j=l+1}^k
\frac{\Z}{p_j(1)\Z}.
\end{eqnarray}
\end{lemma}

{\bf Proof:} 
Because of \eqref{7.11}, 
formula \eqref{12.8} in Theorem \ref{t12.4} can be
rewritten as the following sum with rational coefficients,
\begin{eqnarray}\nonumber
\divv p_{Mil}^f &=& 
\sum_{s=0}^n\sum_{\{i_1,...,i_s\}\subset\{1,...,n\}}\\
&&\frac{(-1)^s}{w_{i_1}\cdot ...\cdot w_{i_s}\cdot
\lcm(t_{i_1},...,t_{i_s})} \cdot 
\Lambda_{\lcm(t_{i_1},...,t_{i_s})}.\label{13.5}
\end{eqnarray}
Each $\Lambda_m$ contains $\Phi_1$ with multiplicity 1.
This shows \eqref{13.1}. 
The uniqueness of the elementary divisors was discussed
in Remark \ref{t2.5} (iv). Their definition gives
immediately \eqref{13.3}. 

Now suppose that $f$ satisfies Conjecture \ref{t1.2}. 
That means
\begin{eqnarray}\label{13.6}
(H_{Mil},h_{Mil})\cong\bigoplus_{j=1}^k 
(H^{[p_j]},h_{[p_j]})
=\bigoplus_{j=1}^k(\frac{\Z[t]}{p_j\Z[t]},
\textup{mult.\,by\,}t).
\end{eqnarray}

The Wang sequence connects in the case $n\geq 3$ the 
homology of the link $K$ with the pair $(H_{Mil},h_{Mil})$ 
of Milnor lattice and monodromy, see \cite{Or72}.
It gives the following short exact sequence:
\begin{eqnarray}\label{13.7}
0\to H_{n-1}(K,\Z)\to H_{Mil}
\stackrel{h_{Mil}-\id}{\longrightarrow} H_{Mil}
\to H_{n-2}(K,\Z)\to 0.
\end{eqnarray}
We find
\begin{eqnarray}\label{13.8}
H_{n-1}(K,\Z)&\cong&\ker\Bigr(h_{Mil}-\id:H_{Mil}\to H_{Mil}
\Bigr) \cong\Z^l,\\
H_{n-2}(K,\Z)&\cong& \frac{H_{Mil}}{(h_{Mil}-\id)(H_{Mil})}
\nonumber \\
&\cong& \bigoplus_{j=1}^k \frac{\Z}{p_j(1)\Z}
\cong\Z^l\oplus\bigoplus_{j=l+1}^k \frac{\Z}{p_j(1)\Z}.
\label{13.9}
\end{eqnarray}
Here $\Z/0\Z\cong \Z$, of course. \hfill$\Box$

\begin{remarks}\label{t13.4}
(i) The first author's interest in Orlik's Conjecture
\ref{t1.2} comes from his work on a Torelli conjecture
for isolated hypersurface singularities (\cite{He92}
and many later papers). 
Proofs of this conjecture for many classes of singularities
consist of two steps, a transcendent step, the calculation of
a period map, and an algebraic step, the determination
of the action of the automorphism group of the pair 
$(H_{Mil},\textup{Seifert form})$
on the image of this period map. 
The Seifert form is a bilinear unimodal form
$H_{Mil}\times H_{Mil}\to\Z$ which determines the
monodromy $h_{Mil}$. In many cases, it is useful to
determine first the automorphism group of the pair
$(H_{Mil},h_{Mil})$, and this becomes easier
if this pair satisfies Orlik's conjecture. 
By Remark \ref{t12.3} (which builds on Theorem \ref{t12.1}
and \cite{He20}), the automorphism group of an
Orlik block $(H^{[p_j]},h_{[p_j]})$, 
where $p_j$ is as in Lemma 
\ref{t13.3}, is simply $\{\pm h_{[p_j]}^j\,|\, j\in\Z\}$. 

\medskip
(ii) It would be useful to generalize Orlik's Conjecture
\ref{t1.2} to a conjecture on the pair 
$(H_{Mil},\textup{Seifert form})$ for a quasihomogeneous
singularity. But it is not at all clear how this 
generalization could look like. 
This is an interesting open question. 
\end{remarks}


\begin{thebibliography}{99}
\bibitem[Ap70]{Ap70} T. Apostol: \quad
      Resultants of cyclotomic polynomials.
      Proc. A.M.S. {\bf 24} (1970), 457--462.
\bibitem[AGV85]{AGV85} Arnold, V.I., S.M. Gusein-Zade, 
      A.N. Varchenko: \quad 
      Singularities of differentiable maps, volume I. 
      Birkh\"auser, Boston 1985.
\bibitem[Bo08]{Bo08} C.P. Boyer: \quad 
      Sasakian geometry: the recent work of 
      Krzysztof Galicki.
      Note di Matematica {\bf 1}, suppl. n {\bf 1} (2008),
      63--105.
\bibitem[Co82]{Co82} B.G. Cooper: \quad
      On the monodromy at isolated singularities of 
      weighted homogeneous polynomials.
      Transactions of the A.M.S. {\bf 269} (1982), 149--166.
\bibitem[He92]{He92} C. Hertling: 
      Analytische Invarianten bei den unimodularen und 
      bimodularen Hyperfl\"achensingularit\"aten. 
      Dissertation. 188 pages. Bonner Math. Schriften 250 , 
      Bonn 1992. 
\bibitem[HK12]{HK12} C. Hertling, R. Kurbel:\quad 
     On the classification of quasihomogeneous singularities.
     J. Singul. {\bf 4} (2012), 131--153. 
\bibitem[He20]{He20} C. Hertling: \quad 
     Automorphisms with eigenvalues in $S^1$ of a $\Z$-lattice 
     with cyclic finite monodromy. 
     Acta Arithmetica {\bf 192.1} (2020), 1--30.
\bibitem[HZ19]{HZ19} C. Hertling, Ph. Zilke: \quad 
     Seven combinatorial problems around isolated 
     quasihomogeneous singularities.
     Journal of Algebraic Combinatorics {\bf 50} (2019),
     447--482.
\bibitem[HM20-1]{HM20-1} C. Hertling, M. Mase: \quad 
     The integral monodromy of cycle type singularities.
     Preprint, arXiv:2009.07533, September 2020, 38 pages.
\bibitem[HM20-2]{HM20-2} C. Hertling, M. Mase: \quad 
     Combinatorial results around isolated quasihomogeneous
     singularities. 
     Manuscript in preparation (autumn 2020).
\bibitem[Ko76]{Ko76} A.G. Kouchnirenko:  \quad
     Poly\`edres de Newton et nombres de Milnor.
     Invent. Math. {\bf 32}, 1--31 (1976)
\bibitem[MW86]{MW86} F. Michel, C. Weber: \quad 
     Sur le r\^ole de la monodromie 
     enti\`ere dans la topologie des singularit\'es. 
     Ann. Inst. Fourier Grenoble {\bf 36} (1986), 183--218. 
\bibitem[Mi68]{Mi68} J. Milnor: \quad
     Singular points of complex hypersurfaces.
     Annals of Mathematics Studies {\bf 61}, 
     Princeton University Press, 1968.
\bibitem[MO70]{MO70} J. Milnor, P. Orlik: \quad 
     Isolated singularities defined by weighted 
     homogeneous polynomials.
     Topology {\bf 9} (1970), 385--393.
\bibitem[Or72]{Or72} P. Orlik: \quad 
    On the homology of weighted homogeneous manifolds. 
    In: Lecture Notes in Math. 298, Springer, Berlin, 
    1972, pp. 260--269.
\bibitem[OR77]{OR77} P. Orlik, R.C. Randell: \quad
     The monodromy of weighted homogeneous singularities.
     Invent. Math. {\bf 39} (1977), 199--211.
\bibitem[OW71]{OW71} P. Orlik, Ph. Wagreich: \quad
     Isolated singularities of algebraic surfaces with 
     $\C^*$ action.
     Ann. of Math. {\bf 93.2} (1971), 205--228.
\bibitem[Ra75]{Ra75} R.C. Randell: \quad
     The homology of generalized Brieskorn manifolds.
     Topology {\bf 14} (1975), 347--355.
\bibitem[ST71]{ST71} M. Sebastiani, R. Thom: \quad 
     Un resultat sur la monodromie.
     Invent. Math. {\bf 13} (1971), 90--96.
\bibitem[vW71]{vW71} B.L. van der Waerden: \quad
     Algebra I, Springer, Berlin, Heidelberg, New York, 
     8. Auflage, 1971.
\bibitem[Wa82]{Wa82} L.C. Washington: \quad 
     Introduction to cyclotomic fields.
     Springer, New York, 1982.
\end{thebibliography}
\end{document}